\documentclass[12pt,a4paper,twoside]{article}

\usepackage[active]{srcltx}
\usepackage[utf8]{inputenc}
\usepackage{hyperref}
\usepackage[english]{babel}
\usepackage{authblk}

\usepackage{xcolor}

\usepackage[a4paper,hcentering,vcentering]{geometry}
\usepackage{url,epsfig,amssymb,amsthm,amsmath}

\numberwithin{equation}{section}

\newtheorem{theorems}{Theorem}
\numberwithin{theorems}{section}

\numberwithin{corollarys}{section}

\newtheorem{definition}{Definition}
\numberwithin{definition}{section}

\newtheorem{lemma}{Lemma}
\numberwithin{lemma}{section}
\newtheorem{proposition}{Proposition}
\numberwithin{proposition}{section}
\newtheorem{theoremx}{Theorem}

\theoremstyle{definition}
\newtheorem{remark}{Remark}
\numberwithin{remark}{section}

\newcommand{\N}{\ensuremath{\mathbb{N}}}

\newcommand{\R}{\ensuremath{\mathbb{R}}}
\newcommand{\T}{\ensuremath{\mathbb{T}}}
\newcommand{\Z}{\ensuremath{\mathbb{Z}}}

\definecolor{JungleGreen}{cmyk}{0.99,0,0.52,0}
\definecolor{RawSienna}{cmyk}{0,0.72,1,0.45}
\definecolor{bulgarianrose}{rgb}{0.28, 0.02, 0.03}
\definecolor{Magenta}{cmyk}{0,1,0,0}
\definecolor{airforceblue}{rgb}{0.36, 0.54, 0.66}
\definecolor{darkpastelgreen}{rgb}{0.01, 0.75, 0.24} 
\definecolor{brightgreen}{rgb}{0.4, 1.0, 0.0} 
\definecolor{burntorange}{rgb}{0.8, 0.33, 0.0} 
\definecolor{amber}{rgb}{1.0, 0.75, 0.0} 
\definecolor{amber(sae/ece)}{rgb}{1.0, 0.49, 0.0} 

\providecommand{\keywords}[1]
{
  \small	
  \textbf{\textit{Keywords---}} #1
  \normalsize
}

\begin{document}

\title{\LARGE{\textbf{Weakly asymptotically quasiperiodic solutions for time-dependent Hamiltonians with a view to celestial mechanics}}}%
\author{Donato Scarcella \thanks{E-mail address: donato.scarcella@upc.edu}}
\affil{Departament de Matemàtiques, Universitat Politècnica de Catalunya, Diagonal 647, 08028 Barcelona, Spain.}
\date{}%
\maketitle


\begin{abstract}
We consider the planar three-body problem perturbed by a celestial body modeled as a time-dependent perturbation that decays in time.  We assume that the motion of the celestial body is given and is unbounded with a non-zero asymptotic velocity. We prove the existence of orbits converging in time to some motions that are ``close'' to the quasiperiodic solutions associated with the Hamiltonian of the planar three-body problem.  

The proof relies on an abstract theorem that contains a substantial portion of the mathematical complexities presented in this work. This theorem is flexible and can be applied to many other physical phenomena. It considers Hamiltonian vector fields that are the sum of two components. The first possesses quasiperiodic solutions, and the second decays polynomially fast as time tends to infinity. We prove the existence of orbits converging in time to some motions that are ``close'' to the quasiperiodic solutions associated with the unperturbed system. 
 It generalizes a previous work where a stronger polynomial decay in time was considered, and solutions converging in time to the quasiperiodic orbits associated with the unperturbed system were proved.  
 In the abstract theorem contained in the present paper, the too-weak decay in time of the perturbation strongly modifies the dynamic at infinity.  This serious difficulty requires a deep modification of the proof. This new strategy relies on the application of a Nash-Moser implicit function theorem (the previous result was proved with the fixed point theorem) and the introduction of weak solutions (in this case, the orbits do not converge to the quasiperiodic solutions associated with the unperturbed system).
\end{abstract}
\keywords{Dynamical systems, Hamiltonian systems, Time-dependence,  Celestial mechanics, $3+1$ body problem.}

\tableofcontents

\section{Introduction}

Many physical phenomena can be described by dynamical systems subjected to a perturbation that decays in time.  For instance, a molecule interacting with another molecule or with a laser pulse~\cite{kawai2007transition, BdlL11}, and in reference~\cite{thieme1995asymptotically}, one can observe the development of some epidemiological models.  In the present paper, we consider the planar three-body problem perturbed by a given celestial body.  The motion of this celestial body is a given smooth function $c(t)$, which depends on time,  and only the planetary system is influenced by $c(t)$.  More specifically, we assume that 
\begin{equation}
\label{chypIntro}
\lim_{t \to +\infty}|c(t)| = + \infty,  \quad \lim_{t \to +\infty}{d \over dt}|c(t)| = v_c>0.
\end{equation}
In other words, $c(t)$ diverges in time with a nonzero asymptotic velocity.  This model can describe many celestial phenomena. For example, the case of the planar three-body problem perturbed by a given comet coming from and going back to infinity asymptotically along a hyperbolic Keplerian orbit.  In this case,  if the comet is very far from the planetary system,  one can assume that, in a first approximation, $c(t)$ is not influenced by the gravitational attraction exerted by the planetary system. Another interesting example is the planar three-body problem perturbed by a big star moving in a uniform rectilinear motion. The star, because of its large mass, in a first approximation, is not affected by the motion of the other three planets. 

On a suitable phase space, the Hamiltonian that describes the planar three-body problem plus the celestial body is a time-dependent Hamiltonian of the form
\begin{equation}
\label{HP3BP+CB}
H^t = H_0 + H^t_c,
\end{equation}
where $H_0$ is the Hamiltonian of the planar three-body problem, and $H^t_c$ is responsible for the interaction of the three bodies with the celestial body $c(t)$. We will call $H^t$ the Hamiltonian of the planar three-body problem plus the celestial body (P3BP+CB).

We recall that, in 1963, Arnold proved the existence of quasiperiodic solutions for the Hamiltonian of the planar three-body problem~\cite{Arn63b}. In this work, however, we follow the setting of Féjoz~\cite{Fe02}, which provides more general solutions. Indeed, in a rotating frame of reference, the author proves the existence of quasiperiodic orbits with three frequencies for the Hamiltonian of the planar three-body problem. Before the symplectic reduction by the symmetry of rotations, these quasiperiodic motions have one additional frequency, namely the angular speed of the simultaneous rotation of the three ellipses. Furthermore, before the symplectic reduction by the symmetry of translations, each of these invariant tori translates into a $1$-parameter family of invariant tori parametrized by the center of mass of the planetary system.  That is an invariant cylinder supporting quasiperiodic solutions.  We refer to Section \ref{QPMP3BP} for a more in-depth description of the result obtained by Féjoz. 

In the present paper, in a neighborhood of the above-mentioned invariant cylinder,  we fix the motion of the celestial body $c(t)$, satisfying~\eqref{chypIntro}, and we study how it perturbes the three-body problem.  The Hamiltonian~\eqref{HP3BP+CB}, which describes this dynamical system, depends explicitly on time.   We will see that,  in an appropriate neighborhood of an invariant cylinder supporting quasiperiodic solutions associated with $H_0$, the Hamiltonian vector field associated with the perturbation $H^t_c$ satisfies the following decay in time
\begin{equation}
\label{IntroDecayHct}
|X_{H^t_c}| \sim {1 \over t^2}
\end{equation}
on a suitable norm.  Roughly speaking, for all $i =0, 1, 2$,  let $x_i$ be the positions of each body.  On a suitable subset of the phase space $H_c^t(x)\sim\sum_{i=0}^2{1 \over |x_i -c(t)|}$.  One can see that~\eqref{chypIntro} implies $|c(t)| \sim \mathrm{cst} + v_ct$ for $t$ sufficiently large.  Hence, initially, we have that $|H^t_c| \sim {1 \over t}$ on a suitable norm.  After an expansion using Legendre polynomials, when we consider the Hamiltonian vector field associated with $H_c^t$, we will see that we obtain an extra power in the decay in time, that is~\eqref{IntroDecayHct}.  We refer to Section \ref{P3C} for more details. 

For $|c(1)|$ and $v_c$ large enough, we prove the existence of an open set of initial conditions giving rise to orbits for the Hamiltonian of the P3BP+CB that converge as time tends to infinity to motions that are ``close'' (in some way we specify later) to some quasiperiodic solutions of frequency vector $\omega$ associated with $H_0$ (see Theorem \ref{Thmcomet}).  We call these motions weakly asymptotically quasiperiodic solutions (see Definition \ref{asymsolI}).  The reason is that the motions do not converge to the quasiperiodic solutions associated with the unperturbed system but rather to orbits that are ``close'' to quasiperiodic solutions.
Indeed, for each of these solutions, the center of mass of the planetary system is attracted by the celestial body $c(t)$ with a zero asymptotic velocity. Whereas,  in a frame of reference attached to the center of mass of the planetary system, the motions of the planets converge in time to some orbits whose dynamic is conjugated to a time-dependent perturbation of the constant vector field $\omega$.

In the sequel, $C^\sigma$ indicates the class of Hölder functions and $|\cdot|_{C^\sigma}$ the Hölder norm (we refer to Appendix \ref{A} for a very brief introduction). The above-mentioned result relies on an abstract theorem that is interesting by itself and encapsulates a substantial portion of the mathematical complexities presented in this work.  It is about the existence of $C^\sigma$-weakly asymptotic cylinders (see Definition \ref{weakasymcyl} below) for time-dependent perturbations decaying polynomially fast in time of Hamiltonian systems having an invariant cylinder supporting quasiperiodic solutions of frequency vector $\omega$ (we refer to Theorem \ref{MD}).  Roughly speaking,  a $ C^\sigma$-weakly asymptotic cylinder is a family of embedded cylinders converging as time tends to infinity to the invariant cylinder associated with the unperturbed system. The motions on this family of embeddings converge in time to some orbits conjugated to the flow of a time-dependent perturbation of the constant vector field on the cylinder $\T^n \times \R^m$ given by $(\omega, 0)$.  We want to point out that the hypotheses on the dynamics associated with the unperturbed system are slightly more general in Definition \ref{weakasymcyl}, but this is the more natural situation.  The same remark works for the hypotheses of Theorem \ref{MD}. For the sake of clarity, in this section, we prefer to consider this more natural and easier setting. 

The definition of $ C^\sigma$-weakly asymptotic cylinder we present generalizes the one of $ C^\sigma$-asymptotic KAM torus (see Definition \ref{AsymKAMtorus}) introduced by Canadell and de la Llave~\cite{CdlL15}.  In words,  let $X^t$ be a time-dependent vector field converging as time tends to infinity to a vector field $X_0$ having an invariant torus $\varphi_0$ supporting quasiperiodic solutions.  A $ C^\sigma$-asymptotic KAM torus is a family of embedded tori converging in time to the invariant torus associated with the unperturbed system $X_0$.  The motions on this family of embeddings converge in time to the quasiperiodic solutions associated with $X_0$ on the invariant torus $\varphi_0$.
We stress that in the present paper, we are interested in establishing the existence of families of embedded cylinders rather than embedded tori. This decision is motivated by the fact that, as mentioned before,  the Hamiltonian of the planar three-body problem has an invariant cylinder supporting quasiperiodic solutions in the phase space before the symplectic reductions of the symmetries of translations and rotations.  The definition of $ C^\sigma$-weakly asymptotic cylinder generalizes the one of $C^\sigma$-asymptotic KAM torus because, in this case,  we lose some information on the dynamic at infinity (the motions do not converge to the dynamic associated to the unperturbed system).  We refer to Section \ref{MR} for a more detailed comparison. 

As mentioned before,  in the present paper, we prove the existence of $C^\sigma$-weakly asymptotic cylinders for time-dependent Hamiltonian vector fields converging polynomially fast in time to Hamiltonian vector field having an invariant cylinder supporting quasiperiodic solutions.  More specifically, we consider time-dependent Hamiltonian vector fields of the form $X^t = X_0 + F^t$ where $X_0$ has an invariant cylinder $\varphi_0$ supporting quasiperiodic solutions.  Moreover,  we assume that $F^t(q,0)$ and $\partial_p F^t(q,0)$ decay as ${1 \over t^2}$ and ${1 \over t}$, respectively.  Then, we prove the existence of a $C^\sigma$-weakly asymptotic cylinder (see Theorem \ref{MD}).  

The present paper is the second part of a previous work~\cite{Sca22} where a stronger time decay is considered.  In reference~\cite{Sca22}, time-dependent Hamiltonian vector fields of the form $X^t = X_0 + F^t$ are analyzed.  Here, $X_0$ has an invariant torus $\varphi_0$ supporting quasiperiodic solutions of frequency vector $\omega$. Moreover, letting $l >2$, we assume that $F^t(q,0)$ and $\partial_p F^t(q,0)$ decay as ${1 \over t^l}$ and ${1 \over t^{l-1}}$, respectively. We prove the existence of a $C^\sigma$-asymptotic KAM torus $\varphi^t$ satisfying $|\varphi^t - \varphi_0|_{C^\sigma} \le C {1 \over t^{l-2}}$ for a suitable constant $C$. This previous paper generalizes the work of Canadell and de la Llave~\cite{CdlL15} where exponential decay in time was considered, and the one of Fortunati and Wiggins~\cite{FW14} where arithmetic, non-degeneracy conditions, and exponential decay in time are assumed. We note that neither in~\cite{Sca22} nor in~\cite{CdlL15} arithmetic conditions on the frequency vector $\omega$ or non-degeneracy assumptions on the unperturbed Hamiltonian are considered.

Unlike the result contained in the present paper, in reference~\cite{Sca22}, the stronger decay in time allows for a stronger conclusion to be proved. The existence of a $C^\sigma$-asymptotic KAM torus establishes the presence of orbits converging as time tends to infinity to the quasiperiodic solutions associated with the unperturbed system.  The cost of this stronger result is a loss of two powers in the decay in time of the found $C^\sigma$-asymptotic KAM torus $\varphi^t$ (we recall that $\varphi^t$ converges in time to the invariant torus associated to the unperturbed system as ${1 \over t^{l-2}}$).
The proof of the result in~\cite{Sca22} relies on the fixed point theorem.  It is based on a different abstract formulation of the dynamical problem with respect to the previous works~\cite{FW14, CdlL15}.
An important point of the proof concerns the solution and estimation of the associated linearized problem. It is solved by integration thanks to a suitable change of coordinates that rectifies the dynamic on the torus.  

Now, we stress that the Hamiltonian vector field associated with the time-dependent perturbation $H_c^t$ of the Hamiltonian of the P3BP+CB (see~\eqref{HP3BP+CB}) decays as ${1 \over t^2}$. Hence, the decay in time is too weak, and it is not sufficient to apply the result in~\cite{Sca22}.
This serious difficulty requires the formulation of another abstract theorem (Theorem \ref{MD}).  In fact, the dynamic at infinity is strongly modified by the non-autonomous perturbation.  For this reason, the proof of Theorem \ref{MD} requires significant modification compared to the one in~\cite{Sca22}.  Moreover, we prove a weaker result because we establish the existence of a $C^\sigma$-weakly asymptotic cylinder rather than a $C^\sigma$-asymptotic KAM torus. Hence, we lose some information about the dynamic at infinity of the found solutions. 
The slower decay in time is at the origin of a loss of derivatives, making fixed point arguments useless.  We need to introduce a new formulation of the dynamical problem and employ a different strategy using tame estimates of the composition and product of Hölder functions in an appropriate scale of Banach spaces (we refer to Remark \ref{CompGMD} for a more detailed comparison between the proof of the result in~\cite{Sca22} and the one of Theorem \ref{MD} contained in this paper).  
The proof relies on the Nash-Moser theorem, and it is based on a refined analysis of the Banach spaces involved in the solutions of non-linear invariant equations.
 In this case, the associated linearized problem is considerably more complicated than the one in~\cite{Sca22}. It relies on the solution of the homological equation analyzed in Section \ref{HENCBD} (we refer to Remark \ref{CompHE} for a comparison between the associated linearized problem studied in the present paper and the one in~\cite{Sca22}). 
 
In another paper~\cite{Sca22b},  using a slightly different version of the result contained in~\cite{Sca22},  we proved the existence of motions converging to quasiperiodic solutions in the future and in the past for time-dependent perturbations of integrable and near-integrable Hamiltonians. On the other hand, in reference~\cite{Sca22c}, the case of time-dependent perturbations of Hamiltonian having an invariant torus supporting arbitrary dynamic is considered. Assuming exponential decay in time for the perturbation, we prove the existence of orbit converging in time to the motions associated with the unperturbed system. Also in this case, the proof is achieved using the fixed point theorem. 

We point out that there are many examples in the literature of planetary systems having planets subjected to unbounded motions.  For example,  we can consider the parabolic-elliptic or hyperbolic-elliptic type defined by Chazy~\cite{Chazy22}.  In this case, one of the bodies goes to infinity with a zero or non-zero asymptotic velocity, respectively, while the others perform bounded motions.  A recent work~\cite{BFMn+1} investigates the existence of solutions for the planar $(n+1)$-body problem defined for all forward time and tending to a parabolic motion.  More specifically, it proves the existence of orbits when one of the bodies goes to infinity with a zero asymptotic velocity while the rest perform a bounded motion. The existence of parabolic motions for the $n$-body problem is also studied in~\cite{MV} using variational methods.  Meanwhile,  in~\cite{BDT} the $n$-center problem is considered, and parabolic trajectories having asymptotic directions are proved.

Finally, this paper contributes to an extensive series of works (\cite{CdlL15, FW14, Sca22, Sca22b, Sca22c}) dedicated to the development of a non-autonomous KAM theory.  The classical KAM theory shows the persistence of quasiperiodic solutions in nearly integrable Hamiltonian systems, where KAM is the acronym of the three mathematicians who are at the origin of this theory Kolmogorov~\cite{Kol54},  Arnold~\cite{Arn63a, Arn63b}, and Moser~\cite{moser1962invariant}.  These first works laid the foundations for the development of a large theory encompassing numerous variations and extensions.  We refer to Bost~\cite{bost1986tores}, Pöschel ~\cite{Posc01}, Chierchia~\cite{Chi03} and Féjoz~\cite{Fe16} for very interesting surveys.  The interest in the KAM theory is also attributed to its application in stability problems in celestial mechanics. Concerning some stability results in the $n$-body problem,  one can see~\cite{Arn63b},~\cite{PC11} and~\cite{Fe04}.

This theory is complex because of the small denominators that appear in the perturbation series, making the convergence argument extremely delicate.  Under an arithmetic assumption on the frequency vector and a non-degeneracy assumption on the unperturbed Hamiltonian, Kolmogorov overcomes this problem thanks to an iterative scheme based on a Newton algorithm characterized by a super-linear convergent method.  We refer to~\cite{BGGS84} for a detailed explanation of Kolmogorov's proof.  Many other interesting approaches have been developed to prove the KAM theorem.  Remarkable is the introduction of a Nash-Moser implicit function theorem in a scale of Banach spaces that replaces the iterative scheme introduced by Kolmogorov.  Some works based on this idea are the ones of Zehnder~\cite{Zeh76, Zeh75}, Herman~\cite{bost1986tores}, Berti-Bolle~\cite{BB15}, and Féjoz~\cite{Fe04}.  For the non-conservative normally hyperbolic case, one can see Massetti~\cite{M18, M19}.

The nonautonomous KAM theory differs from the classical KAM theory because of the lack of arithmetic conditions on the frequency vector and non degeneracy assumptions on the unperturbed Hamiltonian. This is due to the absence of small denominators when considering the homological equations appearing in time-dependence problems.  The interest of the present paper relies on the existence of solutions that do not have their counterpart in the classical KAM theory. We found invariant families of embedded cylinders without controlling the dynamic at infinity.  It is an unexplored direction in the classical KAM theory, where the existence of invariant tori is strictly linked with a specific dynamic. 

The paper is organized as follows.  In Section \ref{MR}, we present the main results of the paper.  Section \ref{FS} contains properties about $C^\sigma$-weakly asymptotic cylinders,  along with a series of notations and definitions that we will use in the rest of this paper.  Section \ref{HENCBD} is devoted to the solution of the homological equation.  In Section \ref{Zehnder}, we recall the version of the Nash-Moser theorem that we will use in the proof of Theorem \ref{MD} contained in Section \ref{ProofMD}.  
In Section \ref{QPMP3BP}, we recall the result of Féjoz~\cite{Fe02} about the existence of quasiperiodic solutions of the Hamiltonian of the planar three-body problem.  Section \ref{ProofTC} contains the proof of Theorem \ref{Thmcomet} concerning the existence of weakly asymptotically quasiperiodic solutions for the Hamiltonian of the P3BP+CB.

\section{Statement of the main results}\label{MR}

This section is devoted to the introduction of the main results of the present paper. For this purpose,  we need to introduce some notations and definitions.  Let $B\subset \R^{n+m}$ be a ball centered at the origin,  and $J=[1, +\infty) \subset \R$.  Given $f:\T^n \times \R^m \times B \times J \to \R$, for a fixed $t \in J$ we define the following function
\begin{equation}
\label{ft}
f^t:\T^n \times \R^m \times B \to \R, \quad f^t(q,p)=f(q,p,t)
\end{equation}
where $q \in \T^n \times \R^m$ and $p \in B$. 
We will use this notation for the rest of this work,  including for vector-valued functions and matrices.

Now,  let us introduce the definition of $C^\sigma$-weakly asymptotic cylinder.  We recall that $C^\sigma$ stands for the class of Hölder functions and $|\cdot|_{C^\sigma}$ the Hölder norm (see Appendix \ref{A}). Furthermore,  as we indicated, we denote the time dependence with an apex $t$.
Given $\sigma \ge 0$, a positive integer $k \ge 0$ and  $\bar \omega = (\omega, 0) \in \R^{n+m}$ with $\omega \in \R^n$, we consider time-dependent vector fields $X^t$, $X^t_0$ of class $C^{\sigma+k}$ on $\T^n \times \R^m \times B$, for all fixed $t \in J$, an embedding $\varphi_0:\T^n \times \R^m \to \T^n \times \R^m \times B$ of class $C^{\sigma}$ and a time-dependent vector field $\gamma^t$ of class $C^{\sigma}$ on $\T^n \times \R^m$, for all fixed $t \in J$, such that 
\begin{align}
\label{WACD1}
&  \lim_{t \to +\infty}  |X^t - X^t_0|_{C^{\sigma +k}} = 0,\\
\label{WACD2}
& X_0^t \circ \varphi_0(q)= \partial_q \varphi_0(q)(\bar \omega + \gamma^t(q)) \hspace{2mm} \mbox{for all $(q, t) \in \T^n \times \R^m \times J$,}\\
\label{WACD3}
&  \lim_{t \to +\infty} |\gamma^t|_{C^\sigma} = 0.
\end{align}
In words, $X^t- X_0^t$ converges to zero when $t \to +\infty$. Furthermore, the vector field $X^t_0$ has an invariant cylinder $\varphi_0$, and the restriction of $X^t_0$ is conjugated to the non-autonomous vector field $\bar \omega + \gamma^t$, which is a time-dependent vector field converging to $\bar \omega$ when  $t \to +\infty$.  The most natural case is when $X_0^t$ does not depend on time and $\gamma^t \equiv 0$, but the result that we will prove works in this general context.

\begin{definition}[$C^\sigma$-weakly asymptotic cylinder]
\label{weakasymcyl}
We assume that $(X^t, X^t_0, \varphi_0)$ satisfy~\eqref{WACD1},~\eqref{WACD2} and~\eqref{WACD3}. A family of $C^\sigma$ embeddings $\varphi^t: \T^n \times \R^m \to \T^n \times \R^m \times B$  is a $C^\sigma$-weakly asymptotic cylinder associated  to $(X^t, X^t_0, \varphi_0)$ if there exists a time-dependent vector field $\Gamma$ such that, for any $t \in J$, $\Gamma^t$ is of class $C^{\sigma}$ on $\T^n \times \R^m$, and $\varphi^t$ and $\Gamma$ satisfy
\begin{align}
\label{hyp2WAC}
&  \lim_{t \to +\infty}  |\varphi^t - \varphi_0|_{C^\sigma} = 0,\\
\label{hyp1WAC}
&   X^t\circ \varphi^t(q) =  \partial_q \varphi^t(q)(\bar \omega + \Gamma^t(q))  + \partial_t \varphi^t(q) \hspace{2mm} \mbox{for all $(q, t) \in \T^n \times \R^m \times J$},\\
\label{hyp3WAC}
&  \lim_{t \to +\infty} |\Gamma^t|_{C^\sigma} = 0,
\end{align}
Moreover,  we say that $\varphi^t$ is Lagrangian if $\varphi^t(\T^n \times \R^m)$ is  Lagrangian for all $t\in J$.  

We will refer to $\Gamma^t$ as the disturbing term associated with the $C^\sigma$-weakly asymptotic cylinder.
\end{definition}
In words, a $C^\sigma$-weakly asymptotic cylinder is a family of embeddings $\varphi^t$ converging in time to the invariant cylinder $\varphi_0$ associated with $X^t_0$. Moreover, the dynamics on this family of embeddings do not converge to the motions associated with $X^t_0$ on $\varphi_0$ but to orbits generated by the time-dependent vector field $\bar \omega + \Gamma^t$, hence the term weakly.  
We now introduce the definition of weakly asymptotically quasiperiodic solutions.
\begin{definition}[Weakly asymptotically quasiperiodic solution]
\label{asymsolI}
We assume that $(X^t, X^t_0, \varphi_0)$ satisfy~\eqref{WACD1},~\eqref{WACD2} and~\eqref{WACD3}. An integral curve $g: J \to \T^n \times \R^m \times B$ of $X^t$ is a weakly asymptotically quasiperiodic solution associated to $(X^t, X^t_0, \varphi_0)$ if there exist  a time-dependent vector field $\Gamma : \T^n \times \R^m \times J \to \R^{n+m}$ and $q \in \T^n \times \R^m$ such that
\begin{equation*}
\lim_{t \to +\infty} |g(t) - \varphi_0 \circ \psi^t_{t_0, \bar \omega + \Gamma}(q)|=0,
\end{equation*}
where $\psi^t_{t_0, \bar \omega + \Gamma}$ is the flow with initial time $t_0$ of the vector field $\bar \omega + \Gamma$. 
\end{definition}
It is straightforward to verify that a $C^\sigma$-weakly asymptotic cylinder implies the existence of weakly asymptotically quasiperiodic solutions (see Proposition \ref{propasymsolI} and  Section \ref{FS} for more details).  
We observe that Definition \ref{weakasymcyl} generalizes the concept of $C^\sigma$-asymptotic KAM torus introduced by Canadell and de la Llave (one can see~\cite{CdlL15,Sca22}).
\begin{definition}[$C^\sigma$-asymptotic KAM torus]
\label{AsymKAMtorus}
We assume that $(X^t, X^t_0, \varphi_0)$ satisfy~\eqref{WACD1},~\eqref{WACD2} and~\eqref{WACD3} and that $\varphi^t$ is a $C^\sigma$-weakly asymptotic cylinder associated to $(X^t, X^t_0, \varphi_0)$ with disturbing term $\Gamma^t$.  When $m=0$, $\gamma^t \equiv 0$ and $\Gamma^t \equiv 0$, we call $\varphi^t$ a $C^\sigma$-asymptotic KAM torus associated to $(X^t, X^t_0, \varphi_0)$.
\end{definition}
Roughly speaking, a $C^\sigma$-asymptotic KAM torus is a family of embedded tori $\varphi^t$ converging in time to the invariant torus $\varphi_0$. Moreover, the dynamic on this family of embeddings converges as time tends to infinity to the quasiperiodic solutions associated with $X_0^t$ on the invariant torus $\varphi_0$.

The subsequent content is subdivided into two parts.  In the first part (Section \ref{AbsThm}), we present the abstract theorem concerning the existence of $C^\sigma$-weakly asymptotic cylinders for some time-dependent Hamiltonians (Theorem \ref{MD}). The second part (Section \ref{P3BP+CThm}) contains the result about the existence of weakly asymptotically quasiperiodic solutions for the Hamiltonian of the P3BP+CB (Theorem~\ref{Thmcomet}).

\subsection{The abstract theorem}\label{AbsThm}

We need to introduce some notations.  Let $\sigma \ge 0$ and $l \ge 0$ be positive real parameters. In the sequel, we denote by $[\sigma]$ the integer part of $\sigma$. For $i \in\N^{2(n+m)}$,  let 
\begin{equation*}
|i| = i_1+...+i_{2(n+m)} \quad \mbox{and} \quad \partial_{(q,p)}^i=\partial_{q_1}^{i_1} \ldots \partial_{q_{n+m}}^{i_{n+m}}\partial_{p_1}^{i_{n+m+1}}\ldots \partial_{p_{n+m}}^{i_{2(n+m)}}
\end{equation*}
be the partial derivatives of order $|i|$ with respect to the variables $(q,p)$.  As a convention, we define $\partial_{(q,p)}^0 f = f$. 

Now, to quantify the regularity and decay in time of smooth functions, we introduce the following Banach space. 
\begin{eqnarray}
\mathcal{S}_{\sigma, l} &=& \Big\{ f:\T^n \times \R^m \times B \times J \to \R \hspace{1mm} | \hspace{1mm} f^t \in C^\sigma( \T^n \times \R^m \times B) \hspace{1mm} \mbox{for all fixed $t \in J$,}\nonumber\\
&& \partial^i_{(q,p)}f \in C( \T^n \times \R^m \times B \times J) \hspace{1mm} \mbox{for all $i \in\N^{2(n+m)}$ with $0 \le |i| \le [\sigma]$},\nonumber\\
\label{SMD}
&& \mbox{and} \hspace{1mm} \sup_{t \in J}|f^t|_{C^\sigma}t^l <\infty \Big\}
\end{eqnarray}
with the norm
\begin{equation}
\label{normMD}
|f|_{\sigma,l} = \sup_{t \in J}|f^t|_{C^\sigma}t^l,
\end{equation}
where we recall that $|\cdot|_{C^\sigma}$ is the Hölder norm (see Appendix \ref{A}).  We refer to Section \ref{FS} for some properties about $|\cdot|_{\sigma, l}$. We emphasize that we will use the same notation also for functions defined on $\T^n \times \R^m  \times J$, vector-valued functions or matrices.  This will be specified by the context.  

Now, we possess all the necessary instruments to state the first result of the present paper. Let $s$, $\lambda$, $\rho$, $\beta$ and $\alpha$ be positive parameters satisfying the following conditions
\begin{equation}
\label{parametersMD}
\begin{cases}
1 \le \rho <\lambda <s,\\
\displaystyle s>\max\left\{{\alpha \over \alpha -1}, \lambda + {\alpha \over \beta-1}\right\}, \quad 1<\beta <2, \quad \alpha >1,\\
\displaystyle \lambda >{2\beta \over 2-\beta}, \quad \rho<{\lambda - \beta \over \beta^2}.
 \end{cases}
\end{equation}
Given $e \in \R$, $\omega \in \R^n$ and real positive parameters $0 \le \delta <1$, $\varepsilon > 0$, and $\Upsilon \ge 1$, we consider the following time-dependent Hamiltonian
\begin{equation}
\label{HMD}
\begin{cases}
H : \T^n   \times \R^m \times B \times J \longrightarrow \R\\
H(q,p,t) = e + \omega \cdot p + a(q,t) + b(q,t) \cdot p + m(q,p,t) \cdot p^2\\
b(q,t) = b_0(q,t) + b_r(q,t)\\
a \in C(\T^n   \times \R^m \times J), \hspace{2mm} \partial_q a \in \mathcal{S}_{s,2}, \hspace{2mm}  b_0, b_r \in \mathcal{S}_{s+1, 1}, \hspace{2mm}  \partial_p^2 H  \in \mathcal{S}_{s+1, 0}\\
|b_0|_{2,1} \le \delta, \quad |b_0|_{s+1,1} \le \Upsilon,\\
|a|_{\lambda +1,0} + |\partial_q a|_{\lambda,2} < \varepsilon,  \quad |b_r|_{\lambda+1,1} < \varepsilon,\\
|a|_{s+1,0} + |\partial_q a|_{s,2} \le \Upsilon,  \quad |b_r|_{s+1,1} \le \Upsilon, \quad |\partial_p^2 H|_{s+1, 0} \le \Upsilon.
\end{cases}
\end{equation}
We also define the following trivial embedding $\varphi_0 : \T^n \times \R^m \to \T^n   \times \R^m \times B$ by $\varphi_0(q) = (q, 0)$ and the Hamiltonian $h : \T^n   \times \R^m \times B \times J \to \R$ 
\begin{equation}
\label{tildehBD}
h(q,p,t) = e + \left(\omega + b_0(q,t)\right) \cdot p  + m(q,p,t) \cdot p^2.
\end{equation}

\begin{theoremx}
\label{MD}
Let H be as in~\eqref{HMD} and we assume that $s$, $\lambda$, $\rho$, $\beta$ and $\alpha$ satisfy~\eqref{parametersMD}. Then, for $\delta$ small enough with respect to $s$, there exists $\varepsilon_0$, depending on $\delta$, $s$, $\lambda$, $\beta$, $\alpha$, $\rho$ and $\Upsilon$, such that for all $\varepsilon \le \varepsilon_0$ we have the existence of $v$, $\Gamma:\T^n \times \R^m \times J \to \R^{n+m}$, with $v$, $\Gamma \in \mathcal{S}_{\rho, 1}$, satisfying that, for all $t \in J$, 
\begin{equation*}
 \varphi^t(q) := (q, v^t(q))
\end{equation*}
 is a Lagrangian $C^\rho$-weakly asymptotic cylinder associated to $(X_H, X_{h}, \varphi_0)$, with $\Gamma$ the disturbing term associated to $\varphi^t$ (see Definition \ref{weakasymcyl}). Moreover
\begin{equation}
\label{StimevGamma}
|v|_{\rho, 1}<1, \quad |\Gamma|_{\rho, 1} <1.
\end{equation}
\end{theoremx}

Here, $b_0$ plays the role of $\gamma$ in Definition \ref{weakasymcyl}, contributing to the dynamics associated with the unperturbed Hamiltonian $h$ on the invariant cylinder $\varphi_0$.  We stress that we fix the parameter $\delta$ small enough in~\eqref{delta} and set $\varepsilon_0$ sufficiently small in Section \ref{PSMD}.

As mentioned before, the proof relies on a version of the Nash-Moser theorem due to Zehnder~\cite{Zeh76} (see Section \ref{Zehnder} for a brief introduction about the result of Zehnder and Theorem \ref{NMZ} for the version of the Nash-Moser theorem that we will use in the proof of Theorem \ref{MD}).  The author introduces a series of parameters satisfying suitable conditions (see~\eqref{T1} below)
in order to control the regularity of smooth functions and ensure the convergence of the iterative scheme employed in the proof of Theorem \ref{NMZ}.  Hypothesis~\eqref{parametersMD} is derived from~\eqref{T1} by setting $\upsilon=1$.  For the sake of clarity,  an example of parameters satisfying~\eqref{parametersMD} is the following
\begin{equation*}
\displaystyle s>\lambda + {7 \over 3}, \quad \beta = {3 \over 2}, \quad \alpha ={7 \over 6}, \quad
\displaystyle \lambda > 6, \quad \rho<{4 \over 9}\lambda - {2 \over 3}.
\end{equation*}
On the other hand, Zehnder proved that the minimum order required by $s$ to satisfy~\eqref{parametersMD} is $s \ge 8$ (see~\cite{Zeh76}).  For this value of $s$ the following set of parameters verifies~\eqref{parametersMD}
\begin{equation}
\label{parametriconsminimo}
s \ge 8, \quad \beta=1 + {7 \over 3s}, \quad  \alpha = {7 \over 6}, \quad \lambda = 2 + {14 \over s}, \quad \rho=1.
\end{equation}

For the regularity assumptions on the Hamiltonian $H$,  we want to emphasize that our proof does not work in the $C^\infty$ and analytic setting. In other words, we are not able to find $C^\infty$ or analytic weakly asymptotic cylinders.  The reason is that we cannot provide $C^\infty$ or analytic solutions to the homological equation~\eqref{HEMD} analyzed in Section \ref{SHE}.  The difficulties are explained in detail in Section \ref{HENCBD}, especially at the end of Section \ref{SHE}.  

For the sake of clarity, we want to point out that the Hamiltonian considered in~\cite{Sca22} is obtained by~\eqref{HMD} letting $m=0$ and $\delta=0$.  Given $l>2$, in this case we assume that $\partial_q a$ and $b$ decay as ${1\over t^l}$ and ${1 \over t^{l-1}}$, respectively, and we prove the existence of a $C^\sigma$-asymptotic KAM torus, which establishes the presence of orbits converging as time tends to infinity to the quasiperiodic solutions associated with the unperturbed system. This is because, in this case, the disturbing term $\Gamma^t$ in Definition \ref{weakasymcyl} is identically equal to zero.  The cost of this stronger result is a loss of two powers in the decay in time of the found $C^\sigma$-asymptotic KAM torus $\varphi^t$ (we recall that $\varphi^t$ converges in time to the invariant torus associated to the unperturbed system as ${1 \over t^{l-2}}$).

The proof of Theorem \ref{MD} is contained in Section \ref{ProofMD} and requires significant modification compared to the one in~\cite{Sca22}. Here, we introduce a suitable functional $\mathcal{F}$ related to the invariant equation~\eqref {WACD2} on specific  Banach spaces (we refer to~\eqref{WACD2} for the definition of the functional and to Section \ref{FS} for the definition of the Banach spaces).  The functional $\mathcal{F}$ depends on the perturbative terms $(a,b)$, the components of the $C^\sigma$-weakly asymptotic cylinder, and the disturbing term $\Gamma$ that we are looking for. It is defined in such a way that if $\mathcal{F}=0$, then~\eqref{WACD2} is satisfied.  In Section \ref{PSMD}, we verify that $\mathcal{F}$ satisfies the hypotheses of a version of the Nash-Moser theorem due to Zehnder~\cite{Zeh76}, which is recalled in Section \ref{Zehnder}. The presence of the disturbing term $\Gamma$ induces a loss of regularity, making the use of the Nash-Moser theorem necessary.  
We refer to Remark \ref{CompGMD} for a detailed comparison between the idea of the proof of the result in~\cite{Sca22} and the one of Theorem \ref{MD}. The associated linearized problem is analyzed in Section \ref{HENCBD}.  It is solved by integration thanks to a suitable change of coordinates that rectifies the dynamic on the cylinder. The most complicated part is verifying that the solution has a good decay in time. It is proved using tame estimates for the product and the composition of functions in Banach spaces of Hölder functions. 
It is considerably more complicated with respect to the linearized problem analyzed in~\cite{Sca22}; we refer to Remark \ref{CompHE} for an in-depth comparison.

\subsection{Planar three-body problem plus celestial body}\label{P3BP+CThm}

We consider three points of fixed masses $m_0$, $m_1$, and $m_2$ undergoing gravitational attraction in the plane and a celestial body of fixed mass $m_c$. The motion of the celestial body is a given continuous function $c(t)$ and only the planetary system is influenced by $c(t)$.  Moreover, we assume that the celestial body describes a hyperbolic motion, that is,
\begin{equation*}
\lim_{t \to +\infty}|c(t)| = + \infty,  \quad \lim_{t \to +\infty}{d \over dt}|c(t)| = v_c>0.
\end{equation*}
We denote $\R^{2*} = \R^2 \setminus \{(0,0)\}$.
Given $0<\varepsilon \le {1 \over 2}$ and $J = [1, +\infty)$, the phase space is 
\begin{equation}
\label{xyC}
\left\{((x_i,y_i) _{0 \le i\le 2}, t) \in \left(\R^2 \times \R^{2*}\right)^3 \times J \hspace{1mm}\Big|\hspace{1mm} \begin{matrix}\forall 0\le i < j \le 2, \hspace{1mm} x_i \ne x_j  \\
\forall 0\le i \le 2, \hspace{1mm}  {|x_i| \over |c(t)|} < \varepsilon \end{matrix}\right\}
\end{equation}
of linear momentum covectors $(y_0, y_1, y_2)$ and position vectors $(x_0, x_1, x_2)$ of each body. The Hamiltonian of the planar three-body problem plus celestial body (P3BP+CB) is
\begin{equation}
\label{Hc}
H(x,y,t) = H_0(x,y) +H_c(x,t),
\end{equation}
with
\begin{equation*}
H_0(x,y) = \sum_{i=0}^2 {|y_i|^2 \over 2m_i} - G\sum_{0 \le i<j\le 2}{m_i m_j \over |x_i - x_j|}, \quad H_c(x,t) = - G\sum_{i=0}^2 {m_i m_c \over |x_i - c(t)|}
\end{equation*}
where $G$ is the universal constant of gravitation that we may, in suitable unities,  suppose equal to $1$.
$H$ is the sum of the Hamiltonian of the planar three-body problem $H_0$ and the Hamiltonian of the interaction with the celestial body $H_c$. Consider $\phi_0$ as a $1$-parameter family of invariant tori for $H_0$ supporting quasiperiodic dynamics with four frequencies, whose existence is guaranteed by~\cite{Fe02}. We denote $M = m_0 + m_1+m_2$ and $\psi^t_{t_0, H}$ as the flow at time $t$ with initial time $t_0$ of $H$.  

\begin{theoremx}
\label{Thmcomet}
Let $H$ be as in~\eqref{Hc} and $\phi_0$ be a $1$-parameter family of invariant tori for $H_0$ supporting quasiperiodic dynamics with four frequencies.  Then, there exist constants $\Upsilon_{2, \phi}$ depending on $\phi_0$ and $\Upsilon_{\phi,M,m_c}$ depending on $\phi_0$, $M$, and $m_c$ such that if
\begin{equation}
\label{conditioncv}
|c(1)| > { \max\{1, 3\Upsilon_{2, \phi}\} \over \varepsilon}, \quad v_c >{12 \over \varepsilon}
\end{equation}
for $\varepsilon$ small enough with respect to $\Upsilon_{\phi,M,m_c}$,  then there exists a non-empty open subset $\mathcal{W} \subset \left(\R^2 \times \R^{2*}\right)^3$ such that, for all $x \in \mathcal{W}$, $\psi^t_{1, H}(x)$ is a weakly asymptotically quasiperiodic solution associated to $(X_H, X_{H_0}, \phi_0)$. 
\end{theoremx}

We stress that the constants $\Upsilon_{2, \phi}$ and $\Upsilon_{\phi,M,m_c}$ are explicit and they are defined by ~\eqref{k1} and~\eqref{UpsilonPhiPi}, respectively.

The proof follows the following steps.  In Sections \ref{P1C} and \ref{SecQPDKC}, we introduce symplectic coordinates $(\theta, \xi, r, \eta) \in \T^4 \times \R^2 \times \R^4 \times \R^2$ to place the invariant cylinder $\phi_0$ onto the zero section.  Section \ref{P3C} defines a subset $\mathcal{U}$ of the phase space (see~\eqref{U}) characterized by the orbits that at each time stay sufficiently far from the celestial body.
We establish that, for $|c(1)|$ and $v_c$ large enough, on the subset $\mathcal{U}$, the perturbation $H_c^t$ satisfies good decay properties (see Lemma \ref{Ux} and Lemma \ref{LemmaHcstime}).  In Section \ref{P4C} we introduce a suitable smooth extension $\hat H^t$ of the Hamiltonian $H^t$ of the P3BP+CB because we need the Hamiltonian to be defined in the entire phase space in order to apply Theorem \ref{MD} (we recall that $H_c^t$ satisfies good decay properties just in the subset $\mathcal{U}$ and not in the entire phase space).    
The Hamiltonian $\hat H^t$ coincides with $H^t$ on a suitable subset of $\mathcal{U}$, where we expect the motion to take place (we refer to~\eqref{HatH} for the definition of $\hat H^t$). The application of Theorem \ref{MD} ensures the existence of a $C^1$-weakly asymptotic cylinder for the Hamiltonian $\hat H^t$.  In the last part of the proof, contained in Section \ref{P5C}, we define an appropriate subset $\mathcal{W} \subset  \left(\R^2 \times \R^{2*}\right)^3$ and we verify that, for $|c(1)|$ and $v_c$ large enough, each initial point $x \in \mathcal{W}$ gives rise to a weakly asymptotically quasiperiodic solution associated to $(X_H, X_{H_0}, \phi_0)$. 
The reason is that, for $|c(1)|$ and $v_c$ large enough,  we will prove that an orbit associated to $\hat H^t$ with initial condition in $\mathcal{W}$ stays sufficiently far from the celestial body for all $t \in J$ and hence it remains in the subset of the phase space where $\hat H^t$ coincides with $H^t$. 

We point out that Theorem \ref{Thmcomet} provides the existence of orbits such that the center of mass of the planetary system is attracted by the celestial body $c(t)$ with a zero asymptotic velocity. Whereas,  in a frame of reference attached to the center of mass of the planetary system, the motions of the planets converge in time to some orbits close to quasiperiodic solutions in the sense of Definition \ref{asymsolI}.

\section{Functional setting}\label{FS}

The first part of this section is dedicated to analyzing the definition of $C^\sigma$-weakly asymptotic cylinder (see Definition \ref{weakasymcyl}), providing a series of properties.  For this purpose, we recall that $B \subset \R^{n+m}$ is an open ball centered at the origin and $J=[1, +\infty) \subset \R$.  Let $X$, $X_0$ and $\varphi_0$ be as in Definition \ref{weakasymcyl} and we denote by $\varphi^t$ a $C^\sigma$-weakly asymptotic cylinder associated to $(X, X_0, \varphi_0)$.  In the following proposition, we will see that we can rewrite the invariant equation~\eqref{hyp1WAC} in terms of the flow of $X$.  Let $\psi^t_{t_0,X}$ and $\psi^t_{t_0,\bar \omega + \Gamma}$ be the flow at time $t$ with initial time $t_0$ of $X$ and $\bar \omega + \Gamma$, respectively.
\begin{proposition}
\label{IntroP1}
If the flows $\psi^t_{t_0,X}$ and $\psi^t_{t_0,\bar \omega + \Gamma}$ are defined for all $t$, $t_0 \in J$, then~\eqref{hyp1WAC} is equivalent to 
\begin{equation}
\label{hyp1WACbiss}
\psi^t_{t_0, X} \circ \varphi^{t_0}(q) = \varphi^t \circ \psi_{t_0, \bar \omega + \Gamma}^t(q),
\end{equation} 
for all $t$, $t_0 \in J$ and $q \in \T^n \times \R^m$.
\end{proposition}
\begin{proof}
The proof is straightforward; it suffices to prove that both sides of~\eqref{hyp1WACbiss} satisfy the same initial value problem. 
\end{proof}
Using the previous proposition, we can see that it is always possible to find a family of embeddings that satisfies~\eqref{hyp1WAC}.

\begin{proposition}
\label{IntroP2}
If $\psi^t_{t_0,X}$ and $\psi^t_{t_0,\bar \omega + \Gamma}$ are defined for all $t$, $t_0 \in J$, it is always possible to find a family of embeddings satisfying~\eqref{hyp1WAC}
\end{proposition}
\begin{proof}
We consider an embedding $\hat \varphi : \T^n \times \R^m \to \T^n \times \R^m \times B$.  Furthermore, for all $t$, $t_0 \in J$ and $q \in \T^n \times \R^m$, we define the following family of embeddings
\begin{equation*}
\varphi^t(q) = \psi_{t_0, X}^t \circ \hat \varphi \circ \psi_{t, \bar \omega + \Gamma}^{t_0}(q).
\end{equation*} 
The latter is a family of embeddings satisfying~\eqref{hyp1WACbiss} and hence~\eqref{hyp1WAC}.
\end{proof}

The proposition below asserts that if we have the existence of a $C^\sigma$-weakly asymptotic cylinder defined for all $t$ large, then we can extend the set of definition for all $t \in \R$.

\begin{proposition}
\label{IntroP3}
We assume that $\psi^t_{t_0,X}$ and $\psi^t_{t_0,\bar \omega + \Gamma}$ are defined for all $t$, $t_0 \in \R$.  Given $T \ge 0$, if there exists a $C^\sigma$-weakly asymptotic KAM torus $\varphi^t$ defined for all $t \ge T$, then we can extend the set of definition for all $t \in \R$.
\end{proposition}
\begin{proof}
For all $q \in \T^n \times \R^m$, we consider
\begin{equation*}
\phi^t(q) = \begin{cases} \varphi^t(q) \hspace{4mm} \mbox{for all $t \ge T$}\\
\psi_{T, X}^t \circ \varphi^T\circ \psi_{t, \bar \omega + \Gamma}^{T}(q) \hspace{4mm} \mbox{for all $t \le T$}.\end{cases}
\end{equation*}
The above family of embeddings verifies~\eqref{hyp2WAC},~\eqref{hyp1WAC} and~\eqref{hyp3WAC}.
\end{proof}

The last proposition of the first part of this section provides some information concerning the dynamics associated with a $C^\sigma$-weakly asymptotic cylinder.  Letting $X$, $X_0$ and $\varphi_0$ be as in Definition \ref{weakasymcyl}, we will see that if there exists a $C^\sigma$-weakly asymptotic cylinder associated to $(X, X_0, \varphi_0)$, then we have the existence of weakly asymptotic solutions associated to $(X, X_0, \varphi_0)$ (see Definition \ref{asymsolI}).

\begin{proposition}
\label{propasymsolI}
Let $\varphi^t$ be a $C^\sigma$-weakly asymptotic cylinder associated to $(X, X_0,\varphi_0)$. 
Then, for all $q \in \T^n \times \R^m$ and $t_0 \in J$, 
\begin{equation*}
g(t) = \psi_{t_0, X}^t \circ \varphi^{t_0}(q)
\end{equation*}
is a weakly asymptotically quasiperiodic solution associated to $(X, X_0, \varphi_0)$. 
\end{proposition}
\begin{proof}
The proof is a straightforward consequence of~\eqref{hyp2WAC} and~\eqref{hyp1WACbiss}.
\end{proof}

The second part of this section is dedicated to a series of properties of the norm~\eqref{normMD} introduced in Section \ref{AbsThm}.  For this reason,  we recall the definition of the Banach space $\left(\mathcal{S}_{\sigma, l}, |\cdot|_{\sigma, l}\right)$. Given positive parameters $\sigma$, $l \ge 0$
\begin{eqnarray*}
\mathcal{S}_{\sigma, l} &=& \Big\{ f:\T^n \times \R^m \times B \times J \to \R \hspace{1mm} | \hspace{1mm} f^t \in C^\sigma( \T^n \times \R^m \times B) \hspace{1mm} \mbox{for all fixed $t \in J$,}\\
&& \partial^i_{(q,p)}f \in C( \T^n \times \R^m \times B \times J) \hspace{1mm} \mbox{for all $i \in\N^{2(n+m)}$ with $0 \le i \le [\sigma]$},\\
&& \mbox{and} \hspace{1mm} \sup_{t \in J}|f^t|_{C^\sigma}t^l <\infty \Big\}.
\end{eqnarray*}
whereas the associated norm is given by
\begin{equation*}
|f|_{\sigma,l} = \sup_{t \in J}|f^t|_{C^\sigma}t^l.
\end{equation*}
The following proposition contains various properties of the latter norm.  To avoid a flow of constants,  in the present paper, we will denote by $C(\cdot)$ constants depending on $n + m$ and the parameters in brackets, while $C$ stands for constants depending only on $n + m$. 

\begin{proposition}
\label{normpropertiesMD}
Given $\sigma$, $l$, $d >0$, for all $f \in \mathcal{S}_{\sigma, l}$ and $g \in \mathcal{S}_{\sigma, d}$ we have the following properties
\begin{enumerate}
\item For all $s>0$ and  $\beta \in \N^{2(n + m)}$, if $|\beta| + s \le \sigma$, then  
\begin{equation*}
\left|{\partial^{\beta} \over \partial{q_1}^{\beta_1}... \partial{q_{2n}}^{\beta_{2n}} \partial{p_1}^{\beta_{2n+1}}... \partial{p_n}^{\beta_{2(n +m)}}} f \right|_{s, l} \le C |f|_{\sigma, l}
\end{equation*}
\item For all $l' \ge 0$, $|f|_{\sigma, l}  \le |f|_{\sigma, l+l'}$, 
\item $|fg|_{\sigma, l+d} \le C(\sigma)\left(|f|_{0,l}|g|_{\sigma,d} + |f|_{\sigma,l}|g|_{0,d}\right)$. 
\item We consider $\sigma \ge 1$,  and we assume that $g:\T^n \times \R^m \times B \times J \to \T^n \times \R^m \times B$.  Letting $\tilde g: \T^n \times \R^m \times B \times J \to \T^n \times \R^m \times B \times J$ such that $\tilde g(q,p,t) = (g(q,p,t), t)$, then $f \circ \tilde g \in \mathcal{S}_{\sigma, l+d}$ and 
\begin{equation*}
|f \circ \tilde g|_{\sigma, l+d} \le C(\sigma) \left(|f|_{\sigma,l}|\partial_{(q,p)} g|^\sigma_{0,d} + |f|_{1,l}|\partial_{(q,p)}  g|_{\sigma-1,d} +  |f|_{0, l+d}  \right)
 \end{equation*}
 where $\partial_{(q,p)}$ stands for the partial derivatives with respect to the coordinates $(q,p) \in \T^n \times \R^m \times B$.
\end{enumerate}
\end{proposition}
\begin{proof}
The proof of property \textit{1} consists in a straightforward application of property \textit{1} of Proposition \ref{Holder}. 
Property \textit{2} is obvious. Now, we prove  \textit{3}. Thanks to property \textit{3} of Proposition \ref{Holder}, one has
\begin{eqnarray*}
 |f g|_{\sigma, l+d} &=& \sup_{t \in J}|f^t  g^t|_{C^\sigma}t^{l+d} \le C(\sigma)\sup_{t \in J}\left(|f^t|_{C^0}|g^t|_{C^\sigma} + |f^t|_{C^\sigma}|g^t|_{C^0}\right)t^{l+d} \\
&\le& C(\sigma) \sup_{t \in J}\left(|f^t|_{C^0}t^l|g^t|_{C^\sigma}t^d + |f^t|_{C^\sigma}t^l|g^t|_{C^0}t^d\right)\\
&\le& C(\sigma)\left(|f|_{0,l}|g|_{\sigma,d} + |f|_{\sigma,l}|g|_{0,d}\right).
\end{eqnarray*}
This concludes the proof of \textit{3}.  In the proof of \textit{4} we use property \textit{5} of Proposition \ref{Holder}
\begin{eqnarray*}
 |f\circ \tilde g|_{\sigma, l+d} &=& \sup_{t \in J}|f^t \circ g^t|_{C^\sigma}t^{l+d} \\
&\le& C(\sigma) \sup_{t \in J}\left(|f^t|_{C^\sigma}|\partial_{(q,p)} g^t|^\sigma_{C^0} + |f^t|_{C^1}|\partial_{(q,p)} g^t|_{C^{\sigma-1}} + |f|_{C^0}\right)t^{l+d} \\
&\le& C(\sigma) \sup_{t \in J}\left(|f^t|_{C^\sigma}t^l|\partial_{(q,p)} g^t|^\sigma_{C^0}t^{\sigma d}{t^{(1-\sigma)d}} + |f^t|_{C^1}t^l|\partial_{(q,p)} g^t|_{C^{\sigma-1}}t^d\right)\\
 &+& C(\sigma) \sup_{t \in J} |f|_{C^0}t^{l+d}\\
&\le& C(\sigma) \left(|f|_{\sigma,l}|\partial_{(q,p)} g|^\sigma_{0,d} + |f|_{1,l}|\partial_{(q,p)} g|_{\sigma-1,l} +  |f|_{0, l+d}  \right)
\end{eqnarray*}
notice that $t \ge 1$ and $1-\sigma \le 0$ imply $t^{(1-\sigma)d} \le 1$. 
\end{proof}

In this third part, we introduce notation and Banach spaces that will be employed in the proofs of the main theorems of this paper. For this purpose,  given $f:\T^n \times \R^m \times B\times J\to \R$ we denote the function $\tilde{f}$ as follows
\begin{equation}
\label{tildef}
\tilde f :\T^n \times \R^m \times B \times J \to \R \times J,\hspace{1mm} \mbox{such that} \hspace{1mm} \tilde f (q,p,t) = (f(q,p,t), t).
\end{equation}
On the other hand,  for a given $\omega \in \R^n$,  we define by $\bar \omega$ and $\bar \Omega$ the following vectors 
\begin{equation}
\label{BarOmega}
\bar \omega = (\omega, 0) \in \R^{n+m}, \quad \bar \Omega = (\bar \omega, 1) \in \R^{n+m+1}.
\end{equation}
We consider $u:\T^n \times \R^m \times J \to \R^{n+m}$, $\omega \in \R^n$ and we denote by $Du$ the differential of $u$. Then,  using the above notation~\eqref{BarOmega}, one can see that 
\begin{equation}
Du(q,t) \bar \Omega = \partial_q u(q,t) \bar \omega + \partial_t u(q,t)
\end{equation}
for all $(q,t) \in \T^n \times \R^m \times J$. 

In the proof of Theorem \ref{MD}, we need specific Banach spaces to control the various elements of the dynamic problem we aim to solve. For this reason,  we fix $\sigma$, $l \ge 0$ and $\omega \in \R^n$.  In order to quantify the regularity of the components of the $C^\sigma$-weakly asymptotic cylinder we will look for, we define the following Banach space
\begin{equation}
\label{V}
\mathcal{V}_{\sigma, \omega, l} = \Big\{ v:\T^n \times \R^m \times J \to \R^{n+m} \hspace{1mm} | \hspace{1mm} v \in \mathcal{S}_{\sigma+1, l}, Dv\bar \Omega \in \mathcal{S}_{\sigma, l+1} \Big \}
\end{equation}
with the norm
\begin{equation}
\label{NormV}
|v|_{\sigma, \omega, l} = \max \{|v|_{\sigma+1, l},  |Dv\bar \Omega|_{\sigma, l+1}\}
\end{equation}
To characterize the regularity of the constant term in the time-dependent perturbation of~\eqref{HMD}, we introduce a Banach space as follows
\begin{eqnarray}
\mathcal{S}_{\sigma, (0,l)} &=& \Big\{ g:\T^n \times \R^m \times J \to \R \hspace{1mm} | \hspace{1mm} g \in C(\T^n \times \R^m \times J) \hspace{1mm},  \partial_q g \in \mathcal{S}_{\sigma, l},  \nonumber\\
\label{S2}
&&\mbox{and} \hspace{1mm} |g|_{\sigma+1, 0}< \infty\Big \}
\end{eqnarray}
whereas the associated norm is given by
\begin{equation}
\label{NormS2}
|g|_{\sigma, (0,l)} = |g|_{\sigma+1, 0} + |\partial_q g|_{\sigma, l}.
\end{equation}
We point out that in the definition of the above norms~\eqref{NormV} and~\eqref{NormS2}, we used the norm $|\cdot|_{\sigma, l}$ defined by~\eqref{normMD}. 

In the final part of this section, we prove some properties of the above Banach spaces. As mentioned before, one of the main instruments in the proof of Theorem \ref{MD} is the version of the Nash-Moser theorem proved by Zehnder~\cite{Zeh76}.  Zehnder consider one-parameter families of Banach spaces $\{(\mathcal{X}^\sigma, |\cdot|_\sigma)\}_{\sigma \ge 0}$ such that for all $0 \le \sigma' \le \sigma < \infty$
\begin{equation}
\label{Xsigma}
\begin{aligned}
\mathcal{X}^0 \supseteq \mathcal{X}^{\sigma'} &\supseteq \mathcal{X}^\sigma \supseteq \mathcal{X}^\infty= \bigcap_{\sigma \ge 0} \mathcal{X}^\sigma\\
&|x|_{\sigma'} \le |x|_{\sigma} 
\end{aligned}
\end{equation}
for all $x \in \mathcal{X}^\sigma$.  The author needs that each of these families of Banach spaces possess a $C^\infty$-smoothing (see Definition \ref{smoothing} below).
\begin{definition}
\label{smoothing}
A $C^\infty$-smoothing in $\{(\mathcal{X}^\sigma, |\cdot|_\sigma)\}_{\sigma \ge 0}$ is a one-parameter family $\{S_\tau\}_{\tau >0}$ of linear mappings $S_\tau :\mathcal{X}^0 \to \mathcal{X}^\infty$ together with constants $C(r,d)$, for positive integers $r$ and $d$, satisfying the following conditions:
\begin{equation}
\label{S1}
|S_\tau x|_r \le \tau^{r-d}C(r,d)|x|_d
\end{equation}
for all $x \in \mathcal{X}^d$ and $0 \le d \le r$, 
\begin{equation}
\label{S2}
|(S_\tau - 1)x|_d \le \tau^{-(r-d)}C(r,d)|x|_r
\end{equation}
for all $x \in \mathcal{X}^r$ and $0 \le d \le r$.
\end{definition}

Now, we fix $l \ge 0$ and $\omega \in \R^n$. Using the Banach spaces introduced in this section, we consider the following families of Banach spaces $\{(\mathcal{S}_{\sigma,l}, |\cdot|_{\sigma,l})\}_{\sigma \ge 0}$, $\{(\mathcal{S}_{\sigma, (0,l)}, |\cdot|_{\sigma, (0,l)})\}_{\sigma \ge 0}$ and $\{(\mathcal{V}_{\sigma, \omega, l}, |\cdot|_{\sigma, \omega, l})\}_{\sigma \ge 0}$. 

It is straightforward to verify that, for all $0 \le \sigma' \le \sigma <\infty$,
\begin{equation*}
\mathcal{S}_{0, l} \supseteq \mathcal{S}_{\sigma', l} \supseteq  \mathcal{S}_{\sigma, l} \supseteq  \mathcal{S}_{\infty, l} = \bigcap_{\sigma \ge 0}  \mathcal{S}_{\sigma, l} , \hspace{6mm} \mathcal{V}_{0, \omega, l} \supseteq \mathcal{V}_{\sigma', \omega, l} \supseteq \mathcal{V}_{\sigma, \omega, l} \supseteq \mathcal{V}_{\infty, \omega, l}= \bigcap_{\sigma \ge 0} \mathcal{V}_{\sigma, \omega, l},
\end{equation*}
\begin{equation*}
|f|_{\sigma', l} \le |f|_{\sigma, l} \hspace{12mm} |v|_{\sigma', \omega, l} \le |v|_{\sigma, \omega, l} 
\end{equation*}
\vspace{3mm}
\begin{equation*}
\mathcal{S}_{0, (0,l)} \supseteq \mathcal{S}_{\sigma', (0,l)} \supseteq  \mathcal{S}_{\sigma, (0,l)} \supseteq  \mathcal{S}_{\infty, (0,l)} = \bigcap_{\sigma \ge 0}  \mathcal{S}_{\sigma, (0,l)},
\end{equation*}
\begin{equation*}
 |g|_{\sigma', (0,l)} \le |g|_{\sigma, (0,l)}
\end{equation*}
for all $f \in \mathcal{S}_{\sigma, l}$, $g \in \mathcal{S}_{\sigma, (0,l)}$ and $v \in \mathcal{V}_{\sigma, \omega, l}$.

The following lemma proves the existence of a $C^\infty$-smoothing for these families of Banach spaces. This is not surprising because the behavior of these norms is very similar to that of the Hölder norms.

\begin{lemma}
\label{St}
There exists a $C^\infty$-smoothing for the latter families of Banach spaces. 
\end{lemma}
\begin{proof}
We begin by proving the existence of a $C^\infty$-smoothing for the family of Banach spaces $\{(\mathcal{S}_{\sigma,l}, |\cdot|_{\sigma,l})\}_{\sigma \ge 0}$. Following the lines of~\cite{Zeh76}, we take a function $\tilde s \in C^\infty_0(\R^{n+m})$  vanishing outside a compact set and identically equal to $1$ in a neighborhood of $0$. Let $s$ be its Fourier transform.  For all $f : \T^n \times \R^m \times B \times J \to \R$ with $f \in \mathcal{S}_{0, l}$,  letting $x = (q, p) \in  \T^n \times \R^m \times B$, we define  
\begin{equation}
\label{Stf}
S_\tau f(x,t) = {1 \over \tau^{2(n+m)}} \int_{\R^{2(n+m)}}s\left({x - \vartheta\over \tau}\right) f( \vartheta,t) d \vartheta. 
\end{equation}
In the first part of this proof, we want to prove that for all $f \in \mathcal{S}_{0, l}$,  one has $S_\tau f \in  \mathcal{S}_{\infty, l}$.  To this end, we observe that,  $S_\tau z^t\in C^\infty(\T^n \times \R^m \times B) = \bigcap_{\sigma \ge 0} C^\sigma(\T^n \times \R^m\times B)$ for all fixed $t\in J$ (we refer to~\cite{Zeh76}). Now, we need to verify that $\partial_x^i\left(S_\tau z\right)\in C(\T^n \times \R^m \times B \times J)$  for all $i \in \N^{2(n+m)}$.  For this reason, we observe that, for every $p >0$ and $j \in\N^{2(n+m)}$ with $|j|>0$, there exists a constant $C(|j|,p)>0$ such that 
\begin{equation}
\label{stimes}
|\partial^j s(x)| \le C(|j|,p) (1 + |x|)^{-p},
\end{equation}
where $\partial^j$ stands for partial derivatives of order $j$ (see~\cite{Zeh76}) and we recall that $|j| = j_1 +...+j_{2(n+m)}$. The claim will be a consequence of the regularity of $f$ and the latter. Indeed, for all $(x_1, t_1)$, $(x_2, t_2) \in \T^n \times \R^m \times B \times J$ and $i \in \N^{2(n+m)}$,
\begin{align*}
&\Big|\partial_x^i\left(S_\tau f\right)(x_1,t_1) - \partial^i_x\left(S_\tau f\right)(x_2, t_2)\Big| \\
&= \Big|{1 \over \tau^{2(n+m)+|i|}} \int_{\R^{2(n+m)}}\partial^is\left({x_1 -  \vartheta\over \tau}\right) f( \vartheta,t_1)d \vartheta\\
&-  {1 \over \tau^{2(n+m)+|i|}} \int_{\R^{2(n+m)}}\partial^i s\left({x_2 -  \vartheta\over \tau}\right) f( \vartheta,t_2)d \vartheta  \Big|\\
&={1 \over \tau^{|i|}} \Big|\int_{\R^{2(n+m)}}\partial^is(\rho) \big(f(x_1 - \rho \tau,t_1) - f(x_2 - \rho \tau,t_2) \big) d\rho\Big|\\
&\le{1 \over \tau^{|i|}}\int_{\R^{2(n+m)}}\left|\partial^is(\rho)\right|\left|\big(f(x_1 - \rho \tau,t_1) - f(x_2 - \rho \tau,t_2) \big)\right| d\rho
\end{align*}
where $|\cdot|$ stands for the standard Euclidean norm and, in agreement with the convention made above,  $|i|= i_1+...+i_{2(n+m)}$ if $i \in \N^{2(n+m)}$. In the last line of the latter, we did the following change of coordinates ${x_i - \vartheta \over \tau} = \rho$ for $i = 1,2$.  
Thanks to the latter and~\eqref{stimes}, one can prove that $\partial_x^iS_\tau f \in C(\T^n \times \R^m \times B \times J)$ for all $i \in \N^{2(n+m)}$.  
This concludes the first part of the proof.  In the second part, we show that $S_\tau$ (see~\eqref{Stf}) satisfies~\eqref{S1} and~\eqref{S2}.  For this purpose,  we observe that for all $f \in \mathcal{S}_{d, l}$, $0 \le d \le r$ and fixed $t \in J$ 
\begin{equation*}
|S_\tau f^t|_{C^r} \le \tau^{r-d} C(r, d) |f^t|_{C^d}
\end{equation*}
(always look at~\cite{Zeh76}). Thanks to the above inequality,  for all $f \in \mathcal{S}_{d, l}$ and $0 \le d \le r$ we have that
\begin{equation*}
|S_\tau f|_{r, l} = \sup_{t \in J}|S_\tau z^t|_{C^r}t^l \le \tau^{r-d} C(r, d) \sup_{t \in J}|f^t|_{C^d}t^l = \tau^{r-d} C(r, d) |f|_{d, l}.
\end{equation*}
Hence~\eqref{S1} is verified.  It remains to prove~\eqref{S2}. For this reason,  we can see that for all $f\in \mathcal{S}_{r, l}$,  $0 \le d \le r$ and fixed $t \in J$
\begin{equation*}
|(S_\tau - 1) f^t|_{C^d} \le \tau^{-(r-d)} C(r, d) |f^t|_{C^r}
\end{equation*}
(always see~\cite{Zeh76}). Similarly to the previous case,  for all $f\in \mathcal{S}_{r, l}$ and  $0 \le d \le r$ we have
\begin{eqnarray*}
|(S_\tau - 1) f|_{d, l}&=&\sup_{t \in J} |(S_\tau - 1) f^t|_{C^d}t^l \le \tau^{-(r-d)} C(r, d) \sup_{t \in J}|f^t|_{C^r}t^l \\
&\le& \tau^{-(r-d)} C(r, d) |f|_{r, l}.
\end{eqnarray*}
Then~\eqref{S2} is also verified. This concludes the proof about the existence of a $C^\infty$-smoothing for $\{(\mathcal{S}_{\sigma,l}, |\cdot|_{\sigma,l})\}_{\sigma \ge 0}$. Remembering that $S_\tau$ commutes with partial differential operators, similarly,  one can prove the existence of a $C^\infty$-smoothing for the family of Banach spaces $\{(\mathcal{S}_{\sigma, (0,l)}, |\cdot|_{\sigma, (0,l)})\}_{\sigma \ge 0}$.

It remains to show the existence of a $C^\infty$-smoothing for  $\{(\mathcal{V}_{\sigma, \omega, l}, |\cdot|_{\sigma, \omega, l})\}_{\sigma \ge 0}$. Similarly to the previous case, using that $S_\tau$ commutes with partial differential operators,  ona can prove that $S_\tau : \mathcal{V}_{0, \omega, l} \to \mathcal{V}_{\infty, \omega,l}$ is well-defined.

Now, we verify~\eqref{S1} and~\eqref{S2}. We begin by remembering that, for all $v \in \mathcal{V}_{\sigma, \omega, l}$,
\begin{equation*}
|v|_{\sigma, \omega, l} =\max \{|v|_{\sigma+1,l},| Dv\bar\Omega|_{\sigma,l+1} \},
\end{equation*}
(see~\eqref{NormV}). Similarly to the previous case, for all $v \in \mathcal{V}_{d, \omega, l}$, $0 \le d \le r$ and fixed $t \in J$ one has that 
\begin{eqnarray*}
|S_\tau v^t|_{C^{r+1}} \le \tau^{r-d} C(r, d) |v^t|_{C^{d+1}},
\end{eqnarray*}
which implies that for all $v \in \mathcal{V}_{d, \omega, l}$ and $0 \le d \le r$
\begin{eqnarray}
|S_\tau v|_{r+1, l} &=& \sup_{t \in J}|S_\tau v^t|_{C^{r+1}}t^l \le \tau^{r-d} C(r, d) \sup_{t \in J}|v^t|_{C^{d+1}} t^l \nonumber\\ \label{SSmoothing1}
&\le& \tau^{r-d} C(r, d) |v|_{d, \omega, l}.
\end{eqnarray}
Noting that $S_\tau$ commutes with partial differential operators,  for all $v \in \mathcal{V}_{d, \omega, l}$, $0 \le d \le r$ and for fixed $t \in J$,
\begin{eqnarray*}
|D\left(S_\tau v^t\right) \bar\Omega|_{C^r} &=& |S_\tau D  v^t \bar\Omega|_{C^r} \le \tau^{r-d} C(r, d)|D  v^t \bar\Omega|_{C^d}.
\end{eqnarray*}
Multiplying both sides of the latter by $t^{l+1}$ and taking the sup for all $t \in J$ we can see that 
\begin{eqnarray}
|D\left(S_\tau v\right) \bar\Omega|_{r, l+1} &=& \sup_{t \in J}|D\left(S_\tau v^t\right) \bar\Omega|_{C^r}t^{l+1} \le \tau^{r-d} C(r, d) \sup_{t \in J}|D  v^t\bar\Omega|_{C^d}t^{l+1}\nonumber \\ \label{SSmoothing2}
&\le&  \tau^{r-d} C(r, d)| v |_{d, \omega, l}
\end{eqnarray}
for all $v \in \mathcal{V}_{d, \omega, l}$ and $0 \le d \le r$.
Now, thanks to~\eqref{SSmoothing1} and~\eqref{SSmoothing2} we can prove that for all $v \in \mathcal{V}_{d, \omega, l}$ and $0 \le d \le r$
\begin{eqnarray*}
| S_\tau v|_{r, \omega, l} &=& \max \{|S_\tau v|_{r+1, l}, |D\left(S_\tau v\right) \bar\Omega|_{r, l+1} \} 
\le  \tau^{r-d} C(r, d)| v |_{d, \omega ,l}
\end{eqnarray*}
and hence~\eqref{S1} is satisfied. 

Concerning~\eqref{S2}, for all $v \in \mathcal{V}_{r, \omega, l}$,  $0 \le d \le r$ and fixed $t \in J$, we have that 
\begin{equation*}
|(S_\tau - 1) v^t|_{C^{d+1}} \le \tau^{-(r-d)} C(r, d) |v^t|_{C^{r+1}}.
\end{equation*}
Multiplying both sides of the latter by $t^l$ and taking the sup for all $t \in J$, we can see that 
\begin{eqnarray}
|(S_\tau - 1) v|_{d+1,l}&=&\sup_{t \in J} |(S_\tau - 1) v^t|_{C^{d+1}}t^l \le \tau^{-(r-d)} C(r, d) \sup_{t \in J}|v^t|_{C^{r+1}}t^l  \nonumber\\
\label{S2Smoothing}
&\le& \tau^{-(r-d)} C(r, d) |v|_{r+1,l} \le \tau^{-(r-d)} C(r, d) |v|_{r, \omega, l}
\end{eqnarray}
for all $v \in \mathcal{V}_{r, \omega, l}$ and $0 \le d \le r$.
On the other hand,  for all $v \in \mathcal{V}_{r, \omega, l}$,  $0 \le d \le r$ and for fixed $t \in J$,
\begin{eqnarray*}
|D(S_\tau -1)v^t \bar\Omega|_{C^d} &=& |(S_\tau -1) D  v^t\bar \Omega|_{C^d} \le \tau^{-(r-d)} C(r, d)|D v^t \bar\Omega|_{C^r}
\end{eqnarray*}
and hence  for all $v \in \mathcal{V}_{r, \omega, l}$,  and $0 \le d \le r$
\begin{eqnarray}
|D \left((S_\tau -1)v\right)  \bar\Omega|_{d,l+1} &=& \sup_{t\in J} |D\left((S_\tau -1)v^t\right)\bar \Omega|_{C^d}t^{l+1} \nonumber\\
&\le& \tau^{-(r-d)} C(r, d)\sup_{t\in J}|D  v^t\bar\Omega|_{C^r}t^{l+1}\nonumber \\
&\le&\tau^{-(r-d)} C(r, d)|D v \bar\Omega|_{r,l+1}\nonumber \\
\label{S2Smoothing2}
&\le&  \tau^{-(r-d)} C(r, d) |v|_{r, \omega, l}.
\end{eqnarray}
This concludes the proof of this lemma because thanks to~\eqref{S2Smoothing} and~\eqref{S2Smoothing2} we have that for all $v \in \mathcal{V}_{r, \omega, l}$,  and $0 \le d \le r$
\begin{eqnarray*}
|\left(S_\tau -1\right) v |_{d, \omega , l} &=&\max\{|(S_\tau - 1) v|_{d+1,l},  |D \left((S_\tau -1)v \right) \bar\Omega|_{d, l+1} \}\\
&\le&  \tau^{-(r-d)} C(r, d) |v|_{r, \omega, l}.
\end{eqnarray*}
\end{proof}


\section{Homological equation}\label{HENCBD}

Given $\sigma \ge 1$, $\mu \ge 0$, and $\omega \in \R^n$, this section is devoted to solving the following equation for the unknown $\varkappa : \T^n \times \R^m \times J \longrightarrow \R^{n + m}$  
\begin{equation}
\label{HEMD}
\begin{cases}
\partial_q\varkappa(q,t) \underbrace{\left(\bar \omega + f(q,t)\right)}_{F(q,t)} + \partial_t \varkappa(q,t) + g(q,t)\varkappa(q,t) = z(q,t)\\
f, g\in \mathcal{S}_{\sigma,1}, \quad z\in \mathcal{S}_{\sigma,2}, \quad |f|_{1,1} \le \mu, \quad |g|_{1,1} \le \mu, 
\end{cases}
\end{equation}
where we recall that $\bar \omega = (\omega, 0) \in \R^{n+m}$ (see~\eqref{BarOmega}), the Banach space $\mathcal{S}_{\sigma, l}$ is defined by~\eqref{SMD} and  $|\cdot|_{\sigma,l}$ is the norm in~\eqref{normMD}. The functions $f: \T^n \times \R^m \times  J \longrightarrow \R^{n+m}$, $z: \T^n \times \R^m \times  J  \longrightarrow \R^{n+m}$ and $g : \T^n \times \R^m \times  J \longrightarrow M_{n+m}$ are given, with $M_{n+m}$ representing the set of $(n+m)$-dimensional matrices.

We obtain a solution $\varkappa$ of the latter thanks to a suitable change of coordinates (see~\eqref{hHE}) that rectifies the dynamic on the cylinder. The most complicated part consists of verifying the regularity and decay in time of $\varkappa$, i.e., the proof that $\varkappa \in \mathcal{V}_{\sigma-1, \omega, 1}$. 
For this purpose,  this section is divided into two parts (Section \ref{SE} and Section \ref{SHE}).  The first part contains several estimates (Lemma \ref{psi} and Lemma \ref{R}).  In the second part, we prove the existence of a solution $\varkappa$ of~\eqref{HEMD} such that $\varkappa \in \mathcal{V}_{\sigma-1, \omega, 1}$ (we refer to Lemma \ref{homoeqlemmaMD}).

\subsection{Several estimates}\label{SE}
We introduce some fundamental Gronwall-type inequalities that we widely use in this section.
\begin{proposition}
\label{Gronwall}
Let $I$ be an interval in $\R$, $t_0 \in I$, and $a$, $b$, $u \in  C(I)$ continuous positive functions. If we assume that 
\begin{equation*}
u(t) \le a(t) + \left|\int_{t_0}^t b(s) u(s) ds \right|, \quad \forall t \in I
\end{equation*}
then it follows that 
\begin{equation}
\label{G1}
u(t) \le a(t) + \left|\int_{t_0}^t a(s)b(s) e^{\left|\int_s^t b(\tau)d\tau \right|} ds \right|, \quad \forall t \in I.
\end{equation}
If $a$ is a non-decreasing function and we assume that 
\begin{equation*}
u(t) \le a(t) + \int_{t_0}^t b(s) u(s) ds  \quad \forall t \ge t_0,
\end{equation*}
then, we obtain the estimate
\begin{equation}
\label{G2}
u(t) \le a(t) e^{\int_{t_0}^t b(s) ds}, \quad \forall t\ge t_0.
\end{equation}
\end{proposition}
\begin{proof}
We refer to~\cite{Amann} for the proof.
\end{proof}

Let $\sigma \ge 1$ be as in~\eqref{HEMD}.  We define $\psi_{t_0, F}^t$ as the flow at time $t$ with initial time $t_0$ of $F$ (where $F=\bar \omega + f$, see~\eqref{HEMD}).  We recall that $C(\cdot)$ stands for constants depending on $n+m$ and the parameters in brackets, while $C$ denotes constants depending on $n+m$.  In this section, we will widely use the properties of the Hölder norm contained in Proposition \ref{Holder}.

\begin{lemma}
\label{psi}
For all $t$, $t_0 \in J$ 
\begin{equation}
\label{estpsi0}
|\partial_q \psi_{t_0, F}^t|_{C^0} \le e^{c_1 \mu |\ln t - \ln t_0| }.
\end{equation}
Moreover, for all $1<s \le \sigma +1$ and $t$, $t_0 \in J$ 
\begin{equation}
\label{estimatepsi1}
|\partial_q \psi^t_{t_0, F}|_{C^{s-1}} \le C(s)\left(1 + |f|_{s, 1} \left|\ln t - \ln t_0\right| \right) e^{c_s \mu |\ln t - \ln t_0|}
\end{equation}
where $c_1 \ge 1$ is a positive constant depending on $n+m$, and $c_s\ge 1$ are constants depending on $n+m$ and $s$ satisfying
\begin{equation}
\label{csc1s}
c_s \ge c_1 s
\end{equation} 
for all $1<s \le \sigma +1$.
\end{lemma}
\begin{proof}
We prove this lemma for $t \ge t_0$.  Similarly, one can prove the other case.  For all $q \in \T^n \times \R^m$, using the fundamental theorem of calculus, we can write $\psi_{t_0, F}^t$ in the following form
\begin{equation*}
\psi_{t_0, F}^t(q) = q + \int_{t_0}^t F^\tau \circ \psi^\tau_{t_0, F}(q)d\tau.
\end{equation*}
Taking the derivative with respect to $q$,  we have that
\begin{equation*}
\partial_q \psi_{t_0, F}^t(q) = \mathrm{Id} + \int_{t_0}^t \partial_q \left(f^\tau \circ \psi_{t_0, F}^\tau(q) \right)d\tau,
\end{equation*}
where $\mathrm{Id}$ stands for the identity matrix and we recall that $F=\bar \omega + f$ (see~\eqref{HEMD}).
Then, utilizing property \textit{1} of Proposition \ref{Holder}, the norm $|\partial_q \psi_{t_0, F}^t|_{C^{s -1}}$ can be estimated as follows
\begin{equation}
\label{basepsi}
|\partial_q \psi_{t_0, F}^t|_{C^{s -1}}  \le 1 + C\int_{t_0}^t \left|f^\tau \circ \psi_{t_0, F}^\tau\right|_{C^s}d\tau
\end{equation}
for all $1 \le s \le \sigma +1$. We will now consider the cases $s = 1$ and $1 < s \le \sigma + 1$ separately.  In both,~\eqref {basepsi} is our starting point. 

\vspace{5mm}
\textit{Case $s = 1$}. By~\eqref{basepsi} and property \textit{3} of Proposition \ref{Holder}, we have that 
\begin{equation*}
|\partial_q \psi^t_{t_0, F}|_{C^0} \le 1+ C\int_{t_0}^t |f^\tau|_{C^0}d\tau +C\int_{t_0}^t |f^\tau|_{C^1}|\partial_q \psi^\tau_{t_0, F}|_{C^0}d\tau.
\end{equation*}
We can estimate the first integral of the latter in the following way
\begin{equation}
\label{EstimateIntf}
 \int_{t_0}^t |f^\tau|_{C^0}d\tau \le \mu\ \int_{t_0}^t {1 \over \tau} d\tau = \ln\left({t \over t_0} \right)^\mu,
\end{equation}
where we recall that $|f|_{1,1} \le \mu$. This implies that
\begin{equation*}
|\partial_q \psi^t_{t_0, F}|_{C^0} \le 1 + \ln \left({t \over t_0} \right)^{C\mu} +C\int_{t_0}^t |f^\tau|_{C^1}|\partial_q\psi^\tau_{t_0, F}|_{C^0}d\tau.
\end{equation*}
Now,  by Gronwall's inequality~\eqref{G2} and $|f|_{1,1} \le \mu$
\begin{eqnarray}
|\partial_q \psi^t_{t_0, F}|_{C^0} &\le&\left(1 + \ln \left({t \over t_0} \right)^{C\mu} \right)e^{C\int_{t_0}^t |f^\tau|_{C^1}d\tau} \nonumber\\
\label{CasePsi0}
&\le& \left(1 + \ln \left({t \over t_0} \right)^{C\mu}\right) e^{\ln \left({t \over t_0} \right)^{C\mu}} \le \left({t \over t_0} \right)^{ c_1 \mu}
\end{eqnarray}
for a suitable constant $ c_1 \ge 1$ depending on $n+m$. This concludes the proof of~\eqref{estpsi0} when $t \ge t_0$. 

\vspace{5mm}
\textit{Case $1 < s \le \sigma + 1$}. Similarly to the previous case, by property \textit{5} of Proposition \ref{Holder} we can estimate the right hand side of~\eqref{basepsi} as follows
\begin{eqnarray*}
|\partial_q \psi^t_{t_0, F}|_{C^{s -1}} &\le& 1 + C(s)\underbrace{\int_{t_0}^t |f^\tau|_{C^0}d\tau}_{I_1} +C(s) \underbrace{\int_{t_0}^t |f^\tau|_{C^s}|\partial_q \psi_{t_0, F}^\tau|^s_{C^0}d\tau}_{I_2}\\
&+&C(s)\int_{t_0}^t |f^\tau|_{C^1}|\partial_q \psi^\tau_{t_0, F}|_{C^{s -1}}d\tau.
\end{eqnarray*}
We have to estimate the integrals $I_1$ and $I_2$ on the right-hand side of the latter. We have already analyzed $I_1$ (see~\eqref{EstimateIntf}),  it remains to consider the other
\begin{eqnarray}
\label{I21F}
 I_2 &=& \int_{t_0}^t |f^\tau|_{C^s}|\partial_q \psi_{t_0, F}^\tau|^s_{C^0}d\tau \le \int_{t_0}^t {|f|_{s,1} \over \tau}\left({\tau \over t_0} \right)^{c_1s\mu}d\tau \\
 \label{I22F}
 &\le& |f|_{s,1} \left({t \over t_0}\right)^{c_1s \mu}\int_{t_0}^t \tau^{ -1}d\tau= |f|_{s,1}\ln\left({t \over t_0}\right)\left({t \over t_0}\right)^{c_1s \mu}.
\end{eqnarray}
We point out that the inequality~\eqref{I21F} is due to~\eqref{CasePsi0}, whereas in~\eqref{I22F} instead of calculating the integral, we prefer using the trivial estimate $\left({\tau \over t_0} \right)^{ c_1s\mu} \le \left({t \over t_0} \right)^{ c_1s\mu}$ to avoid a division by $\mu$ since we do not assume it is not zero.
Thus,  by~\eqref{CasePsi0} and the latter,  we can estimate $|\partial_q \psi^t_{t_0, F}|_{C^{s -1}}$ as follows
\begin{eqnarray*}
|\partial_q \psi^t_{t_0, F}|_{C^{s-1}} &\le& 1 + C(s)\ln\left({t \over t_0}\right)^{\mu} + C(s)|f|_{s,1}\ln\left({t \over t_0}\right)\left({t \over t_0}\right)^{c_1s \mu}\\
&+&C(s)\int_{t_0}^t |f^\tau|_{C^1}|\partial_q \psi^\tau_{t_0, F}|_{C^{s-1}}d\tau.
\end{eqnarray*}
Then, thanks to Gronwall's inequality~\eqref{G2} 
\begin{eqnarray*}
|\partial_q \psi^t_{t_0, F}|_{C^{s-1}} &\le& \left(1 + C(s)\ln\left({t \over t_0}\right)^{\mu} + C(s)|f|_{s,1}\ln\left({t \over t_0}\right)\left({t \over t_0}\right)^{c_1s \mu}\right) e^{C(s) \int_{t_0}^t  |f^\tau|_{C^1}d\tau}\\
&\le& \left(1 + C(s)\ln\left({t \over t_0}\right)^{\mu} + C(s)|f|_{s,1}\ln\left({t \over t_0}\right)\left({t \over t_0}\right)^{c_1s \mu}\right) e^{\ln \left({t \over t_0} \right)^{C(s)\mu}}\\
&\le& C(s)\left(\left({t \over t_0}\right)^{\mu} + |f|_{s,1}\ln\left({t \over t_0}\right)\left({t \over t_0}\right)^{c_1s \mu}\right)  \left({t \over t_0} \right)^{C(s)\mu}\\
&\le& C(s)\left(1 + |f|_{s, 1} \ln \left({t \over t_0}\right)\right) \left({t \over t_0} \right)^{c_s \mu},
\end{eqnarray*}
for a suitable constant $c_s \ge c_1 s$.  This concludes the proof of~\eqref{estimatepsi1} when $t \ge t_0$.
\end{proof}

Let $M_{n+m}$ be the set of the $(n+m)$-dimensional matrices.  We consider the matrix $R : \T^n \times \R^m \times J \times J \to M_{n+m}$ with elements $r_{ij}(q, t,\tau)$ for all $1 \le i,j \le n+m$ and $(q,\tau,t) \in \T^n \times \R^m \times J \times J$. In other words,  $R(q, t,\tau) = \{r_{ij}(q, t,\tau) \}_{1 \le i,j \le n+m}$ for all $(q,\tau,t) \in \T^n \times \R^m \times J \times J$.  To obtain a more elegant form, we introduce the following notation.  For fixed $t$, $\tau \in J$, we define $R^t_{\tau}$ as the function
\begin{equation}
\label{Rtaut}
R^t_{\tau} : \T^n \times \R^m \to M_{n+m}, \quad R^t_{\tau} (q) = R(q, t,\tau).
\end{equation}
We will use this notation in the rest of this section. 
Therefore,  in agreement with the notation made above, for all fixed $t$, $\tau \in J$ and positive real parameters $s \ge 0$,  we define the following family of norms
\begin{equation*}
|R^t_\tau|_{C^s} = \max_{1 \le i,j \le n}|r_{ij}(\cdot,t,\tau)|_{C^s}.
\end{equation*}
Now,  let $g$ and $F$ be as in~\eqref{HEMD} and $\psi_{t_0, F}^t$ be the flow at time $t$ with initial time $t_0$ of $F$. 
The final part of this section is dedicated to the analysis of the following system, which plays a very important role in the solution of the homological equation~\eqref{HEMD}
\begin{equation}
\label{Rs}
\begin{cases}
\partial_t R(q, t,\tau) = -g(\psi^t_{t_0, F}(q),t) R(q, t,\tau) \\
R(q,\tau,\tau) = \mathrm{Id}
\end{cases}
\end{equation}

\begin{lemma}
\label{R}
There exists a unique solution $R : \T^n \times \R^m \times J \times J \to M_{n+m}$ of~\eqref{Rs} such that, for all fixed $t$, $\tau \in J$, $R^t_{\tau} \in C^\sigma(\T^n \times \R^m)$ and $\partial_q^iR \in C( \T^n \times \R^m \times J \times J)$ for all $i \in \N^{n+m}$ with $0 \le |i| \le [\sigma]$.  Moreover,  for all $\tau$,  $t_0$, $t \in J$ with $\tau \ge t$, we define $\tilde R(q,t,\tau) = R(\psi_{\tau, F}^{t_0}(q),t,\tau)$.  For all $1<s \le \sigma$, we have the following estimates
\begin{eqnarray}
\label{stimeR0}
|R^t_\tau|_{C^0} &\le& \left(\tau \over t\right)^{c^R_0\mu}\\
\label{stimeR1}
|\tilde R^t_\tau|_{C^1} &\le& C\left(\tau \over t\right)^{c^R_1\mu}\\
\label{stimeR}
|\tilde R^t_\tau|_{C^s} &\le& C(s )\left(1 + \left(|f|_{s,1} + |g|_{s,1} \right) \ln \left(\tau \over t\right)\right)\left({\tau \over t}\right)^{c^R_s \mu}
\end{eqnarray}
where $c^R_0$ and $c^R_1$ are positive constants depending on $n+m$ and $c^R_s$ is a positive constant depending on $n+m$ and $s$ satisfying $c_1^R \ge c_1$ and $c_s^R \ge c_s$.  We recall that $\sigma$ is the positive parameter in~\eqref{HEMD}, and $c_1$ and $c_s$ are the positive constants introduced in Lemma \ref{psi}.
\end{lemma}
\begin{proof}
For all $q \in \T^n \times \R^m$, by the theorem of existence and uniqueness, a unique solution $R : \T^n \times \R^m \times J \times J \to M_{n+m}$ of~\eqref{Rs} exists. Moreover, for all fixed $t$, $\tau \in J$, $R^t_{\tau} \in C^\sigma(\T^n \times \R^m)$ and $\partial_q^iR \in C( \T^n \times \R^m \times J \times J)$ for all $i \in \N^{n+m}$ with $0 \le |i| \le [\sigma]$.
It remains to prove the estimates~\eqref{stimeR0},~\eqref{stimeR1} and~\eqref{stimeR}.  First, we verify~\eqref{stimeR0}.

Using the fundamental theorem of calculus, we can write $R$ in the following form
\begin{equation}
\label{formaR}
R^t_\tau(q) = \mathrm{Id} + \int_t^\tau \left(g^\delta \circ \psi_{t_0}^\delta (q) \right)R^\delta _\tau(q) d\delta 
\end{equation}
for all $q \in \T^n \times \R^m$ and $t$, $t_0$,  $\tau \in J$ with $\tau \ge t$.  Using Proposition \ref{Holder} and the latter, we can estimate $|R^t_\tau|_{C^0}$ as follows
\begin{equation*}
|R^t_\tau|_{C^0}\le 1 + \int_t^\tau |g^\delta  \circ \psi_{t_0}^\delta  R^\delta_\tau|_{C^0}d\delta  \le 1 + C\int_t^\tau |g^\delta |_{C^0} | R^\delta_\tau|_{C^0}d\delta,
\end{equation*} 
for a suitable constant $C$ depending on $n+m$.
Thanks to Gronwall's inequality~\eqref{G2} and remembering that  $|g|_{1,1} \le \mu$,
\begin{equation*}
|R^t_\tau|_{C^0}\le e^{C\mu \int_t^\tau {1 \over \delta }d\delta } \le \left({\tau \over t}\right)^{c_0^R \mu}
\end{equation*}
for a suitable positive constant $c_0^R$ depending on $n+m$. Hence,~\eqref{stimeR0} is proved.  Now, we verify~\eqref{stimeR1} and~\eqref{stimeR}. To this end, we note that the composition of $R^t_\tau$ with $\psi_\tau^{t_0}(q)$ takes the form
\begin{equation}
\label{formatildeR}
R^t_\tau\circ \psi_\tau^{t_0}(q) = \tilde R^t_\tau (q) = \mathrm{Id} + \int_t^\tau \left(g^\delta \circ \psi_{\tau}^\delta(q)\right) \tilde R^\delta_\tau(q) d\delta
\end{equation}
for all $q \in \T^n \times \R^m$ and $t$, $\tau \in J$ with $\tau \ge t$.  For all $1 \le s\le \sigma$, thanks to the latter, we can estimate the norm $|\tilde R^t_\tau|_{C^s}$ as follows
\begin{equation}
\label{partialestRMD}
|\tilde R^t_\tau|_{C^s}\le 1 + \int_t^\tau |\left(g^\delta \circ \psi_\tau^\delta\right) \tilde R^\delta_\tau|_{C^s}d\delta.
\end{equation}
First of all, we estimate the norm into the above integral. To this end, we use the properties in Proposition \ref{Holder}. 
\begin{eqnarray*}
|\left(g^\delta \circ \psi_\tau^\delta \right) \tilde R^\delta_\tau|_{C^s} &\le& C(s) \left(|g^\delta \circ \psi_\tau^\delta|_{C^s}|R^\delta_\tau|_{C^0} + |g^\delta \circ \psi_\tau^\delta|_{C^0}|\tilde R^\delta_\tau|_{C^s}\right)\\
|g^\delta \circ \psi_\tau^\delta|_{C^s} &\le& C(s)\left(|g^\delta|_{C^s}|\partial_q \psi_\tau^\delta|^s_{C^0} + |g^\delta|_{C^1}|\partial_q \psi_\tau^\delta|_{C^{s-1}} + |g^\delta|_{C^0} \right).
\end{eqnarray*}
Substituting these estimates into \eqref{partialestRMD}, we obtain the following upper bound for $|\tilde R^t_\tau|_{C^s}$
\begin{eqnarray}
|\tilde R^t_\tau|_{C^s} &\le& 1 + C(s)\underbrace{\int_t^\tau |g^\delta|_{C^s} |\partial_q \psi_\tau^\delta|^s_{C^0} |R^\delta_\tau|_{C^0}d\delta}_{I_1^R} + C(s)\underbrace{\int_t^\tau |g^\delta|_{C^1} |\partial_q \psi_\tau^\delta|_{C^{s-1}} |R^\delta_\tau|_{C^0}d\delta}_{I_2^R}\nonumber\\
\label{stimaRIntermedia}
&+& C(s)\underbrace{\int_t^\tau |g^\delta|_{C^0} |R^\delta_\tau|_{C^0}d\delta}_{I_3^R} + C(s)\int_t^\tau |g^\delta|_{C^0} |\tilde R^\delta_\tau|_{C^s}d\delta.
\end{eqnarray}
Now, thanks to~\eqref{stimeR0} and Lemma \ref{psi}, we can estimate the first three integrals $I_1^R$, $I_2^R$ and $I_3^R$ on the right-hand side of the previous inequality
\begin{eqnarray}
I^R_1 &=& \int_t^\tau |g^\delta|_{C^s} |\partial_q \psi_\tau^\delta|^s_{C^0} |R^\delta_\tau|_{C^0}d\delta\le C(s) |g|_{s, 1}\int_t^\tau \delta^{-1}\left({\tau \over \delta}\right)^{c_1s\mu}\left({\tau \over \delta}\right)^{c_0^R\mu}d\delta\nonumber\\
\label{I11}
&\le& C(s)|g|_{s, 1}\left({\tau \over t}\right)^{\left(c_1s + c_0^R\right)\mu}\int_t^\tau \delta^{-1}d\delta\\
&=&C(s)|g|_{s, 1}\ln \left({\tau \over t}\right) \left({\tau \over t}\right)^{\left(c_1s + c_0^R\right)\mu}\nonumber\\
I^R_2 &=& \int_t^\tau |g^\delta|_{C^1} |\partial_q \psi_\tau^\delta|_{C^{s-1}} |R^\delta_\tau|_{C^0}d\delta\nonumber\\
&\le& C(s) \int_t^\tau {|g|_{1,1} \over \delta}\left(1 + |f|_{s, 1} \ln \left({\tau \over \delta}\right)\right)\left({\tau \over \delta} \right)^{\left(c_s + c^R_0\right)\mu}d\delta\nonumber\\
&\le& C(s)\mu \int_t^\tau {1 \over \delta} \left({\tau \over \delta} \right)^{\left(c_s + c^R_0\right)\mu}d\delta\nonumber\\
\label{I21}
&+& C(s)|f|_{s, 1}  \mu \ln \left({\tau \over t}\right)\int_t^\tau {1 \over \delta}\left({\tau \over \delta} \right)^{\left(c_s + c^R_0\right)\mu}d\delta\\
\label{I22}
&\le& C(s)\mu \left({\tau \over t}\right)^{\left(c_s + c_0^R\right) \mu}\int_t^\tau {d\delta \over \delta}\\
&+& C(s)|f|_{s, 1}  \ln \left({\tau \over t}\right)^\mu\left({\tau \over t} \right)^{\left(c_s + c^R_0\right)\mu}\int_t^\tau {d \delta \over \delta}\nonumber\\
&\le& C(s) \left({\tau \over t}\right)^{\left(c_s + c_0^R\right) \mu}\ln \left({\tau \over t}\right)^\mu + C(s)|f|_{s, 1}  \ln \left({\tau \over t}\right)\left({\tau \over t} \right)^{\left(c_s + c^R_0+1\right)\mu}\nonumber\\
&\le& C(s) \left({\tau \over t}\right)^{\left(c_s + c_0^R+1\right) \mu} + C(s)|f|_{s, 1}  \ln \left({\tau \over t}\right)\left({\tau \over t} \right)^{\left(c_s + c^R_0+1\right)\mu}\nonumber\\
\label{I31}
I^R_3 &=& \int_t^\tau |g^\delta|_{C^0} |R^\delta_\tau|_{C^0}d\delta \le \mu \int_t^\tau {1 \over\delta} \left({\tau \over \delta} \right)^{c_0^R\mu}d\delta \le \mu  \left({\tau \over t} \right)^{c_0^R\mu}\int_t^\tau {d \delta \over\delta}\\
&\le&  \left({\tau \over t} \right)^{c_0^R\mu} \ln\left({\tau \over t} \right)^\mu \le \left({\tau \over t} \right)^{\left(c_0^R+1\right)\mu}. \nonumber
\end{eqnarray}
In the aforementioned estimates, we widely used the simple inequality $\ln \left({\tau \over \delta}\right) \le \ln \left({\tau \over t}\right)$ and the following properties $\mu\ln \left({\tau \over t}\right) = \ln \left({\tau \over t}\right)^\mu \le \left({\tau \over t}\right)^\mu$.  Additionally, in lines~\eqref{I11},~\eqref{I22},~\eqref{I31} we employed the trivial estimates $\left({\tau \over \delta}\right)^{\left( c_1s + c_0^R\right)\mu} \le \left({\tau \over t}\right)^{\left( c_1s + c_0^R\right)\mu}$,  $\left({\tau \over \delta}\right)^{\left( c_s + c_0^R\right)\mu} \le \left({\tau \over t}\right)^{\left( c_s + c_0^R\right)\mu}$,  and $\left({\tau \over \delta}\right)^{c_0^R\mu} \le \left({\tau \over t}\right)^{c_0^R\mu}$, respectively. 

Now,  we replace the latter estimates for $I_1^R$, $I_2^R$,  and $I_3^R$ in~\eqref{stimaRIntermedia} and we use that $c_s \ge c_1 s$ for all $1 < s\le \sigma$ (we refer to~\eqref{csc1s} and we observe that when $s=1$, the inequality $c_s \ge c_1 s$ is trivially verified). Then,  for all $1 \le s \le \sigma$
\begin{eqnarray}
|\tilde R^t_\tau|_{C^s} &\le& 1 + C(s)|g|_{s, 1}\ln \left({\tau \over t}\right) \left({\tau \over t}\right)^{\left(c_1s + c_0^R\right)\mu} + C(s) \left({\tau \over t}\right)^{\left(c_s + c_0^R+1\right) \mu}\nonumber\\
&+& C(s)|f|_{s, 1}  \ln \left({\tau \over t}\right)\left({\tau \over t} \right)^{\left(c_s + c^R_0+1\right)\mu} + C(s)\left({\tau \over t} \right)^{\left(c_0^R+1\right)\mu}.\nonumber\\
&+& C(s)\int_t^\tau |g^\delta|_{C^0} |\tilde R^\delta_\tau|_{C^s}d\delta\nonumber\\
&\le& C(s) \left(1 + \left(|f|_{s, 1} + |g|_{s, 1}\right) \ln \left({\tau \over t}\right) \right)\left( {\tau \over t}\right)^{\left(c_s  + c^R_0 + 1\right)\mu}\nonumber\\
\label{RQuasiFinale}
&+& \left|C(s)\int_\tau^t |g^\delta|_{C^0} |\tilde R^\delta_\tau|_{C^s}d\delta\right|.
\end{eqnarray}
In order to rewrite the above estimate in a more compact form,  we introduce the following function
\begin{equation}
\label{aR}
a_{s, \tau}(t) = C(s) \left(1 + \left(|f|_{s, 1} + |g|_{s, 1}\right) \ln \left({\tau \over t}\right) \right)\left( {\tau \over t}\right)^{\left(c_s  + c^R_0 + 1\right)\mu}
\end{equation}
for all $t \in J$ with $\tau \ge t$. Using the latter, we can rewrite~\eqref{RQuasiFinale} in the following form
\begin{equation}
\label{RPrimaGronwall}
|\tilde R^t_\tau|_{C^s} \le a_{s, \tau}(t) +  \left|C(s)\int_\tau^t |g^\delta|_{C^0} |\tilde R^\delta_\tau|_{C^s}d\delta\right|. 
\end{equation}
We observe that $a$ is a decreasing function for all $t \in J$ with $\tau \ge t$ and hence thanks to the Gronwall inequality~\eqref{G1} we obtain that for all $1 \le s \le \sigma$
\begin{eqnarray}
|\tilde R^t_\tau|_{C^s} &\le& a_{s, \tau}(t) +C(s)\left|\int_\tau^t a_{s, \tau}(\delta) |g^\delta|_{C^0} e^{\left|C(s)\int^t_\delta |g^\rho|_{C^0}d\rho\right|}d\delta\right|\nonumber\\\
&\le& a_{s, \tau}(t) + C(s)a_{s, \tau}(t)  \int_t^\tau {\mu \over \delta} e^{\ln \left({\delta \over t}\right)^{C(s) \mu}}d\delta\nonumber\\ &=& a(t)\left(1 + C(s)  \int_t^\tau {\mu \over \delta}  \left({\delta \over t}\right)^{C(s) \mu} d\delta\right)\nonumber\\\
&\le& a_{s, \tau}(t)\left(1 + C(s)   \left({\tau \over t}\right)^{C(s) \mu} \int_t^\tau {\mu \over \delta} d\delta\right)\nonumber\\
&=&a_{s, \tau}(t)\left(1 + C(s)   \left({\tau \over t}\right)^{C(s) \mu} \ln\left({\tau \over t}\right)^\mu\right)\nonumber\\
\label{FinalR}
&\le& C(s) a_{s, \tau}(t) \left({\tau \over t}\right)^{C(s) \mu}.
\end{eqnarray}
Now, we will study the cases $s=1$ and $1 < s \le \sigma$ separately.  In both~\eqref{FinalR} will be our starting point.  If $s=1$, we have that 
\begin{eqnarray*}
|\tilde R^t_\tau|_{C^1} &\le& C a_{1, \tau}(t) \left({\tau \over t}\right)^{C \mu}\\
&=& C \left(1 + \left(|f|_{1, 1} + |g|_{1, 1}\right) \ln \left({\tau \over t}\right) \right)\left( {\tau \over t}\right)^{\left(c_1  + c^R_0 + 1\right)\mu}\left({\tau \over t}\right)^{C \mu}\\
&\le&C \left(1 + \ln \left({\tau \over t}\right)^\mu \right)\left( {\tau \over t}\right)^{\left(c_1  + c^R_0 + 1\right)\mu}\left({\tau \over t}\right)^{C \mu}\le C \left({\tau \over t}\right)^{c_1^R \mu}
\end{eqnarray*}
for a suitable constant $c_1^R \ge c_1+c_0^R +1$.  We point out that in the last line of the latter, we used the simple estimate $\ln \left({\tau \over t}\right)^\mu \le \left({\tau \over t}\right)^\mu$.  This concludes the proof of~\eqref{stimeR1}.  It remains to verify~\eqref{stimeR}.  By~\eqref{FinalR} and~\eqref{aR}, for all $1<s\le \sigma$, we have that 
\begin{eqnarray*}
|\tilde R^t_\tau|_{C^s} &\le& C(s) a_{s, \tau}(t) \left({\tau \over t}\right)^{C(s) \mu}\\
&=& C(s) \left(1 + \left(|f|_{s, 1} + |g|_{s, 1}\right) \ln \left({\tau \over t}\right) \right)\left( {\tau \over t}\right)^{\left(c_s  + c^R_0 + 1\right)\mu}\left({\tau \over t}\right)^{C(s) \mu}\\
&\le& C(s) \left(1 + \left(|f|_{s, 1} + |g|_{s, 1}\right) \ln \left({\tau \over t}\right) \right)\left( {\tau \over t}\right)^{c^R_s \mu}
\end{eqnarray*}
for a suitable constant $c_s^R \ge c_s  + c^R_0 + 1$. This proves~\eqref{stimeR} and concludes the proof of this lemma. 
\end{proof}

\subsection{Solution of the homological equation}\label{SHE}

The aim of this section is to solve the equation~\eqref{HEMD}.  To simplify the reading, let us recall this problem.  We consider $\sigma \ge 1$, $\mu \ge 0$,  $\omega \in \R^n$ and the following functions $f: \T^n \times \R^m \times  J \longrightarrow \R^{n+m}$, $z: \T^n \times \R^m \times  J  \longrightarrow \R^{n+m}$ and $g : \T^n \times \R^m \times  J \longrightarrow M_{n+m}$ such that $f, g\in \mathcal{S}_{\sigma,1}$,  and $z\in \mathcal{S}_{\sigma,2}$ with $|f|_{1,1} \le \mu$, and $|g|_{1,1} \le \mu$.  We are looking for a solution  $\varkappa : \T^n \times \R^m \times J \longrightarrow \R^{n + m}$ of the following equation 
\begin{equation*}
\partial_q\varkappa(q,t) \underbrace{\left(\bar \omega + f(q,t)\right)}_{F(q,t)} + \partial_t \varkappa(q,t) + g(q,t)\varkappa(q,t) = z(q,t)
\end{equation*}
in such a way that $\varkappa \in \mathcal{V}_{\sigma-1, \omega, 1}$.  We recall that $\bar \omega = (\omega, 0) \in \R^{n+m}$ (see~\eqref{BarOmega}).  The following lemma provides a solution to the above problem. We conclude this section with a comparison between the equation~\eqref{HEMD} and the homological equation studied in our previous work~\cite{Sca22}.
\begin{lemma}
\label{homoeqlemmaMD} 
There exists a solution $\varkappa : \T^n \times \R^m \times J \longrightarrow \R^{n + m}$ of equation~\eqref{HEMD}. Moreover, letting  $c_\sigma^\varkappa = \max\{c_0^R + c_\sigma, c_\sigma^R + c_1 \sigma,  c_1^R + c_\sigma\}$, if
\begin{equation}
\label{mu}
\mu < {1 \over c^\varkappa_\sigma}
\end{equation}
then $\varkappa \in \mathcal{V}_{\sigma-1, \omega, 1}$ and 
\begin{equation}
\label{varkappa}
|\varkappa|_{\sigma,1} \le C(\sigma){|z|_{\sigma, 2} \over 1-c^\varkappa_\sigma \mu} + C(\sigma){|f|_{\sigma,1}  + |g|_{\sigma,1} \over \left(1-c^\varkappa_\sigma \mu\right)^2}|z|_{1,2}.
\end{equation}
We recall that the constants $c_0^R$,  $c_1^R$, $c_\sigma^R$, $c_1$ and $c_\sigma$ are defined in Lemma \ref{psi} and Lemma \ref{R}.
\end{lemma}
\begin{proof}
The proof of this lemma is divided into two parts. First, we show the existence of a formal solution to equation~\eqref{HEMD}. In the second part, we verify that, under the smallness assumption~\eqref{mu}, the found solution satisfies a certain regularity and a suitable decay property (see~\eqref{varkappa}).

\textit{Existence of a formal solution}: For fixed $t_0 \in J$, let us define the following transformation 
\begin{align}
\label{hHE}
&h:\T^n \times \R^m \times J \longrightarrow \T^n \times \R^m \times J\\
&h(q,t) = (\psi_{t, F}^{t_0}(q) , t)\nonumber
\end{align}
 where $\psi_{t_0, F}^t$ is the flow at time $t$ with initial time $t_0$ of $F$ previously defined. We claim that it is enough to prove the first part of this lemma for the much simpler equation in the unknown $\kappa : \T^n \times \R^m \times J \to \R^{n+m}$
\begin{equation}
\label{HE2MD}
\partial_t \kappa(q,t) + g \circ h^{-1}(q,t) \kappa(q,t) = z \circ h^{-1}(q,t).
\end{equation} 
Indeed, if $\kappa$ is a solution of the latter, then $\varkappa = \kappa \circ h$ is a solution of~\eqref{HEMD} and vice versa. For the sake of clarity, we prove this claim. Let $\varkappa$ be a solution of~\eqref{HEMD}, then
\begin{small}
\begin{eqnarray*}
\partial_t(\varkappa \circ h^{-1}) + \left(g\circ h^{-1} \right)\left( \varkappa \circ h^{-1}\right) &=& \left( \partial_q\varkappa \circ h^{-1}\right)\partial_t \psi^t_{t_0, F} +  \partial_t\varkappa \circ h^{-1}\\
&+& \left(g\circ h^{-1}\right)\left(\varkappa \circ h^{-1}\right)\\
&=& \left(\partial_q\varkappa \circ h^{-1}\right)\left(F \circ h^{-1}\right)+  \partial_t\varkappa \circ h^{-1} \\
&+& \left(g\circ h^{-1}\right) \left(\varkappa \circ h^{-1}\right)\\
&=& \left(\left(\partial_q \varkappa \right) F+ \partial_t \varkappa +g \varkappa \right)\circ h^{-1} = z \circ h^{-1}.
\end{eqnarray*}
\end{small}
In the latter, we used the following equalities $\partial_t \psi^t_{t_0, F}(q) = F(\psi^t_{t_0, F}(q), t) = F \circ h^{-1}(q,t)$ for all $(q,t) \in \T^n \times \R^m \times J$.  Moreover, in the last line, we used that $\varkappa$ is a solution of~\eqref{HEMD}.  We proved that $\varkappa \circ h^{-1}$ is a solution for~\eqref{HE2MD}.  It remains to verify the other implication. 
But first, we need to show that
\begin{equation}
\label{KeyEqExClaim}
\partial_q \psi_{t, F}^{t_0}(q) F(q,t)+ \partial_t\psi_{t, F}^{t_0}(q)=0
\end{equation}
for all $q \in \T^n \times \R^m$ and $t$, $t_0 \in J$.  To this end,  we want to rewrite $\partial_q \psi_{t, F}^{t_0}$ and $\partial_t\psi_{t, F}^{t_0}$ in a more convenient form (see~\eqref{HEBDpsiq} and~\eqref{HEBDpsit}, below).
We consider the following trivial equality
\begin{equation}
\label{relazionepsiNC}
\psi^t_{t_0, F}\circ \psi_{t, F}^{t_0}(q) = q,
\end{equation}
for all $q \in \T^n \times \R^m$ and $t$, $t_0 \in J$.
Differentiating both sides of the latter with respect to the variable $q \in \T^n \times \R^m$, we obtain
\begin{equation*}
\partial_q \psi^t_{t_0, F}\circ \psi_{t, F}^{t_0}(q)\partial_q \psi_{t, F}^{t_0}(q) = \mathrm{Id}
\end{equation*}
for all $q \in \T^n \times \R^m$ and $t$, $t_0 \in J$, where $\mathrm{Id}$ stands for the identity matrix.  By the above equation, one can verify that 
\begin{equation}
\label{HEBDpsiq}
\partial_q \psi_t^{t_0}(q) = \left(\partial_q \psi^t_{t_0, F}\circ \psi_t^{t_0}(q)\right)^{-1}.
\end{equation}
On the other hand, taking the derivative with respect to $t$ on both sides of~\eqref{relazionepsiNC}, we can see that
\begin{equation*}
0 = {d \over dt}\left(\psi^t_{t_0, F}\circ \psi_{t, F}^{t_0}(q) \right) =  \partial_q \psi^t_{t_0, F}\circ \psi_{t, F}^{t_0}(q) \partial_t \psi^{t_0}_{t, F}(q) + \partial_t \psi_{t_0, F}^t\circ  \psi^{t_0}_{t, F}(q)
\end{equation*}
for all $q \in \T^n \times \R^m$ and $t$, $t_0 \in J$.  Thanks to the latter, $\partial_t \psi^{t_0}_t(q)$ is equal to
\begin{equation}
\label{HEBDpsit}
\partial_t \psi^{t_0}_{t, F}(q) =  -\left(\partial_q \psi^t_{t_0, F}\circ \psi_{t, F}^{t_0}(q)\right)^{-1} \partial_t \psi_{t_0, F}^t\circ  \psi^{t_0}_{t, F}(q).
\end{equation}
Now, using~\eqref{HEBDpsiq} and~\eqref{HEBDpsit}, we can rewrite $\partial_q \psi_t^{t_0} F+ \partial_t\psi_t^{t_0}$ in the following form 
\begin{equation*}
\partial_q \psi_{t, F}^{t_0}(q) F(q,t)+ \partial_t\psi_{t, F}^{t_0}(q)=\left(\partial_q \psi^t_{t_0, F}\circ \psi_{t, F}^{t_0}(q) \right)^{-1} \left(F(q,t) -  \partial_t \psi_{t_0, F}^t\circ \psi^{t_0}_{t, F}(q)\right)
\end{equation*}
for all $t$, $t_0 \in J$ and $q \in \T^n \times \R^m$. This implies~\eqref{relazionepsiNC} because
\begin{equation*}
F(q,t) -  \partial_t \psi_{t_0, F}^t\circ \psi^{t_0}_{t, F}(q) = F(q,t) - F(\psi_{t_0, F}^t\circ \psi^{t_0}_{t, F}(q), t) = 0
\end{equation*}
for all $t$, $t_0 \in J$ and $q \in \T^n \times \R^m$.  

Now that we have proven~\eqref{relazionepsiNC}, we can use it to verify that if $\kappa$ is a solution of~\eqref{HE2MD}, then $\kappa \circ h$ is a solution of~\eqref{HEMD}.  For this purpose,  let $\kappa$ be a solution of~\eqref{HE2MD}, we have that
\begin{eqnarray*}
\partial_q(\kappa \circ h)  F+ \partial_t(\kappa \circ h) + g ( \kappa \circ h)&=& \left(\partial_q\kappa\circ h \right) \left(\partial_q \psi_t^{t_0} F+ \partial_t\psi_t^{t_0}\right) + \partial_t \kappa \circ h\\
&+& g ( \kappa \circ h)\\
&=&\partial_t \kappa \circ h + g ( \kappa \circ h)= z.
\end{eqnarray*}
In the last equality of the latter, we used~\eqref{relazionepsiNC} and the fact that $\kappa$ is a solution of~\eqref{HE2MD}.  This proves that $\kappa \circ h$ is a solution of~\eqref{HEMD}, and hence the claim stated at the beginning of this proof. 

Now, we can reduce the proof of the first part of this lemma by studying the existence of a solution for the easier equation~\eqref{HE2MD}.  Let $R : \T^n \times \R^m \times J \times J \to M_{n+m}$ be the unique solution of~\eqref{Rs}. Then, a solution $\kappa : \T^n \times \R^m \times J \to \R^{n+m}$ of~\eqref{HE2MD} exists and 
\begin{eqnarray*}
\kappa (q,t) &=& R(q, t, 1)  e(q) - \int_{1}^t R(q, t, \tau)  z\circ h^{-1}(q, \tau) d\tau\\
&=& R(q, t,1)  \left(e(q) - \int_{1}^t R(q, 1, \tau)  z \circ h^{-1}(q, \tau) d\tau  \right)
\end{eqnarray*}
with a free function $e:\T^n \times \R^m \to \R^{n+m}$.  This concludes the proof of the first part of this lemma. 

\textit{Regularity of the found solution}: We choose $e$ in such a way that 
\begin{equation}
\label{e}
e(q) = \int_{1}^{+\infty}R(q, 1, \tau)  z \circ h^{-1}(q, \tau) d\tau.
\end{equation}
It is well defined because by Proposition \ref{Holder}, Lemma \ref{R} and~\eqref{mu},
\begin{eqnarray*}
\left|\int_{1}^{+\infty}R(q, 1, \tau)  z\circ h^{-1}(q, \tau) d\tau \right| &\le& C\int_{1}^{+\infty}|R^{1}_\tau|_{C^0} |z^\tau|_{C^0} d\tau\\
&\le&C \int_{1}^{+\infty}\tau^{c_0^R\mu} {|z|_{0,2} \over \tau^2} d\tau \\
&=& C |z|_{0,2}\int_{1}^{+\infty}\tau^{c_0^R\mu-2} d\tau\\
&=& C{|z|_{0,2}  \over 1-c_0^R\mu}.
\end{eqnarray*}
Now, letting $e$ as in~\eqref{e}, we can rewrite $\varkappa=\kappa \circ h$ in the following form
\begin{eqnarray}
\varkappa(q,t) &=& \kappa \circ h(q,t) = -\int_t^{+\infty} R^t_\tau\circ \psi_t^{t_0}(q)  z^\tau\circ\psi_t^\tau(q)d\tau\nonumber\\
&=&-\int_t^{+\infty} R^t_\tau \circ \psi_\tau^{t_0} \circ \psi_t^\tau(q)  z^\tau\circ\psi_t^\tau(q)d\tau\nonumber\\
\label{varkappa2}
&=& -\int_t^{+\infty} \tilde R^t_\tau \circ \psi_t^\tau(q)  z^\tau\circ\psi_t^\tau(q)d\tau
\end{eqnarray}
for all $(q, t) \in \T^n \times \R^m \times J$.
In what follows, we will verify that, under the assumption~\eqref{mu},  $\varkappa \in \mathcal{V}_{\sigma-1, \omega, 1}$ and it satisfies~\eqref{varkappa}.We begin by proving~\eqref{varkappa}.

For all $t \in J$, by~\eqref{varkappa2} and Proposition \ref{Holder}, we can estimate $|\varkappa^t|_{C^\sigma}$ as follows
\begin{equation}
\label{varkappaPrimaStima}
|\varkappa^t|_{C^\sigma} \le C(\sigma) \int_t^{+\infty} |R^t_\tau|_{C^0}|z^\tau \circ \psi_t^\tau|_{C^\sigma} + |\tilde R^t_\tau \circ \psi_t^\tau |_{C^\sigma}|z^\tau|_{C^0}d\tau.
\end{equation}
Now, we want to estimate the norms into the above integral.  For this purpose, thanks to property \textit{5} of Proposition \ref{Holder}, we have that 
\begin{eqnarray*}
|z^\tau \circ \psi_t^\tau|_{C^\sigma} &\le& C(\sigma)\left( |z^\tau|_{C^\sigma}|\partial_q \psi_t^\tau|^\sigma_{C^0} + |z^\tau|_{C^1}|\partial_q\psi_t^\tau|_{C^{\sigma-1}} + |z^\tau|_{C^0}\right)\\
 |\tilde R^t_\tau \circ \psi_t^\tau |_{C^\sigma} &\le& C(\sigma)\left( |\tilde R^t_\tau|_{C^\sigma}|\partial_q \psi_t^\tau|^\sigma_{C^0} + |\tilde R^t_\tau|_{C^1}|\partial_q\psi_t^\tau|_{C^{\sigma-1}} + | R^t_\tau|_{C^0}\right)
\end{eqnarray*}
and replacing the latter into~\eqref{varkappaPrimaStima}, we obtain the following upper bound for $|\varkappa^t|_{C^\sigma}$
\begin{eqnarray}
|\varkappa^t|_{C^\sigma} &\le& C(\sigma) \int_t^{+\infty} |R^t_\tau|_{C^0}|z^\tau|_{C^\sigma}|\partial_q \psi_t^\tau|^\sigma_{C^0} d\tau + C(\sigma)\int_t^{+\infty} |R^t_\tau|_{C^0}|z^\tau|_{C^1}|\partial_q\psi_t^\tau|_{C^{\sigma-1}}d\tau \nonumber\\
&+&C(\sigma)\int_t^{+\infty} |\tilde R^t_\tau|_{C^\sigma}|\partial_q \psi_t^\tau|^\sigma_{C^0}|z^\tau|_{C^0}d\tau + C(\sigma)\int_t^{+\infty} |\tilde R^t_\tau|_{C^1}|\partial_q \psi_t^\tau|_{C^{\sigma -1}}|z^\tau|_{C^0}d\tau\nonumber\\
\label{varkappaQuasiFinale}
&+& C(\sigma)\int_t^{+\infty} | R^t_\tau|_{C^0} |z^\tau|_{C^0}d\tau.
\end{eqnarray}
Now, we have to estimate each integral on the right-hand side of the latter. First, we observe that, for all $t \in J$ and $0 \le x <1$
\begin{equation*}
\int_t^{+\infty} \tau^{x-2} \ln \left(\tau \over t\right)d\tau = {1 \over 1-x}\int_t^{+\infty} \tau^{x-2} d\tau.
\end{equation*}
It is obtained by integrating by part. Then, using Lemma \ref{psi}, Lemma \ref{R},~\eqref{mu} and the latter
\begin{eqnarray*}
\int_t^{+\infty} |R^t_\tau|_{C^0}|z^\tau|_{C^\sigma}|\partial_q \psi_t^\tau|^\sigma_{C^0} d\tau &\le&C(\sigma) \int_t^{+\infty} {|z|_{\sigma, 2}\over \tau^2}\left({\tau \over t}\right)^{\left(c_0^R + c_1\sigma\right)\mu}d\tau\\
&=& C(\sigma){|z|_{\sigma, 2} \over t^{\left(c^R_0 + c_1\sigma\right)\mu}} \int_t^{+\infty}\tau^{\left(c^R_0 + c_1\sigma \right)\mu-2}d\tau\\
&=&C(\sigma){|z|_{\sigma, 2} \over 1-\left(c^R_0 + c_1\sigma\right)\mu}{1\over t}
\end{eqnarray*}
\begin{eqnarray*}
\int_t^{+\infty} |R^t_\tau|_{C^0}|z^\tau|_{C^1}|\partial_q\psi_t^\tau|_{C^{\sigma-1}}d\tau &\le& C(\sigma) \int_t^{+\infty}{|z|_{1, 2}\over \tau^2}\left(1 + |f|_{\sigma,1} \ln\left(\tau \over t\right) \right)\left({\tau \over t}\right)^{\left(c_0^R + c_\sigma\right)\mu}d\tau\\
&=&C(\sigma)\int_t^{+\infty}{|z|_{1, 2}\over \tau^2}\left({\tau \over t}\right)^{\left(c_0^R + c_\sigma\right)\mu}d\tau\\
&+&C(\sigma)\int_t^{+\infty}{|z|_{1, 2}\over \tau^2}|f|_{\sigma,1} \ln\left(\tau \over t\right)\left({\tau \over t}\right)^{\left(c_0^R + c_\sigma\right)\mu}d\tau\\ 
&=&C(\sigma) {|z|_{1,2} \over 1 - (c_0^R + c_\sigma)\mu}{1\over t}\\
&+& C(\sigma) {|z|_{1,2} |f|_{\sigma,1} \over 1 - (c_0^R + c_\sigma)\mu}{1\over t^{(c_0^R + c_\sigma)\mu}}\int_t^{+\infty}\tau^{\left(c^R_0 +  c_\sigma \right)\mu-2}d\tau\\
&=&C(\sigma) {|z|_{1,2} \over 1 - (c_0^R + c_\sigma)\mu}{1\over t}+C(\sigma) {|z|_{1,2} |f|_{\sigma,1} \over \left(1 - (c_0^R + c_\sigma)\mu\right)^2}{1\over t}\\
\int_t^{+\infty} |R^t_\tau|_{C^0} |z^\tau|_{C^0}d\tau &\le& C \int_t^{+\infty}{|z|_{0,2} \over \tau^2}\left({\tau \over t} \right)^{c_0^R\mu}d\tau = C{|z|_{0,2} \over 1-c_0^R\mu}{1\over t}\\
\int_t^{+\infty} |\tilde R^t_\tau|_{C^\sigma}|\partial_q \psi_t^\tau|^\sigma_{C^0}|z^\tau|_{C^0}d\tau &\le& C(\sigma )\int_t^{+\infty}\left(1 + \left(|f|_{\sigma,1} + |g|_{\sigma,1}\right)\ln \left(\tau \over t\right)\right){|z|_{0,2} \over \tau^2}\left({\tau \over t}\right)^{\left(c^R_\sigma + c_1\sigma\right) \mu}d\tau\\
&=&C(\sigma )\int_t^{+\infty}{|z|_{0,2} \over \tau^2}\left({\tau \over t}\right)^{\left(c^R_\sigma + c_1\sigma\right) \mu}d\tau\\
&+& C(\sigma) \left(|f|_{\sigma,1} + |g|_{\sigma,1}\right)\int_t^{+\infty}\ln \left(\tau \over t\right){|z|_{0,2} \over \tau^2}\left({\tau \over t}\right)^{\left(c^R_\sigma + c_1\sigma\right) \mu}d\tau\\
&=&C(\sigma) {|z|_{0,2} \over 1 - (c_\sigma^R + c_1\sigma)\mu}{1\over t}+C(\sigma) {|z|_{0,2}\left(|f|_{\sigma,1} + |g|_{\sigma,1}\right) \over \left(1 - (c_\sigma^R + c_1\sigma)\mu\right)^2}{1\over t}
\end{eqnarray*}
\begin{eqnarray*}
\int_t^{+\infty} |\tilde R^t_\tau|_{C^1}|\partial_q \psi_t^\tau|_{C^{\sigma -1}}|z^\tau|_{C^0}d\tau &\le&C(\sigma ) \int_t^{+\infty}{|z|_{0, 2}\over \tau^2}\left(1 + |f|_{\sigma,1} \ln\left(\tau \over t\right) \right)\left({\tau \over t}\right)^{\left(c_1^R + c_\sigma\right)\mu}d\tau\\
&=&C(\sigma)\int_t^{+\infty}{|z|_{0, 2}\over \tau^2}\left({\tau \over t}\right)^{\left(c_1^R + c_\sigma\right)\mu}d\tau\\
&+&C(\sigma)\int_t^{+\infty}{|z|_{0, 2}\over \tau^2}|f|_{\sigma,1} \ln\left(\tau \over t\right)\left({\tau \over t}\right)^{\left(c_1^R + c_\sigma\right)\mu}d\tau\\ 
&=&C(\sigma) {|z|_{0,2} \over 1 - \left(c_1^R + c_\sigma\right)\mu}{1\over t}+C(\sigma) {|z|_{0,2} |f|_{\sigma,1} \over \left(1 - \left(c_1^R + c_\sigma\right)\mu\right)^2}{1\over t}.
\end{eqnarray*}
Then, replacing the above estimate into~\eqref{varkappaQuasiFinale}, we obtain that
\begin{eqnarray*}
|\varkappa^t|_{C^\sigma}t &\le& C(\sigma) \Bigg({|z|_{\sigma, 2} \over 1-\left(c^R_0 + c_1\sigma\right)\mu} + {|z|_{1,2} \over 1 - (c_0^R + c_\sigma)\mu} +  {|z|_{1,2} |f|_{\sigma,1} \over \left(1 - (c_0^R + c_\sigma)\mu\right)^2} \\
&+& {|z|_{0,2} \over 1-c_0^R\mu}+{|z|_{0,2} \over 1 - (c_\sigma^R + c_1\sigma)\mu}+{|z|_{0,2}\left(|f|_{\sigma,1} + |g|_{\sigma,1}\right) \over \left(1 - (c_\sigma^R + c_1\sigma)\mu\right)^2}\\
&+&{|z|_{0,2} \over 1 - \left(c_1^R + c_\sigma\right)\mu}+{|z|_{0,2} |f|_{\sigma,1} \over \left(1 - \left(c_1^R + c_\sigma\right)\mu\right)^2}\Bigg)\\
&\le&C(\sigma){|z|_{\sigma, 2} \over 1-c_\sigma^\varkappa\mu} + C(\sigma){\left(|f|_{\sigma,1} + |g|_{\sigma,1}\right) \over \left(1 - c_\sigma^\varkappa\mu\right)^2}|z|_{1,2}
\end{eqnarray*}
for all $t\in J$, where we recall that $c_\sigma^\varkappa = \max\{c_0^R + c_\sigma, c_\sigma^R + c_1 \sigma,  c_1^R + c_\sigma\}$ . Taking the sup for all $t \in J$ on the left-hand side of the latter, we conclude the proof of~\eqref{varkappa}.  It remains to verify that $\varkappa \in \mathcal{V}_{\sigma-1, \omega, 1}$ (we refer to~\eqref{V} for the definition of $\mathcal{V}_{\sigma-1, \omega, 1}$).  Thanks to~\eqref{varkappa},  and the regularity of $F$ and $R$ (see Lemma \ref{R}),  one can verify that $\varkappa \in \mathcal{S}_{\sigma, 1}$. We need to show that $D\varkappa \bar\Omega \in \mathcal{S}_{\sigma-1, 2}$, we refer to Section \ref{FS} for the notation.  We know that $\varkappa$ is a solution of~\eqref{HEMD}, hence $D\varkappa \bar\Omega$ satisfies the following equation
\begin{equation*}
D\varkappa(q,t)\bar \Omega = z(q,t) - \partial_q\varkappa(q,t) f(q,t) -  g(q,t)\varkappa(q,t)
\end{equation*}
for all $(q,t) \in \T^n \times \R^m\times J$. Thanks to the latter, and remembering that $\varkappa \in \mathcal{S}_{\sigma, 1}$, one can verify that $D\varkappa \bar\Omega \in \mathcal{S}_{\sigma-1, 2}$. This concludes the proof of this lemma. 
\end{proof}

The homological equation~\eqref{HEMD} studied in this section is considerably more complicated than the one analyzed in the previous work~\cite{Sca22}.  As previously mentioned, in the present paper, we relax the decay in time of the perturbative terms.  However, this comes at the cost of losing information regarding the dynamics at infinity of the obtained solutions.

\begin{remark}\label{CompHE}
 The homological equation studied in~\cite{Sca22}  is derived from~\eqref{HEMD} by setting $m=0$ and $\mu=0$.  More specifically,  given $\sigma \ge 1$, $l>1$ and $\omega \in \R^n$, we investigated the following problem for the unknown $\varkappa:\T^n \times J \to \R$
\begin{equation}
\label{HEGD}
\begin{cases}
\partial_q\varkappa(q,t) \omega + \partial_t \varkappa(q,t) = z(q,t)\\
z\in \mathcal{S}_{\sigma,l}
\end{cases}
\end{equation}
where $z:\T^n \times J \to \R$ is given.  We solve the above equation by integration thanks to a change of coordinates that rectifies the dynamic on the torus.  Consequently, a unique solution of~\eqref{HEGD} satisfying $\varkappa \in \mathcal{S}_{\sigma, l-1}$ exists and it is given by
\begin{equation}
\varkappa(q,t) = -\int_t^{+\infty} z(q + \omega(\tau- t), \tau)d\tau.
\end{equation}
Unlike the problem~\eqref{HEGD}, to solve~\eqref{HEMD}, we need to impose a smallness assumption on $f$ and $g$. Additionally, the proof of estimate~\eqref{varkappa}, and consequently $\varkappa \in \mathcal{V}_{\sigma-1, \omega, 1}$,  requires significantly more effort compared to our previous work in~\cite{Sca22}. We point out that we do not prove the existence of a unique solution of~\eqref{HEMD}, whereas there exists a unique solution $\varkappa$ of~\eqref{HEGD} such that $\varkappa \in \mathcal{S}_{\sigma, l-1}$.
\end{remark}

We want to emphasize that Lemma \ref{homoeqlemmaMD} does not provide $C^\infty$ solutions of the homological equation~\eqref{HEMD}.  The reason lies in the definition of the constant $c_\sigma^\varkappa$ introduced in Lemma \ref{homoeqlemmaMD}.  More specifically, we recall that
\begin{equation}
\label{HEInfty1}
c_\sigma^\varkappa = \max\{c_0^R + c_\sigma, c_\sigma^R + c_1 \sigma,  c_1^R + c_\sigma\}
\end{equation}
where $c_0^R$,  $c_1^R$, $c_\sigma^R$, $c_1$ and $c_\sigma$ are the constants defined in Lemma \ref{psi} and Lemma \ref{R}. We have established that
\begin{equation}
\label{HEInfty2}
c_\sigma^R \ge c_\sigma \ge c_1 \sigma \ge 1
\end{equation}
as proven in Lemma \ref{R} and~\eqref{csc1s} contained in Lemma \ref{psi}.  Using~\eqref{HEInfty1} and~\eqref{HEInfty2}, one can see that 
\begin{equation*}
\lim_{\sigma \to + \infty}c_\sigma^\varkappa = +\infty.
\end{equation*}
Hence, the bigger $\sigma$ is, the more we have to take $\mu$ small in order to find a solution $\varkappa \in \mathcal{V}_{\sigma-1, \omega, 1}$ of~\eqref{HEMD}. 

In the analytic setting, we are also unable to provide holomorphic solutions for~\eqref{HEMD}.  In this case,  for all fixed $t \in J$, $F^t = \bar \omega + f^t$ is defined in a complex neighborhood of $\T^n \times \R^m$ and $f^t$ decay as ${1 \over t}$.  Thus, the flow $\psi^t_{t_0, F}$ diverges at least as $\ln t$.  This prevents us from well-defining the change of coordinates~\eqref{hHE} and finding a solution to the equation~\eqref{HEMD}.

\section{The Nash-Moser theorem}\label{Zehnder}

In this section, we recall the Nash-Moser theorem proved by Zehnder~\cite{Zeh76}. For this purpose, we consider three one-parameter families of Banach spaces $\{(\mathcal{X}^\sigma, |\cdot|_\sigma)\}_{\sigma \ge 0}$, $\{(\mathcal{Y}^\sigma, |\cdot|_\sigma)\}_{\sigma \ge 0}$ and $\{(\mathcal{Z}^\sigma, |\cdot|_\sigma)\}_{\sigma \ge 0}$ satisfying the property~\eqref{Xsigma}.  Moreover, we assume that each of these families of Banach spaces has a $C^\infty$-smoothing, which we denote by the same letter $\{ S_\tau\}_{\tau >0}$.  Now,  let $\mathcal{F}$ be the following functional 
\begin{equation*}
\mathcal{F} : \mathcal{X}^0 \times \mathcal{Y}^0 \longrightarrow \mathcal{Z}^0
\end{equation*}
such that
\begin{equation*}
 \mathcal{F}(x_0, y_0) =0
\end{equation*}
for some $(x_0, y_0) \in \mathcal{X}^0 \times \mathcal{Y}^0$. Given a positive parameter $0< \zeta \le 1$, for all $\sigma \ge 0$, we define 
\begin{equation}
\label{OMD}
\mathcal{O}_\zeta^\sigma = \{(x,y) \in \mathcal{X}^\sigma \times \mathcal{Y}^\sigma  : |x - x_0|_\sigma, |y - y_0|_\sigma <\zeta \}
\end{equation}
and we assume $\mathcal{F} : \mathcal{O}^0_\zeta  \to \mathcal{Z}^0$ to be continuous. 
For given $x \in \mathcal{X}^0 \cap \mathcal{O}_\zeta^0$, the aim of the Zehnder theorem is to solve the equation $\mathcal{F}(x,y) =0$ assuming $x$ sufficiently close to $x_0$. The author makes the following hypotheses \textit{H1-H4}.

\vspace{5mm}
\textit{H.1 Smoothness:} We assume that $\mathcal{F}(x, \cdot):\mathcal{Y}^0 \to \mathcal{Z}^0$ is two times differentiable with the uniform estimate
\begin{equation*}
|D_y \mathcal{F}(x,y)|_0, |D_y^2 \mathcal{F}(x,y)|_0 \le C
\end{equation*} 
for all $(x,y) \in \mathcal{O}_\zeta^0$ and for some constant $C \ge 1$, where $D_y$ is the differential with respect to the second component.  

\vspace{3mm}
\textit{H.2 $\mathcal{F}$ is uniformly Lipschitz in $\mathcal{X}^0$:} For all $(x_1, y)$, $(x_2, y) \in \mathcal{O}_\zeta^0$, 
\begin{equation*}
|\mathcal{F}(x_1, y) - \mathcal{F}(x_2, y)|_0 \le C|x_1 - x_2|_0.
\end{equation*}

\vspace{3mm}
\textit{H.3 Existence of a right-inverse of loss $\upsilon$, $1 \le \upsilon < s$ ($s$ will be specified later):} For every $(x,y) \in \mathcal{O}_\zeta^\upsilon$ there exists a linear map $\eta(x,y): \mathcal{Z}^\upsilon \to \mathcal{Y}^0$ such that, for all $z \in \mathcal{Z}^\upsilon$,
\begin{align}
& D_y \mathcal{F}(x,y) \circ \eta (x,y)z = z \nonumber\\
\label{eta1}
& |\eta(x,y)z|_0 \le C |z|_\upsilon.
\end{align} 
Moreover, for all $\upsilon \le \sigma \le s$, if $(x,y) \in \mathcal{O}_\zeta^\upsilon \cap  (\mathcal{X}^\sigma \times \mathcal{Y}^\sigma)$, then the linear map $\eta : \mathcal{Z}^\sigma \to \mathcal{Y}^{\sigma -\upsilon}$  is well defined and if $|x - x_0|_\sigma$,  $|y - y_0|_\sigma\le K$, then 
\begin{equation}
\label{eta2}
|\eta(x,y)\mathcal{F}(x,y)|_{\sigma-\upsilon} \le C(\sigma)K.
\end{equation}

\vspace{3mm}
\textit{H.4 Order:} The triple $(\mathcal{F}, x_0, v_0)$ is of order $s$, $s > \upsilon \ge 1$. Here, Zehnder uses the following 
\begin{definition}
\label{order}
$(\mathcal{F}, x_0, y_0)$ is called of order $s$, $1\le s<\infty$, if the following three conditions are satisfies:
\begin{enumerate}
\item $(x_0, y_0) \in \mathcal{X}^s \times \mathcal{Y}^s$,
\item $\mathcal{F}(\mathcal{O}_\zeta^0 \cap (\mathcal{X}^\sigma \times \mathcal{Y}^\sigma)) \subset \mathcal{Z}^\sigma$, \quad $1 \le \sigma\le s$
\item there exist constants $C(\sigma)$, $1 \le \sigma \le s$, such that if $(x,y) \in (\mathcal{X}^\sigma \times \mathcal{Y}^\sigma)\cap \mathcal{O}_\zeta^1$ satisfies $|x - x_0|_\sigma$, $|y - y_0|_\sigma \le K$ then 
\begin{equation*}
|\mathcal{F}(x,y)|_\sigma \le C(\sigma) K.
\end{equation*}
\end{enumerate}
\end{definition}

Zehnder, in his paper, assumes the existence of an approximate right-inverse. The reason for this assumption is evident in his works~\cite{Zeh76} and~\cite{Zeh75}, where he applies generalized implicit function theorems to solve small divisor problems. This includes proving Arnold's normal form theorem for vector fields on the torus and the KAM theorem. 
In the proof of these theorems, Zehnder defines a functional $\mathcal{F}$, which does not admit a right-inverse but just an approximate right-inverse. We do not have this problem; hence, we prefer to write H.4 in this form. In his paper~\cite{Zeh76}, Zehnder proved the following

\begin{theorems}
\label{NMZ}
Let $\alpha$, $\beta$, $\lambda$, $\rho$, $\upsilon$ and $s$ be positive real numbers satisfying the following set of inequalities:
\begin{equation}
\begin{aligned}
\label{T1}
&1  < \beta <2, \quad 1 < \alpha, \quad 1 \le \upsilon \le  \rho < \lambda < s,\\
&\hspace{5mm} \lambda > \max\left\{{2 \beta \upsilon \over 2-\beta}, \beta(\upsilon  + \rho\beta)  \right\}\\
&\hspace{5mm} s > \max\left\{{\alpha \upsilon \over \alpha - 1}, \lambda + {\alpha \upsilon \over \beta - 1} \right\}.
\end{aligned}
\end{equation}
Let $(\mathcal{F}, x_0, v_0)$ be of order $s$ and satisfy H.1-H.4 with a loss of $\upsilon$. Then there exists $\varepsilon_0$, depending on $\alpha$, $\beta$, $\lambda$, $\mu$, $s$ and $\zeta$, such that for all $\varepsilon \le \varepsilon_0$ we have the existence of an open neighborhood $\mathcal{D}^\lambda \subset \mathcal{X}^\lambda$ of $x_0$, $\mathcal{D}^\lambda = \{x \in \mathcal{X}^\lambda : |x-x_0|_\lambda \le \varepsilon \}$ and a mapping $\psi : \mathcal{D}^\lambda \to \mathcal{V}^\rho$ such that
\begin{align*}
&\mathcal{F}(x, \psi(x))=0, \quad x\in \mathcal{D}^\lambda \\
&|\psi(x) - v_0|_\rho \le \zeta,
\end{align*}
where $\zeta$ is the positive parameter defined by~\eqref{OMD}.
\end{theorems}
The statement of the above theorem is slightly different from the original.  In his work, Zehnder~\cite{Zeh76} considers $\zeta=1$.  It is simply a quantitative difference, and it does not change the proof, but in this case, one needs to take $\varepsilon_0$ small also with respect to $\zeta$.  This change does not warrant further justifications. However, one can see~\cite{Sca22T} for more details.

\section{Proof of Theorem \ref{MD}}\label{ProofMD}

\subsection{Outline of the proof of Theorem \ref{MD}}\label{OotPBD}
Let $\rho$ and $\omega$ be as in~\eqref{HMD}. We are looking for a $C^\rho$-weakly asymptotic cylinder $\varphi^t$ associated to $(X_H, X_{h}, \varphi_0)$, where $H$ and $h$ are the Hamiltonians defined by~\eqref{HMD} and~\eqref{tildehBD}, respectively,  and $\varphi_0$ is the trivial embedding $\varphi_0 : \T^n \times \R^m  \to \T^n \times \R^m \times  B$,  such that $\varphi_0(q) = (q,0)$. More concretely, for given $H$, we are searching for $v$, $\Gamma :\T^n \times \R^m \times J \to  \R^{n+m}$  such that
\begin{equation}
\label{VarphiThmA}
\varphi:\T^n \times \R^m \times J \to \T^n \times \R^m \times B \times J, \quad \varphi(q,t) = (q, v(q,t))
\end{equation}
is a family of embedded cylinder and $\varphi$, $v$ and $\Gamma$ satisfy 
\begin{align}
\label{hyp1MD2}
&   X_H(\varphi(q, t), t) -  \partial_q \varphi(q, t)(\bar \omega + \Gamma(q,t))  - \partial_t \varphi(q, t) = 0,\\
\label{hyp2MD2}
&  v, \Gamma \in \mathcal{S}_{\rho, 1}.
\end{align}
We refer to~\eqref{BarOmega} for the notation $\bar \omega$. We recall that if  $\varphi$, $v$ and $\Gamma$ satisfy~\eqref{hyp1MD2} and~\eqref{hyp2MD2}, then $\varphi^t$ is a $C^\rho$-weakly asymptotic cylinder $\varphi^t$ associated to $(X_H, X_{\tilde h}, \varphi_0)$ (one can see Definition \ref{weakasymcyl}).  In fact,~\eqref{hyp1MD2} is the invariant equation that a $C^\rho$-weakly asymptotic cylinder has to satisfy. Moreover, ~\eqref{hyp2MD2} implies that $\varphi^t$ converges at time tends to infinity to $\varphi_0$ and that $\Gamma$ decays in time as ${1 \over t}$. 

Let us recall the Hamiltonian in~\eqref{HMD} in order to give an idea of the proof of Theorem \ref{MD}. We assume that the following positive parameters $s$, $\lambda$, $\rho$, $\beta$, and $\alpha$,  satisfy the conditions~\eqref{parametersMD}.  For a given $\omega \in \R^n$ and real positive parameters $0 \le \delta < 1$, $\varepsilon > 0$ and $\Upsilon \ge 1$, we consider the following time-dependent Hamiltonian
\begin{equation*}
\begin{cases}
H : \T^n   \times \R^m \times B \times J \longrightarrow \R\\
H(q,p,t) =  \omega \cdot p + a(q,t) + b(q,t) \cdot p + m(q,p,t) \cdot p^2\\
b(q,t) = b_0(q,t) + b_r(q,t)\\
 a \in \mathcal{S}_{s,(0,2)}, \hspace{2mm}  b_0, b_r \in \mathcal{S}_{s+1, 1}, \hspace{2mm}  \partial_p^2 H  \in \mathcal{S}_{s+1, 0}\\
|b_0|_{2,1} \le \delta, \quad |b_0|_{s+1,1} \le \Upsilon,\\
|a|_{\lambda,(0, 2)} < \varepsilon,  \quad |b_r|_{\lambda+1,1} < \varepsilon,\\
|a|_{s,(0, 2)} \le \Upsilon,  \quad |b_r|_{s+1,1} \le \Upsilon, \quad |\partial_p^2 H|_{s+1, 0} \le \Upsilon
\end{cases}
\end{equation*}
where we refer to Section \ref{FS} for the definition of the Banach spaces and the norms we used.
Here, unlike in~\eqref{HMD}, we have set $e=0$. We can do it without loss of generality. 
We need to introduce the following notation.  Let $H$ be the previous Hamiltonian, we define
\begin{equation}
\label{barm}
\bar m(q,p,t) = \int_0^1 \partial_p^2 H(q, \tau p, t) d\tau
\end{equation}
for all $(q,p,t) \in \T^n \times \R^m \times B \times J$.  In what follows, we introduce a suitable functional $\mathcal{F}$ depending on the perturbative terms $(a,b) \in \mathcal{S}_{0, (0,2)} \times \mathcal{S}_{1, 2}$, the quadratic terms $(m, \bar m) \in \mathcal{S}_{s+1, 0} \times \mathcal{S}_{s+1, 0}$ and the components of the $C^\rho$-weakly asymptotic cylinder $v \in \mathcal{V}_{0, \omega,1}$ we are looking for.  We refer to Section \ref{FS} for the above notation. 

 The functional $\mathcal{F}$ is defined in such a way that if $\mathcal{F}(a,b,m, \bar m, v)=0$ then, letting $\Gamma = b + \left(\bar m \circ \tilde\varphi\right)v$, the found $v$ and $\Gamma$ satisfy~\eqref{hyp1MD2}. We refer to~\eqref{tildef} for the notation $\tilde \varphi$ and to~\eqref{F} for the definition of the functional $\mathcal{F}$.  In Section \ref{DF}, we will see that $\mathcal{F}(0,b,m, \bar m, 0)=0$ for all $(b, m, \bar m) \in \mathcal{S}_{1, 2} \times \mathcal{S}_{s+1, 0} \times \mathcal{S}_{s+1, 0}$.  In particular, 
\begin{equation*}
\mathcal{F}(0,b_0,m, \bar m, 0)=0
\end{equation*}
for all $(m, \bar m) \in \mathcal{S}_{s+1, 0} \times \mathcal{S}_{s+1, 0}$, where $b_0$ is the term defined by~\eqref{HMD} that determines, together with $\omega$,  the dynamic associated with $h$ (see~\eqref{tildehBD}) on the invariant cylinder $\varphi_0$. We can reformulate our dynamical problem in the following terms.  For fixed $b_0$, $m$, $\bar m$ and for $(a,b)$ sufficiently close to $(0,b_0)$, we are looking for $v$ in such a way that $\mathcal{F}(a,b,m, \bar m, v)=0$.  Section \ref{DF} is dedicated to the definition of the functional $\mathcal{F}$.

The proof relies on the Nash-Moser theorem proved by Zehnder~\cite{Zeh76} (see Theorem \ref{NMZ} in Section \ref{Zehnder}).  For this reason, in the first part of Section \ref{PSMD},  we fix $m$ and $\bar m$ as in~\eqref{HMD} and~\eqref{barm}, respectively, and we verify that the above-mentioned functional $\mathcal{F}(\cdot, \cdot, m, \bar m, \cdot)$ satisfies the hypotheses of Theorem \ref{NMZ} (see Lemma \ref{LemmaregolaritaFMD}).  Let $D_v\mathcal{F}$ be the differential of $\mathcal{F}$ with respect to $v$.  Among the hypotheses of Theorem \ref{NMZ}, a key requirement is the existence of a right inverse for $D_v\mathcal{F}$. The proof of this condition relies on the solution of the homological equation~\eqref{HEMD}, particularly on Lemma \ref{homoeqlemmaMD}.
We know that the existence of a solution to the problem~\eqref{HEMD} is contingent upon a smallness assumption (see~\eqref{mu}).  For this reason,  in the proof of Lemma \ref{LemmaregolaritaFMD}, we need to fix $\delta$ to be suitably small (see~\eqref{delta}) in order to apply Lemma \ref{homoeqlemmaMD} and prove the existence of a right inverse of $D_v\mathcal{F}$. 

The application of Theorem \ref{NMZ}, provides the existence of a $C^\rho$-weakly asymptotic cylinder $\varphi^t$ associated to $(X_H, X_{h}, \varphi_0)$. 
In the second part of Section \ref{PSMD}, we verify that the found $C^\rho$-weakly asymptotic cylinder is Lagrangian (we refer to Lemma \ref{LagrBDLemma} for the proof and to Definition \ref{weakasymcyl} for the notion of Lagrangian $C^\rho$-weakly asymptotic cylinder). This concludes the proof of Theorem \ref{MD}.

\subsection{Definition of the functional $\mathcal{F}$}\label{DF}

This section is dedicated to the definition of the above-mentioned functional $\mathcal{F}$. First, we need to verify that 
\begin{equation}
\label{barmp=pmp2}
\bar m(q,p,t) p = \partial_p \Big(m(q,p,t) \cdot p^2 \Big)
\end{equation}
for all $(q,p,t) \in \T^n \times \R^m \times B \times J$.  To this end, letting $H$ be the Hamiltonian in~\eqref{HMD}, thanks to the Taylor formula, we can see that 
\begin{eqnarray}
H(q,p,t) &=& H(q,0,t) + \partial_p H(q,0,t) \cdot p \nonumber\\
\label{Provabarm1}
&+&  \int_0^1 (1 -\tau) \partial_p^2 H(q,\tau p, t) d\tau \cdot p^2\\
\label{Provabarm2}
\partial_p H(q,p,t) &=& \partial_p H(q,0,t) + \int_0^1 \partial_p^2 H(q,\tau p, t) d\tau \cdot p,
\end{eqnarray}
whereas differentiating~\eqref{Provabarm1} with respect to $p$, we obtain
\begin{equation}
\label{Provabarm3}
\partial_p H(q,p,t) = \partial_p H(q,0,t) + \partial_p\left( \int_0^1 (1 -\tau) \partial_p^2 H(q,\tau p, t) d\tau\cdot p^2\right).
\end{equation}
Now,  we compare~\eqref{Provabarm2} with~\eqref{Provabarm3} and we have that
\begin{equation*}
\int_0^1 \partial_p^2 H(q,\tau p, t) d\tau \cdot p = \partial_p\left(\int_0^1 (1 -\tau) \partial_p^2 H(q,\tau p, t) d\tau \cdot p^2\right).
\end{equation*}
This concludes the proof of~\eqref{barmp=pmp2}. Now,  in order to define the functional $\mathcal{F}$, we observe that the Hamiltonian system associated with the Hamiltonian $H$  is equal to 
\begin{equation*}
X_H(q,p,t)  = \begin{pmatrix} \bar\omega + b(q,t) + \bar m(q,p,t) p  \\
-\partial_q a(q, t) - \partial_q b(q,t)p - \partial_q m(q,p,t) p^2\end{pmatrix},
\end{equation*}
for all $(q,p,t) \in \T^n \times \R^m \times B \times J$, where we refer to~\eqref{BarOmega} for the notation $\bar \omega$.  Letting, $\varphi$ the family of embeddings introduced by~\eqref{VarphiThmA}, the composition $X_H \circ \tilde \varphi$ takes the following form
\begin{equation}
\label{DefF1}
X_H \circ \tilde\varphi(q,t) = \begin{pmatrix}\bar \omega + b(q,t) + \Big(\bar m \circ \tilde \varphi (q,t)\Big) v(q, t)  \\
-\partial_q a(q, t) - \partial_q b (q,t) v(q,t) - \partial_qm \circ \tilde \varphi(q,t) v^2(q,t)\end{pmatrix}
\end{equation}
for all $(q, t) \in \T^n \times \R^m \times J$, where we refer to~\eqref{tildef} for the notation $\tilde \varphi$.  On the other hand,  for all $(q, t) \in \T^n \times \R^m \times J$,
\begin{equation}
\label{DefF2}
\partial_q \varphi(q,t)\left(\bar \omega + \Gamma(q,t) \right) +\partial_t \varphi(q,t) = \begin{pmatrix} \bar \omega + \Gamma(q,t) \\
 \partial_q v(q,t)(\omega + \Gamma(q,t)) + \partial_t v(q,t) \end{pmatrix}.
\end{equation}
Hence,  using~\eqref{DefF1} and~\eqref{DefF2}, we can rewrite~\eqref{hyp1MD2} in the following form
\begin{eqnarray}
\label{hyp1MD3}
\begin{pmatrix} \Gamma - b - \left(\bar m \circ \tilde \varphi \right) v \\
\partial_q a + \left(\partial_q b \right)  v + \partial_qm \circ \tilde \varphi\cdot  v^2-\partial_q v(\omega + \Gamma)- \partial_t v
\end{pmatrix} = \begin{pmatrix} 0\\0\end{pmatrix}.
\end{eqnarray}
The latter is composed of sums and products of functions defined on $(q,t) \in \T^n \times \R^m\times J$. We have omitted the arguments $(q,t)$ in order to achieve a more elegant form. We keep this notation for the rest of the proof. 

We need to introduce the following family of Banach spaces  $\left\{\left(\mathcal{X}_\sigma, |\cdot |_{\sigma}\right)\right\}_{\sigma \ge 0}$ such that,  for all $\sigma \ge 0$, 
\begin{equation*}
\mathcal{X}_{\sigma} = \mathcal{S}_{\sigma, (0,2)} \times \mathcal{S}_{\sigma+1, 2}
\end{equation*}
and for all $x=(a,b) \in \mathcal{X}_{\sigma}$, 
\begin{equation*}
|x|_\sigma = \max\{|a|_{\sigma, (0,2)}, |b|_{\sigma+1,1}\}.
\end{equation*}
Let $\mathcal{F}$ be the following functional
\begin{equation}
\label{F}
\mathcal{F} : \mathcal{X}_0 \times \mathcal{S}_{s+1,0} \times \mathcal{S}_{s+1,0} \times \mathcal{V}_{0, \omega, 1} \longrightarrow \mathcal{S}_{0, 2}
\end{equation}
\begin{equation*}
\mathcal{F}(x, m, \bar m,v) = Dv\bar \Omega + \partial_q v  \left(b + \left(\bar m \circ \tilde \varphi\right) v \right) + \partial_q a + \left(\partial_q b \right) v + \partial_q m \circ \tilde \varphi \cdot v^2
\end{equation*}
where $s$ is the positive parameter in~\eqref{HMD} and, for the notation used in the definition of the latter,  we refer to Section \ref{FS}.  It is obtained by the second equation of~\eqref{hyp1MD3}, wherein we substitute $\Gamma$ with $b + \left(\bar m \circ \tilde \varphi \right) v$.  Using  Proposition \ref{normpropertiesMD}, it is straightforward to verify that $\mathcal{F}$ is well defined.  

We observe that, for all $(b, m, \bar m) \in \mathcal{S}_{0, 2} \times \mathcal{S}_{s+1,0} \times \mathcal{S}_{s+1,0}$, letting $x_{0,b} = (0, b)$, we have that $\mathcal{F}(x_{0,b}, m, \bar m, 0) = 0$.  More specifically, let $x_0 =(0,b_0)$, where $b_0$ is the unperturbed term in~\eqref{HMD}, we have that
\begin{equation*}
\mathcal{F}(x_0,m,\bar m,0) = 0.
\end{equation*}
Now, we fix $x_0 = (0,b_0)$. For all $\sigma \ge 0$ and for a suitable parameter $0 < \zeta <1$, that we will speficy later,
we define the following subset of $\mathcal{X}_\sigma \times \mathcal{V}_{\sigma, \omega,1}$
\begin{equation*}
\mathcal{O}_\sigma^\zeta  = \{(x, v) \in \mathcal{X}_\sigma \times \mathcal{V}_{\sigma, \omega,1} : |x - x_0|_\sigma, |v|_{\sigma, \omega, 1} < \zeta \}.
\end{equation*}
We fix $(m,\bar m) \in \mathcal{S}_{s+1,0} \times \mathcal{S}_{s+1,0}$ as in~\eqref{HMD}, we consider the following continuous functional
\begin{equation}
\label{Fmbarm}
\mathcal{F}_{m, \bar m} : \mathcal{O}_0^\zeta  \longrightarrow \mathcal{S}_{0, 2}, \quad \mathcal{F}_{m, \bar m}(x,v) = \mathcal{F}(x,m, \bar m,v).
\end{equation}
In the following section, we will prove that $\mathcal{F}_{m, \bar m}$ satisfies the hypotheses of Theorem \ref{NMZ}.  As mentioned before, an essential point of the proof consists in showing the existence of a right inverse of the differential of $\mathcal{F}_{m, \bar m}$ with respect to the variable $v$.  Let 
\begin{eqnarray*}
f &=& b + \left(\bar m \circ \tilde \varphi \right) v,\\
g &=& \partial_q b + \partial_q v \left( \partial_p \bar m \circ \tilde \varphi \right) v+ \partial_q v \left( \bar m \circ \tilde \varphi \right) + v^T \left(\partial^2_{pq} m \circ \tilde \varphi  \right) v + 2 \left(\partial_q m \circ \tilde \varphi \right) v,
\end{eqnarray*}
where $T$ denotes the transpose.  In Section \ref{PSMD}, we will see that the differential of $\mathcal{F}_{m, \bar m}$ with respect to the variable $v$ is equal to 
\begin{eqnarray*}
D_v \mathcal{F}_{m, \bar m} (x,v) \hat v = D\hat v\bar \Omega+  \left(\partial_q \hat v \right) f+ g \hat v.
\end{eqnarray*}
As mentioned before, we can find a right inverse for the latter. The proof relies on the solution of the equation~\eqref{HEMD}, and hence on Lemma  \ref{homoeqlemmaMD}.  For this purpose, in the proof of Lemma \ref{LemmaregolaritaFMD},  first, we fix $\delta$ small enough (see~\eqref{delta}), and then we fix $\zeta$ sufficiently small with respect to $\delta$ (see~\eqref{zeta}) in order to make $f$ and $g$ small enough to apply Lemma \ref{LemmaregolaritaFMD} and prove the existence of a right inverse of $D_v \mathcal{F}_{m, \bar m}$.

\subsection{Existence of a Lagrangian $C^\rho$-weakly asymptotic cylinder}\label{PSMD}

This section is divided into two parts. First, we prove that the functional $\mathcal{F}_{m, \bar m}$, defined by~\eqref{Fmbarm}, satisfies the hypotheses of Theorem \ref{NMZ} (we refer to Lemma \ref{LemmaregolaritaFMD}).  As mentioned before, the application of Theorem \ref{NMZ} ensures the existence of a $C^\rho$-weakly asymptotic cylinder associated to $(X_H, X_{h}, \varphi_0)$.  In the second part, we will verify that the found $C^\rho$-weakly asymptotic cylinder is Lagrangian (see Lemma \ref{LagrBDLemma}). 

In the proof of the following lemma, we widely use the notation in Section \ref{FS} and the properties contained in Proposition \ref{normpropertiesMD} without specifying each time.  

\begin{lemma}
\label{LemmaregolaritaFMD}
$\mathcal{F}_{m, \bar m}$ satisfies hypotheses H.1-H.4 of Theorem \ref{NMZ}. 
\end{lemma}
\begin{proof}
\textit{H.1. Smoothness}: For all fixed $x \in \mathcal{X}_0$, one can see that,  the functional $\mathcal{F}_{m, \bar m}(x, \cdot) : \mathcal{V}_{0, \omega, 1} \to \mathcal{S}_{0,2}$ is two times differentiable with respect to the variable $v$ and for all $(x,v) \in \mathcal{X}_0 \times \mathcal{V}_{0, \omega, 1}$ 
\begin{equation*}
D_v \mathcal{F}_{m, \bar m}(x,v) :\mathcal{V}_{0, \omega, 1} \to \mathcal{S}_{0,2}, \quad D^2_v \mathcal{F}_{m, \bar m}(x,v) :\mathcal{V}_{0, \omega, 1} \times \mathcal{V}_{0, \omega, 1} \to \mathcal{S}_{0,2}
\end{equation*}
with
\begin{eqnarray}
\label{DiffF}
D_v \mathcal{F}_{m, \bar m}(x,v) \hat v &=& D\hat v \bar\Omega+  \left(\partial_q \hat v \right) f+ g \hat v,\\
\label{Diff2F}
D^2_v \mathcal{F}_{m, \bar m}(x,v)(\hat v_1, \hat v_2) &=&\left(\partial_q \hat v_1\right) h_1  \hat v_2+ \left(\partial_q \hat v_2 \right) h_2 \hat v_1 + \hat v_2^T h_3  \hat v_1,
\end{eqnarray}
where, letting $T$ be the transpose, the elements $f$, $g$, $h_1$, $h_2$ and $h_3$ are defined by
\begin{eqnarray}
f &=& b + \left(\bar m \circ \tilde \varphi \right) v,\nonumber\\
g &=& \partial_q b + \partial_q v \left( \partial_p \bar m \circ \tilde \varphi \right) v+ \partial_q v \left( \bar m \circ \tilde \varphi \right) + v^T \left(\partial^2_{pq} m \circ \tilde \varphi  \right) v + 2 \left(\partial_q m \circ \tilde \varphi \right) v,\nonumber\\
\label{fgh123}
h_1 &=& \bar m \circ \tilde \varphi  + \left(\partial_p \bar m \circ \tilde \varphi \right) v,\\
h_2&=& \left(\partial_p \bar m \circ \tilde \varphi \right) v + \bar m \circ \tilde \varphi,\nonumber\\
h_3 &=& \partial_q v \left(\partial^2_p m \circ \tilde \varphi \right) v + 2\partial_q v \left( \partial_p \bar m \circ \tilde \varphi \right)  + v^T \left( \partial_p^2 \partial_q m \circ \tilde \varphi \right) v \nonumber\\
&+& 4 \left(\partial_{pq}^2 m \circ \tilde \varphi \right) v + 2 \partial_q m \circ \tilde \varphi.\nonumber
\end{eqnarray}
For all fixed $(x,v) \in \mathcal{O}_0^{\zeta}$, we have to verify the existence of a constant $C_* \ge 1$ such that
\begin{equation*}
|D_v \mathcal{F}_{m, \bar m} (x,v)\hat v|_{0,2} \le C_*|\hat v|_{0,\omega, 1}, \quad |D^2_v \mathcal{F}_{m, \bar m} (x,v)(\hat v_1, \hat v_2)|_{0,2} \le C_*|\hat v_1|_{0,\omega, 1}|\hat v_2|_{0,\omega, 1},
\end{equation*}
for all $\hat v$, $\hat v_1$, $\hat v_2 \in \mathcal{V}_{0, \omega, 1}$.  For this purpose,  for all fixed $(x,v) \in \mathcal{O}_0^{\zeta}$ and for all $\hat v \in \mathcal{V}_{0, \omega,1}$,  
\begin{eqnarray*}
|D_v \mathcal{F}_{m, \bar m} (x,v)\hat v|_{0,2} &\le& | D\hat v \bar\Omega|_{0, 2} + C\left(|f|_{0,1}|\hat v|_{1,1} + |g|_{0,1}|\hat v|_{0,1}\right)\\
&\le& |\hat v|_{0, \omega, 1} \left(1 + C\left(|f|_{0,1} +  |g|_{0,1} \right)\right).
\end{eqnarray*}
Now, we have to estimates $|f|_{0,1}$ and $|g|_{0,1}$. Hence, for all $(x,v) \in \mathcal{O}_0^{\zeta}$,  thanks to~\eqref{HMD} and~\eqref{barm}, we have the following upper bounds
\begin{eqnarray*}
|f|_{0,1} &\le& C \left(|b|_{0,1} + |v|_{0,1}|\bar m|_{0,0}\right)\\
&\le& C \left(|b_0|_{0,1} + |b - b_0|_{0,1}\right) + C\left(|v|_{0,1}|\bar m|_{0,0}\right)\le C(\delta + \zeta) + C\Upsilon \zeta \\
|g|_{0,1} &\le& C \left(|b|_{1,1} + |v|_{0,1}|\bar m|_{1,0}|v|_{1,1} +  |v|_{1,1}|\bar m|_{0,0}  +  \left(|v|_{0,1}\right)^2|\bar m|_{2,0} +  |v|_{0,1}|m|_{1,0}\right)\\
&\le& C \left(|b|_{1,1} + (|v|_{0, \omega, 1} + \left(|v|_{0, \omega, 1}\right)^2)\Upsilon \right) \le C (\delta + \zeta) + C \Upsilon \zeta.
\end{eqnarray*}
where we recall that $b_0$ is defined by~\eqref{HMD}. This implies the claim for $D_v \mathcal{F}_{m, \bar m} (x,v)$. Similarly, we have the claim for $D^2_v \mathcal{F}_{m, \bar m} (x,v)$.

\vspace{5mm}
\textit{H.2. $\mathcal{F}_{m, \bar m}$ is uniformly Lipschitz in $\mathcal{X}_0$}: For all $(x_1,v)$, $(x_2, v) \in  \mathcal{O}_0^\zeta$, remembering that $|x|_0 = \max\{|a|_{0, (0,2)}, |b|_{1,1}\}$,
\begin{eqnarray*}
|\mathcal{F}_{m, \bar m}(x_1,v) - \mathcal{F}_{m, \bar m}(x_2,v)|_{0,2} &=& | \partial_q v (b_1 - b_2)+ (\partial_q a_1 - \partial_q a_2) + (\partial_q b_1 - \partial_q b_2) v|_{0,2}\\
&\le&C|b_1 - b_2|_{0,1}|v|_{1,1}\\
&+& C\left(|\partial_q a_1 - \partial_q a_2|_{0,2} + |b_1 - b_2|_{1,1}|v|_{0,1} \right)\\
&\le&C(1 + \zeta)|x_1 - x_2|_0,
\end{eqnarray*}
which proves H.2. Now, we verify H.4 before H.3.

\vspace{5mm}
\textit{H.4. Order}: One can see that the first two conditions of Definition \ref{order} are satisfied, meaning $(x_0,0) \in \mathcal{X}_s \times \mathcal{V}_{s, \omega, 1}$ and $\mathcal{F}_{m, \bar m}\left(\mathcal{O}_0^\zeta \bigcap\left(\mathcal{X}_\sigma \times \mathcal{V}_{\sigma, \omega, 1} \right)\right)\subset \mathcal{S}_{\sigma, 2}$ for all $1 \le \sigma \le s$. We verify the tame estimate. 

For all $1 \le \sigma \le s$ and  $(x,v) \in \mathcal{O}_1^\zeta \bigcap \left(\mathcal{X}_\sigma \times \mathcal{V}_{\sigma, \omega, 1} \right)$, we rewrite the functional $\mathcal{F}_{m, \bar m}$ in the following form 
\begin{eqnarray*}
\mathcal{F}_{m, \bar m}(x,v) &=& Dv \bar \Omega + \partial_q v  \left(b_0 + (b-b_0) + \left(\bar m \circ \tilde \varphi\right) v \right) + \partial_q a\\
&+& \left(\partial_q b_0 \right) v + \partial_q\left(b -b_0\right) v + \partial_q m \circ \tilde \varphi \cdot v^2.
\end{eqnarray*}
For all $1 \le \sigma \le s$, we assume $|x - x_0|_\sigma$, $|v|_{\sigma, \omega, 1} \le K$ (we recall that $x_0 = (0, b_0)$). Then, we can estimate $|\mathcal{F}_{m, \bar m}(x,v)|_{\sigma,2}$ as follows
\begin{eqnarray}
|\mathcal{F}_{m, \bar m}(x,v)|_{\sigma,2} &\le& |Dv \bar\Omega|_{\sigma,2} + | \left( \partial_qv \right)b_0|_{\sigma,2} + | \left( \partial_qv \right)\left(b - b_0\right)|_{\sigma,2}\nonumber\\
&+& |\partial_q v \left(\bar m \circ \tilde \varphi \right) v|_{\sigma,2} + |a|_{\sigma,2} + |\left(\partial_q b_0 \right) v|_{\sigma,2} \nonumber\\
\label{H3FTameEstimate}
&+& |\partial_q\left(b - b_0\right) v|_{\sigma,2} + |\partial_q m \circ \tilde \varphi \cdot v^2|_{\sigma,2}.
\end{eqnarray}
We have to estimate each term on the right-hand side of the latter. We already known that $|a|_{\sigma,2}$ and $|Dv\bar \Omega|_{\sigma,2}$ are bounded by $K$. We recall that $|b_0|_{s+1,1} \le \Upsilon$, $|x - x_0|_\sigma$, $|v|_{\sigma, \omega, 1} \le K$, and $0< \zeta <1$ and we estimate the others

\begin{eqnarray*}
 |b_0 \left( \partial_qv \right)|_{\sigma,2} &\le& C(\sigma) |b_0|_{s+1,1}|v|_{\sigma+1,1} \le  C(\sigma) \Upsilon|v|_{\sigma, \omega, 1} \le C(\sigma) \Upsilon K\\
 | \left( \partial_qv \right)\left(b - b_0\right)|_{\sigma,2} &\le& C(\sigma)\left(|\partial_qv |_{0,1}|b - b_0|_{\sigma,1}  + | \partial_qv|_{\sigma,1} |b - b_0|_{0,1} \right)\\
&\le& C(\sigma)\left(|v |_{1,1}|b - b_0|_{\sigma,1}  + |v|_{\sigma + 1,1} |b - b_0|_{0,1} \right)\\
&\le& C(\sigma)\zeta K  \le C(\sigma) K\\
 |\partial_q v \left(\bar m \circ \tilde \varphi \right) v|_{\sigma,2} &\le& C(\sigma) \left(|\partial_qv\left(\bar m \circ \tilde \varphi \right) |_{\sigma,1} |v|_{0,1}+ |\partial_qv\left(\bar m \circ \tilde \varphi \right) |_{0,1}|v|_{\sigma,1}\right)\\
&\le& C(\sigma)|v|_{0, \omega, 1}\left(|\bar m \circ \tilde \varphi|_{\sigma,0}|v|_{0, \omega, 1} + |\bar m|_{0,0}|v|_{\sigma, \omega, 1}\right)\\
&+&C(\sigma)|v|_{\sigma, \omega, 1}|\bar m|_{0,0}|v|_{0, \omega, 1}\\
&\le& C(\sigma) \zeta^2 |\bar m|_{\sigma,0}(1 + |v|^\sigma_{0, \omega, 1} +  |v|_{\sigma, \omega, 1}) + C(\sigma)\zeta\Upsilon K \\
&\le& C(\sigma)\Upsilon K\\
|\left(\partial_q b_0 \right) v|_{\sigma,2} &\le& C(\sigma)|b_0|_{s+1,1}|v|_{\sigma,1} \le C(\sigma)\Upsilon K \\
|\partial_q\left(b - b_0\right) v|_{\sigma,2} &\le& C(\sigma)\left(|v |_{0,1}|b - b_0|_{\sigma+1,1}  + |v|_{\sigma,1} |b - b_0|_{1,1} \right)\\
&\le& C(\sigma)\zeta K  \le C(\sigma) K\\
|\partial_q m \circ \tilde \varphi \cdot v^2|_{\sigma,2}&\le& C(\sigma) \left(|\left(\partial_q m \circ \tilde \varphi\right) v|_{0,1} |v|_{\sigma,1} + |\left(\partial_q m \circ \tilde \varphi\right) v|_{\sigma,1} |v|_{0,1} \right)\\
&\le&C(\sigma)|m|_{1,0}|v|_{0, \omega, 1}|v|_{\sigma, \omega, 1} \\
&+&C(\sigma)\zeta\left(|\partial_q m \circ \tilde \varphi|_{\sigma,0}|v|_{0,1} + |\partial_q m \circ \tilde \varphi|_{0,0}|v|_{\sigma,1} \right)\\
&\le& C(\sigma) \zeta\Upsilon K \\
&+& C(\sigma) \zeta^2|m|_{\sigma+1,0}(1 + |v|^\sigma_{0, \omega, 1} + |v|_{\sigma, \omega, 1})+C(\sigma) \Upsilon K\\
&\le& C(\sigma) \Upsilon K.
\end{eqnarray*}
Replacing the above estimates in~\eqref{H3FTameEstimate}, one has that 
\begin{equation}
\label{EstimateFH4}
|\mathcal{F}_{m, \bar m}(x,v)|_{\sigma,2} \le C(\sigma) \Upsilon K.
\end{equation}
This concludes the proof of \textit{H.4}. 

\vspace{5mm}
\textit{H.3. Existence of a right-inverse of loss 1}: In this part, we prove the existence of a right inverse for the differential $D_v \mathcal{F}_{m, \bar m}$.  First,  let $C_1$ and $C_2$ be the positive constant defined by~\eqref{f1} and Lemma \ref{phiMDfineGamma}, respectively. We stress that $C_1$ and $C_2$ are explicit and depend only on $n+m$.  Furthermore, let $c^\chi_\sigma$ the constant introduced in Lemma \ref{homoeqlemmaMD}.  We consider the positive parameter $0 \le \delta <1$ in~\eqref{HMD}. We fix $\delta$ in such a way that 
\begin{equation}
\label{delta}
2 \max\{C_1, C_2\} \delta < \inf_{1 \le \sigma \le s}{1 \over c^\chi_\sigma},
\end{equation}
where $s$ is the positive parameter defined by~\eqref{HMD}.  Let $\rho$ be the parameter in~\eqref{HMD} and $K_\rho$ be the constant defined by~\eqref{GammaEst}.  We stress that $K_\rho$ is explicit and only depends on $\rho$ and $n+m$. We choose $0 < \zeta < 1$ small enough with respect to $\delta$, $\Upsilon$, $\rho$ and $s$ such that
\begin{align}
\label{zeta}
&2\max\{C_1, C_2\}\left(\delta + \Upsilon \zeta\right)< \inf_{1 \le \sigma \le s}{1 \over c^\chi_\sigma},\\
\label{zeta2}
& \delta + K_\rho \Upsilon \zeta < 1.
\end{align}
We point out that~\eqref{zeta} is a necessary condition in order to prove the existence of a right inverse for $D_v \mathcal{F}_{m, \bar m}$.  We will also use~\eqref{zeta} in the proof of Lemma \ref{LagrBDLemma}.  On the other hand,~\eqref{zeta2} is a necessary condition to prove~\eqref{StimevGamma}.  To avoid confusion, we prefer to fix the parameters $\delta$ and $\zeta$ at this stage of the proof.

Now, we want to prove that for all $(x,v) \in \mathcal{O}^\zeta_1 \cap\left(\mathcal{X}_\sigma \times \mathcal{V}_{\sigma, \omega, 1}\right)$ with $1 \le \sigma \le s$, a right-inverse of loss $1$ of $D_v \mathcal{F}_{m, \bar m} (x,v)$ exists and is well defined. This means that,  for all $1 \le \sigma \le s$ and $(x,v) \in \mathcal{O}^\zeta_1 \cap\left(\mathcal{X}_\sigma \times \mathcal{V}_{\sigma, \omega, 1}\right)$,  we want to verify the existence of a liner map $\eta_{m, \bar m}(x,v) : \mathcal{S}_{\sigma, 2} \to \mathcal{V}_{\sigma-1, \omega, 1}$ in such a way that $D_v \mathcal{F}_{m, \bar m}(x,v)  \eta_{m, \bar m} (x,v)z = z$ for all $z \in \mathcal{S}_{\sigma, 2}$.  For all $1 \le \sigma \le s$ and $z \in \mathcal{S}_{\sigma, 2}$, it consists in solving the following equation in the unknown $\hat v:\T^n \times \R^m \times J \to \R^{n+m}$
\begin{equation}
\label{eqinv}
D\hat v\bar\Omega+  \left(\partial_q \hat v \right) f+ g \hat v = z,
\end{equation}
in such a way that $\hat v \in \mathcal{V}_{\sigma-1, \omega, 1}$, where $f$ and $g$ are defined by~\eqref{fgh123}. 

Thanks to  Lemma \ref{homoeqlemmaMD},  if 
\begin{equation}
\label{fgH41}
|f|_{1,1}, |g|_{1,1}  < \inf_{1 \le \sigma \le s}{1 \over c^\chi_\sigma}
\end{equation}
then, for all $1 \le \sigma \le s$, a solution $\hat v$ of~\eqref{eqinv} exists and $\hat v \in \mathcal{V}_{\sigma-1, \omega, 1}$.  It remains to verify~\eqref{fgH41}.  To this end, for all $1 \le \sigma \le s$, $(x,v) \in \mathcal{O}^\zeta_1 \cap\left(\mathcal{X}_\sigma \times \mathcal{V}_{\sigma, \omega, 1}\right)$, and remembering that $0<\zeta<1$
\begin{eqnarray}
|f|_{\sigma,1} &\le& |b|_{\sigma, 1} + | \left(\bar m \circ \tilde \varphi \right) v|_{\sigma, 1} \le |b|_{\sigma, 1} + C(\sigma) \left( |\bar m \circ \tilde \varphi|_{\sigma, 0} |v|_{0,1} +  |\bar m \circ \tilde \varphi|_{0, 0} |v|_{\sigma ,1} \right)\nonumber\\
&\le& |b|_{\sigma, 1} +C(\sigma)\left( \Upsilon\left(1 + |v|^\sigma_{1, 1} + |v|_{\sigma, 1}\right)|v|_{0,1} + \Upsilon |v|_{\sigma , 1}\right)\nonumber\\ \label{fsigma}
&\le&  |b|_{\sigma, 1} +C(\sigma)\Upsilon |v|_{\sigma, 1} \\
|g|_{\sigma,1} &\le& |b|_{\sigma+1,1} + |\partial_q v \left( \partial_q\bar m \circ \tilde \varphi \right) v|_{\sigma,1} + |\partial_q v \left(\bar m \circ \tilde \varphi \right)|_{\sigma,1} + |v^T \left( \partial^2_{pq}\bar m \circ \tilde\varphi \right) v|_{\sigma,1}\nonumber\\
&+& |\left(\partial_q m \circ \tilde \varphi \right) v|_{\sigma,1}\nonumber\\ \label{gsigma}
&\le& |b|_{\sigma+1,1} + C(\sigma)\Upsilon|v|_{\sigma+1,1}.
\end{eqnarray}
Hence, taking $\sigma=1$, by~\eqref{HMD} and $(x,v) \in \mathcal{O}^\zeta_1$ we obtain 
\begin{eqnarray}
|f|_{1,1} &\le& |b|_{1,1} +C\Upsilon |v|_{1,1} \le |b_0|_{1,1} + |b-b_0|_{1,1} + C\Upsilon |v|_{1,1}\nonumber\\
\label{f1}
&\le& (\delta + \zeta) + C\Upsilon \zeta \le \delta + C_1\Upsilon \zeta\\
|g|_{1,1} &\le& |b|_{2,1} + C\Upsilon|v|_{2,1} \le |b_0|_{2,1} + |b-b_0|_{2,1} + C\Upsilon |v|_{2,1},\nonumber \\
\label{g1}
&\le& (\delta + \zeta) + C\Upsilon \zeta \le \delta + C_1\Upsilon \zeta
\end{eqnarray}
for a suitable constant $C_1\ge 1$ depending on $n+m$.  Thanks to~\eqref{delta},~\eqref{zeta} and the above estimates,~\eqref{fgH41} is satisfied. Then,  for all $1 \le \sigma \le s$ and $z \in \mathcal{S}_{\sigma, 2}$, the application of Lemma \ref{homoeqlemmaMD} implies the existence of a solution $\hat v \in \mathcal{V}_{\sigma-1, \omega, 1}$ of~\eqref{eqinv} and hence the existence of a right inverse of loss $1$ of $D_v \mathcal{F}_{m, \bar m} (x,v)$ for all $(x,v) \in \mathcal{O}^\zeta_1 \cap\left(\mathcal{X}_\sigma \times \mathcal{V}_{\sigma, \omega, 1}\right)$ with $1 \le \sigma \le s$.

The second part of this proof is dedicated to verifying~\eqref{eta1} and~\eqref{eta2}. In what follows, we drop the indexes ${m, \bar m}$  from $\mathcal{F}$ and $\eta$ to achieve a more elegant proof.

We begin with the proof of~\eqref{eta1}.  For this reason, for all $(x ,v) \in \mathcal{O}_\zeta^1$ and $z \in \mathcal{S}_{1, 2}$, we recall that 
\begin{equation*}
|\eta(x,v) z|_{0, \omega, 2} = \max\{|\eta(x,v) z|_{1,1},|D\left(\eta(x,v)z\right) \bar \Omega|_{0,2} \}.
\end{equation*}
We refer to Section \ref{FS}.
In order to verify~\eqref{eta1}, we need to estimate the two norms $|\eta(x,v) z|_{1,1}$ and $|D\left(\eta(x,v)z\right) \bar \Omega|_{0,2}$ on the right hand side of the latter.
By Lemma \ref{homoeqlemmaMD} (more specifically~\eqref{varkappa}),~\eqref{zeta}, ~\eqref{f1} and~\eqref{g1}, one can see that 
\begin{equation}
\label{etaz11}
|\eta(x,v) z|_{1,1} \le C(\delta, \zeta)|z|_{1,2}.
\end{equation}
On the other hand, as a consequence of~\eqref{eqinv},~\eqref{etaz11},~\eqref{f1},~\eqref{g1}  and~\eqref{zeta}, we have
\begin{eqnarray}
|D\left(\eta(x,v) z \right) \bar \Omega|_{0,2} &=& |z -  \partial_q \left(\eta(x,v)z  \right) f - g \left( \eta(x,v)z\right)|_{0,2}\nonumber\\
&\le&|z|_{0, 2} + C|f|_{0,1}|\eta(x,v)z|_{1,1} + C|g|_{0,1}|\eta(x,v)z|_{0,1}\nonumber\\
&\le&|z|_{0, 2} + C|f|_{1,1}|\eta(x,v)z|_{1,1} + C|g|_{1,1}|\eta(x,v)z|_{0,1}\nonumber\\
\label{etaz02}
&\le& |z|_{0, 2} + C|\eta(x,v) z|_{1,1} \le C(\delta, \zeta)|z|_{1,2}.
\end{eqnarray}
In the last line, we utilized the inequalities $|f|_{1,1} \le 1$ and $|g|_{1,1} \le 1$, which are consequences of~\eqref{fgH41}.
We proved~\eqref{eta1} because~\eqref{etaz11} and~\eqref{etaz02} imply
\begin{equation*}
|\eta(x,v) z|_{0, \omega, 2} = \max\{|\eta(x,v) z|_{1,1},|D\left(\eta(x,v)z\right) \bar \Omega|_{0,2} \} \le C(\delta, \zeta)|z|_{1,2}.
\end{equation*}

 Conserning~\eqref{eta2}, for all $1 \le \sigma \le s$ and $(x,v) \in \mathcal{O}_\zeta^1 \cap \left(\mathcal{X}^\sigma \times \mathcal{V}^\sigma\right)$, we assume $|x-x_0|_\sigma$, $|v|_{\sigma, \omega, 1} \le K$ and we recall that 
\begin{equation*}
|\eta(x,v) \mathcal{F}(x,v)|_{\sigma-1, \omega, 1} = \max\{|\eta(x,v)  \mathcal{F}(x,v)|_{\sigma, 1},| D \left(\eta(x,v)  \mathcal{F}(x,v)\right) \bar \Omega|_{\sigma -1, 2} \}.
\end{equation*}
Similarly to the previous case, we shall prove that the two norms on the right-hand side of the latter are smaller or equal to $K$ multiplied by a suitable constant. To this end,  thanks to~\eqref{fsigma} and~\eqref{gsigma}
\begin{eqnarray}
|f|_{\sigma,1} &\le&  |b|_{\sigma, 1} +C(\sigma)\Upsilon |v|_{\sigma, 1} \nonumber\\
\label{fsigma2}
&\le& |b_0|_{\sigma, 1} + |b-b_0|_{\sigma, 1}+C(\sigma)\Upsilon |v|_{\sigma, 1} \le \Upsilon + C(\sigma)\Upsilon K \\
|g|_{\sigma,1} &\le&   |b|_{\sigma+1,1} + C(\sigma)\Upsilon|v|_{\sigma+1,1},\nonumber\\
\label{gsigma2}
&\le& |b_0|_{\sigma+1,1} + |b-b_0|_{\sigma+1, 1}+C(\sigma)\Upsilon |v|_{\sigma+1, 1} \le \Upsilon + C(\sigma)\Upsilon K,
\end{eqnarray}
where we used $|b_0|_{s+1,1} \le \Upsilon$, $\Upsilon \ge 1$ and the assumption $|x-x_0|_\sigma$, $|v|_{\sigma, \omega, 1} \le K$. Furthermore, thanks to H.4 one has that
\begin{equation}
\label{F1Fsigma}
 |\mathcal{F}(x,v)|_{1,2} \le C\Upsilon\zeta, \quad  |\mathcal{F}(x,v)|_{\sigma,2} \le C(\sigma)\Upsilon K.
\end{equation}
We point out that the first estimate in~\eqref{F1Fsigma} is a consequence of hypothesis \textit{H.4} (especially~\eqref{EstimateFH4}) and the fact that $(x,v) \in \mathcal{O}_\zeta^1 \cap \left(\mathcal{X}^\sigma \times \mathcal{V}^\sigma\right)$ implies $|x-x_0|_1$, $|v|_{1, \omega, 1} \le \zeta$. Whereas, the second estimate in~\eqref{F1Fsigma} follows by~\eqref{EstimateFH4} and the assumption $|x-x_0|_\sigma$, $|v|_{\sigma, \omega, 1} \le K$.

Now, by~\eqref{varkappa},~\eqref{fsigma2}, ~\eqref{gsigma2}, and~\eqref{F1Fsigma} we obtain that
\begin{eqnarray}
|\eta(x,v)  \mathcal{F}(x,v)|_{\sigma,1} &\le& C(\sigma,\Upsilon,\delta, \zeta)|\mathcal{F}(x,v)|_{\sigma,2}\nonumber\\
&+& C(\sigma, \Upsilon, \delta, \zeta)\left(|f|_{\sigma, 1} + |g|_{\sigma ,1}\right)|\mathcal{F}(x,v)|_{1,2}\nonumber\\
&\le&C(\sigma, \Upsilon, \delta, \zeta)\left(K + \left(1 + K\right)|\mathcal{F}(x,v)|_{1,2}\right)\nonumber\\
\label{H3StimeB3}
&\le&  C(\sigma,\Upsilon,\delta, \zeta)K,
\end{eqnarray}
where in the last line of the latter we used the following trivial estimate \newline $\left(1 + K\right)|\mathcal{F}(x,v)|_{1,2} \le |\mathcal{F}(x,v)|_{\sigma,2} + K |\mathcal{F}(x,v)|_{1,2} \le C(\sigma) \Upsilon K + KC\Upsilon\zeta$.

On the other hand,  by~\eqref{eqinv}
\begin{eqnarray*}
|D\left(\eta(x,v)  \mathcal{F}(x,v)\right)\bar \Omega|_{\sigma -1,2} &=& |\mathcal{F}(x,v) - \partial_q \left(\eta(x,v)\mathcal{F}(x,v) \right) f - g \left( \eta(x,v)\mathcal{F}(x,v)\right)|_{\sigma -1,2}\\
&\le&|\mathcal{F}(x,v)|_{\sigma-1, 2} + |\partial_q \left(\eta(x,v)\mathcal{F}(x,v) \right) f|_{\sigma-1, 2}\\
&+& |g \left( \eta(x,v)\mathcal{F}(x,v)\right)|_{\sigma-1, 2}.
\end{eqnarray*}
It remains to prove that each term on the right-hand side of the latter can be estimated by $K$ multiplied by a suitable constant independent of $K$. Thanks to \textit{H.4}, this is true for the first term $|\mathcal{F}(x,v)|_{\sigma-1, 2}$ of the latter.  Now, we provide the proof for the remaining terms 
\begin{eqnarray*}
|\partial_q \left(\eta(x,v)\mathcal{F}(x,v) \right) f|_{\sigma-1, 2}&\le& C(\sigma)\left(|\eta(x,v)\mathcal{F}(x,v) |_{\sigma, 1} |f|_{0,1}+ |\eta(x,v)\mathcal{F}(x,v) |_{1, 1}|f|_{\sigma,1}\right)\\
&\le& C(\sigma, \Upsilon, \delta, \zeta)K |f|_{1,1} + C(\delta, \zeta)|\mathcal{F}(u,v)|_{1, 2}|f|_{\sigma,1}\\
&\le& C(\sigma, \Upsilon, \delta, \zeta)K 
\end{eqnarray*}
where the second line of the latter is due to~\eqref{H3StimeB3} and~\eqref{etaz11}, whereas in the last line we use $|f|_{1,1} \le 1$,~\eqref{F1Fsigma} and~\eqref{fsigma2}.  Similarly, one has 
\begin{eqnarray}
|g \left( \eta(x,v)\mathcal{F}(x,v)\right)|_{\sigma-1, 2} &\le& C(\sigma)\left(|g|_{0,1}|\eta(x,v)\mathcal{F}(x,v) |_{\sigma-1, 1} + |g|_{\sigma-1,1}|\eta(x,v)\mathcal{F}(x,v) |_{0, 1}\right)\nonumber\\
&\le& C(\sigma)\left(|g|_{1,1}|\eta(x,v)\mathcal{F}(x,v) |_{\sigma, 1} + |g|_{\sigma,1}|\eta(x,v)\mathcal{F}(x,v) |_{1, 1}\right)\nonumber\\
&\le& C(\sigma)\left(|\eta(x,v)\mathcal{F}(x,v) |_{\sigma, 1} + |g|_{\sigma,1}C(\delta, \zeta)|\mathcal{F}(x,v) |_{1, 2}\right)\nonumber\\
&\le& C(\sigma, \Upsilon, \delta, \zeta)K + C(\delta, \zeta) \Upsilon|\mathcal{F}(x,v) |_{\sigma, 2} \nonumber\\
\label{StimeB3QF}
&+& C(\sigma, \delta, \zeta)\Upsilon K|\mathcal{F}(x,v) |_{1, 2}\\
&\le& C(\sigma, \Upsilon, \delta, \zeta)K.\nonumber
\end{eqnarray}
We point out that in the third line of the latter, we used $|g|_{1,1} \le 1$, and~\eqref{etaz11}. The inequality~\eqref{StimeB3QF} is due to~\eqref{gsigma2} and the trivial estimate $|\mathcal{F}(x,v) |_{1, 2}  \le C|\mathcal{F}(x,v) |_{\sigma, 2} $.  The last line of the latter is a consequence of~\eqref{F1Fsigma}.

Now, thanks to the above estimates, we have that 
\begin{equation}
\label{H3StimeB32}
|D\left(\eta(x,v)  \mathcal{F}(x,v)\right)\bar \Omega|_{\sigma -1,2} \le C(\sigma, \Upsilon, \delta, \zeta)K
\end{equation}
and hence, by~\eqref{H3StimeB3} and~\eqref{H3StimeB32}, for all $1 \le \sigma \le s$ and $(x,v) \in \mathcal{O}_\zeta^1 \cap \left(\mathcal{X}^\sigma \times \mathcal{V}^\sigma\right)$, with $|x-x_0|_\sigma$, $|v|_{\sigma, \omega, 1} \le K$ we obtain that 
\begin{eqnarray*}
|\eta(x,v) \mathcal{F}(x,v)|_{\sigma-1, \omega, 1} &=& \max\{|\eta(x,v)  \mathcal{F}(x,v)|_{\sigma, 1},| D \left(\eta(x,v)  \mathcal{F}(x,v)\right) \bar \Omega|_{\sigma -1, 2} \}\\
&\le&C(\sigma, \Upsilon, \delta, \zeta)K.
\end{eqnarray*}
This concludes the proof of \textit{H.3} and of this lemma. 
\end{proof}

Let $\alpha$, $\beta$, $\lambda$, $\rho$, $s$ be the positive paramters satisying~\eqref{parametersMD}, and $\varepsilon$ the positive parameter in~\eqref{HMD}.  We fix $x = (a,b_r) \in \mathcal{S}_{s, (0,2)} \times \mathcal{S}_{s+1, 1}$, where $a$ and $b_r$ are the perturbative terms in~\eqref{HMD}.  

By Lemma \ref{LemmaregolaritaFMD}, the functional $\mathcal{F}_{m, \bar m}$ satisfies the hypotheses of Theorem \ref{NMZ}.  Then, there exists $\varepsilon_0$, depending on $\alpha$, $\beta$, $\lambda$, $s$ and $\zeta$ such that for all $\varepsilon \le \varepsilon_0$ there exists $v \in \mathcal{V}_{\rho, \omega, 1}$ in such a way that 
\begin{equation}
\label{F=0foundv}
\mathcal{F}_{m, \bar m}(x, v)=0, \quad |v|_{\rho, \omega, 1} \le \zeta.
\end{equation}
For simplicity, without loss of generality,  we can assume $\varepsilon_0 \le \zeta$. Furthermore, we recall that the parameter $\zeta$ depends on $\delta$, $\Upsilon$, $\rho$ and $s$. Thus, $\varepsilon_0$ depends on $\alpha$, $\beta$, $\lambda$, $\delta$, $\Upsilon$, $\rho$,  and $s$.  We recall that $0 < \zeta <1$, then~\eqref{F=0foundv} implies that first estimate in~\eqref{StimevGamma}.

Let $\varphi :\T^n \times \R^m \times J \to \T^n \times \R^m \times B \times J$ such that $\varphi(q,t) = (q, v(q,t))$ and $\Gamma(q,t) = b_0(q,t) + b_r(q,t) + \bar m \circ \tilde \varphi(q,t)v(q,t)$ for all $(q,t) \in \T^n \times \R^m \times J$.
Then, $\varphi$ is a $C^\rho$-weakly asymptotic cylinder associated to $(X_H, X_{h}, \varphi_0)$ with disturbing term $\Gamma$. We recall that $H$ is the Hamiltonian defined by~\eqref{HMD}, $h$ is the Hamiltonian in~\eqref{tildehBD} and $\varphi_0$ is the trivial embedding $\varphi_0 : \T^n \times \R^m  \to \T^n \times \R^m \times  B^{n+m}$, $\varphi_0(q) = (q,0)$. Moreover,  by~\eqref{HMD},~\eqref{F=0foundv} and the properties in Proposition \ref{normpropertiesMD}
\begin{eqnarray}
|\Gamma|_{\rho,1} &\le&  |b_0|_{\rho,1} +  |b_r|_{\rho,1}  + | \left(\bar m \circ \tilde \varphi \right) v|_{\rho, 1}\nonumber\\
&\le&|b_0|_{\rho,1} + C\left(|b_r|_{\lambda+1,1}  + | \left(\bar m \circ \tilde \varphi \right) v|_{\rho, 1}\right)\nonumber\\
&\le&  |b_0|_{\rho,1} + C\varepsilon + C(\rho)\Upsilon \zeta\nonumber\\
\label{GammaEst}
&\le&  |b_0|_{\rho,1} + K_\rho\Upsilon \zeta
\end{eqnarray}
for a suitable constant $K_\rho$ depending on $\rho$ and $n+m$.  We recall that we have chosen $\varepsilon_0 \le \zeta$.  Thanks to the latter and~\eqref{zeta2}, one has the second estimate in~\eqref{StimevGamma}.
 The second part of this section is dedicated to proving that the found $C^\rho$-weakly asymptotic cylinder $\varphi^t$ is Lagrangian (see Lemma \ref{LagrBDLemma} below). First, we observe that, thanks to~\eqref{GammaEst} and~\eqref{HMD}
\begin{eqnarray}
\label{StimeGamma1}
|\Gamma|_{1,1} &\le& |b_0|_{1,1} + K_1\Upsilon \zeta \le \delta + K_1\Upsilon \zeta.
\end{eqnarray}
Letting $\psi^t_{t_0, \bar \omega + \Gamma}$ be the flow at time $t$ with initial time $t_0$ of $\bar \omega + \Gamma$, we have the following quantitative lemma
\begin{lemma}
\label{phiMDfineGamma}
For all $t$, $t_0 \in J$,  if $t \ge t_0$ 
\begin{equation*}
|\partial_q \psi^t_{t_0, \bar \omega + \Gamma}|_{C^0} \le \left({t \over t_0} \right)^{C_2\left(\delta + \Upsilon \zeta\right)},
\end{equation*}
for a positive contant $C_2 \ge 1$ depending on $n+m$.
\end{lemma}
\begin{proof}
The proof of this lemma is similar to that of Lemma \ref{psi}. It relies on~\eqref{StimeGamma1} and the Gronwall inequality~\eqref{G2}. For this reason, it is omitted.
\end{proof}

Let $\psi^{t}_{t_0,H}$ be the flow at times $t$ with initial time $t_0$ of the Hamiltonian $H$ in~\eqref{HMD}, the following lemma concludes the proof of Theorem \ref{MD}.
\begin{lemma}
\label{LagrBDLemma}
$\varphi^{t_0}$ is Lagrangian for all $t_0 \in J$
\end{lemma}
\begin{proof}
Let $\alpha = dp \wedge dq$ be the standard symplectic form associated to $(q,p) \in \T^n \times \R^m \times B$. For all fixed $t$, $t_0 \in J$, the flow $\psi^t_{t_0,H}$  is a symplectomorphisms. This means that, for all fixed $t$, $t_0 \in J$, $(\psi_{t_0, H}^t)^*\alpha = \alpha$.  By~\eqref{hyp1WACbiss}, 
\begin{equation}
\label{hyp1bissbissMD}
\psi_{t_0,H}^{t_0 + t} \circ \varphi^{t_0} = \varphi^{t_0 + t} \circ \psi_{t_0, \bar \omega + \Gamma}^{t_0+t}
\end{equation}
and taking the pull-back with respect to the standard form $\alpha$ on both sides of the latter, we obtain
\begin{equation*}
(\varphi^{t_0})^*(\psi_{t_0,H}^{t_0 + t})^* \alpha =  ( \psi_{t_0, \bar \omega+ \Gamma}^{t_0+t})^*(\varphi^{t_0 + t})^* \alpha.
\end{equation*}
We know that  $\psi_{t_0, H}^{t_0 + t}$ is symplectic, then letting $(\psi_{t_0}^{t_0 + t})^* \alpha = \alpha$ on the left-hand side of the above equation, we have
\begin{equation*}
(\varphi^{t_0})^* \alpha =  ( \psi_{t_0, \bar \omega+ \Gamma}^{t_0+t})^*(\varphi^{t_0 + t})^* \alpha.
\end{equation*}
We want to prove that, for all $q\in \T^n\times \R^m$, $\left((\varphi^{t_0})^* \alpha\right)_q = 0$, where $\left((\varphi^{t_0})^* \alpha\right)_q$ stands for the symplectic form calculated on $q\in \T^n\times \R^m$. The idea is to prove that, for all $q\in \T^n\times \R^m$, the limit for $t \to +\infty$ on the right-hand side of the above equation converges to zero.  For this purpose, we recall that $\varphi^{t_0+t}(q) = (q, v^{t_0+t}(q))$,  and we observe that for all $q\in \T^n\times \R^m$
\begin{equation*}
\left((\psi_{t_0, \bar \omega+ \Gamma}^{t_0+t})^*(\varphi^{t_0 + t})^* \alpha\right)_q = \sum_{1 \le i <j \le n + m} \sum_{1 \le k < d \le n + m} \alpha^t_{i,j,k,d}(q) dq_k \wedge dq_d
\end{equation*}
where
\begin{eqnarray*}
\alpha^t_{i,j,k,d}(q)  &=& \left(\partial_{q_i}v^{t_0 + t}_j \circ \psi_{t_0, \bar \omega+ \Gamma}^{t_0+t}(q) - \partial_{q_j}v^{t_0 + t}_i \circ \psi_{t_0, \bar \omega+ \Gamma}^{t_0+t}(q) \right) \\
&\times&\left(\partial_{q_k}\psi_{t_0, \bar \omega+ \Gamma, i}^{t_0+t}(q)\partial_{q_d}\psi_{t_0, \bar \omega+ \Gamma,j}^{t_0+t}(q) - \partial_{q_d}\psi_{t_0, \bar \omega+ \Gamma, i}^{t_0+t}(q)\partial_{q_k}\psi_{t_0, \bar \omega+ \Gamma,j}^{t_0+t}(q)\right).
\end{eqnarray*}
In the latter $\times$ stands for the usual multiplication in $\R$. Then, for $t >0$ and fixed $1 \le i <j \le n + m$, $1 \le k < d \le n + m$, by Lemma \ref{phiMDfineGamma}
\begin{eqnarray*}
\left|\alpha^t_{i,j,k,d}\right|_{C^0}  &\le& \left|\partial_{q_i}v^{t_0 + t}_j \circ \psi_{t_0, \bar \omega+ \Gamma}^{t_0+t} - \partial_{q_j}v^{t_0 + t}_i \circ \psi_{t_0, \bar \omega+ \Gamma}^{t_0+t} \right|_{C^0} \\
&\times&\left|\partial_{q_k}\psi_{t_0, \bar \omega+ \Gamma, i}^{t_0+t}\partial_{q_d}\psi_{t_0, \bar \omega+ \Gamma,j}^{t_0+t} - \partial_{q_d}\psi_{t_0, \bar \omega+ \Gamma, i}^{t_0+t}\partial_{q_k}\psi_{t_0, \bar \omega+ \Gamma,j}^{t_0+t}\right|_{C^0}\\
&\le& \left(\left|\partial_{q_i}v^{t_0 + t}_j \right|_{C^0} +\left|\partial_{q_j}v^{t_0 + t}_i \right|_{C^0} \right)\\
&\times&\left(\left|\partial_{q_k}\psi_{t_0, \bar \omega+ \Gamma, i}^{t_0+t}\right|_{C^0} \left|\partial_{q_d}\psi_{t_0, \bar \omega+ \Gamma,j}^{t_0+t}\right|_{C^0} + \left|\partial_{q_d}\psi_{t_0, \bar \omega+ \Gamma, i}^{t_0+t}\right|_{C^0}\left|\partial_{q_k}\psi_{t_0, \bar \omega+ \Gamma,j}^{t_0+t}\right|_{C^0}\right)\\
&\le& C |v^{t_0+t}|_{C^1}\left(|\partial_q \psi^{t_0+t}_{t_0, \bar \omega + \Gamma}|_{C^0}\right)^2 \le C{\zeta \over t_0 + t} \left({t_0 + t \over t_0}\right)^{2C_2\left(\delta + \Upsilon \zeta\right)}.
\end{eqnarray*}
Thanks to~\eqref{zeta},  we have that $2C_2\left(\delta + \Upsilon \zeta\right) < 1$ and hence taking the limit for $t \to +\infty$ on both sides of the latter, the term in the last line converges to zero. This concludes the proof of this lemma.
\end{proof}

Now that the proof of Theorem \ref{MD} is complete, we want to conclude this section by comparing the result in~\cite{Sca22} with Theorem \ref{MD}.  Let $B^n \subset \R^n$ be an open ball centered at the origin and $\omega \in \R^n$.  We have the following
\begin{remark}\label{CompGMD}
In our previous work~\cite{Sca22}, we consider the following smooth time-dependent Hamiltonian 
\begin{equation*}
H:\T^n \times B^n \times J \to \R, \quad H(q,p,t) = \omega \cdot p + a(q,t) + b(q,t)\cdot p + m(q,p,t) \cdot p^2.
\end{equation*}
Given $l >2$, in this case, the perturbative terms $\partial_q a$ and $b$ decay as ${1 \over t^l}$ and ${1 \over t^{l-1}}$, respectively. 
We establish the existence of $u$, $v:\T^n \times J \to \R^n$ in such a way that, letting
\begin{equation}
\label{varphiGD}
\varphi:\T^n \times J \to \T^n \times B^n, \quad \varphi(q,t) = (q + u(q,t) , v(q,t)),
\end{equation}
$\varphi$ is a $C^\sigma$-asymptotic KAM torus associated to the previous Hamiltonian $H$. Additionally,  $u$ and $v$ decay as ${1 \over t^{l-2}}$ and ${1 \over t^{l-1}}$, respectively.  Roughly speaking, the component $u$ empowers us to control the dynamics at infinity. Hence, we prove the existence of orbit converging, as time tends to infinity, to the quasiperiodic solutions associated with the unperturbed system $h(q,p,t) =\omega \cdot p + m(q,p,t) \cdot p^2$.  However, the trade-off for controlling the dynamics at infinity involves a loss of two powers in the time decay of $\varphi^t$. 

 In the proof of the above-mentioned result~\cite{Sca22},  we consider a suitable functional $\mathcal{G}$ depending on the perturbative terms $(a,b)$ and the components of the $C^\sigma$-asymptotic KAM torus $(u,v)$ we are looking for.  Let $X_H$ be the Hamiltonian system associated with the Hamiltonian $H$. We define $\mathcal{G}$ in such a way that $\mathcal{G}(a,b,u,v)=0$ if and only if
 \begin{equation*}
 X_H (\varphi(q,t), t) -\partial_q \varphi(q,t) \omega - \partial_t \varphi(q,t) =0
 \end{equation*}
 where we recall that the latter is the invariant equation that a $C^\sigma$-asymptotic KAM torus has to satisfy (see Definition \ref{AsymKAMtorus}).  The functional $\mathcal{G}$ is defined on suitable Banach spaces of functions decaying polynomially fast in time. 
 The proof relies on the fixed point theorem. 

In Theorem \ref{MD},  we consider the limit case $l=2$.  The idea is to prove the existence of a family of embedded cylinders $\varphi:\T^n \times \R^m \times J \to \T^n \times \R^m\times B$ converging as time tends to infinity to the invariant cylinder $\varphi_0: \T^n \times \R^m \to \T^n \times \R^m\times B$, $\varphi_0(q) = (q,0)$ without controlling the dynamic at infinity.  In the definition of the functional $\mathcal{F}$ (see~\eqref{F}) the term 
\begin{equation*}
\partial_q v \left(b + \left(\bar m \circ \varphi\right)v\right)
\end{equation*}
induces a loss of $1$ derivative because of the presence of the derivative $\partial_q v$. We do not have this complication in~\cite{Sca22} because, in this case, the disturbing term $\Gamma$ is equal to zero (we recall that $\Gamma = b + \left(\bar m \circ \varphi\right)v$).  The presence of the term $\partial_q v$ in the definition of $\mathcal{F}$ leads us to define the Banach spaces $\mathcal{V}_{\sigma, \omega, 1}$ in such a way that the differential $D_v\mathcal{F}$ (see~\eqref{DiffF}) admits a right inverse of loss $1$.  For this reason,  we prove Theorem \ref{MD} using the Nash-Moser Theorem. 
\end{remark}

\section{Quasiperiodic motions in the planar three-body problem}\label{QPMP3BP}

This part is dedicated to a very brief introduction to the work of J. Féjoz~\cite{Fe02} concerning the existence, in a rotating frame of reference, of quasiperiodic motions with three frequencies for the Hamiltonian of the planar three-body problem. This result is an important element for the proof of Theorem \ref{Thmcomet}.

In this work, the author splits the dynamic into two parts: a fast, called Keplerian dynamic, and a slow, called secular dynamic. The first describes the motion of the bodies along three ellipses as if each body underwents the attraction of only one fictitious center of attraction. The slow dynamic describes how the mutual attraction of each planet deforms these Keplerian ellipses. There is a natural splitting
\begin{equation*}
H_0 = H_{\mathrm{Kep}} + H_{\mathrm{per}}
\end{equation*}
of the Hamiltonian when one uses the well-known Jacobi coordinates $\{(X_i, Y_i)\}_{i =0,1,2}$. Here, $H_{Kep}$ is the degenerate Hamiltonian of two decoupled two-body problems and $H_{per}$ is the perturbation.

The author defines the perturbing region contained in the direct product of the phase and parameter spaces. In this region, the Hamiltonian of the planar three-body problem is $C^k$-close  to the dynamically degenerate Hamiltonian of two decoupled two-body problems. To this end, we introduce some notations concerning the Keplerian dynamics. For the $i$th fictitious body, with $i=1$ or $2$, the mean longitude will be designated by $\lambda_i$, the semi-major axis by $a_i$, the eccentricity by $e_i$, the "centricity" $\sqrt{1 - e_i^2}$ by $\epsilon_i$, the argument of the pericenter by $g_i$, the mean motion by $\upsilon_i$ and the difference of the arguments of the pericenter by $g = g_1 - g_2$. We also introduce the well-known Poincaré coordinates $(\Lambda_i, \lambda_i, \xi_i, \eta_i)$, where we refer for example to the notes of A. Chenciner and J. Laskar~\cite{{Che89},{Las89}} or the work of J. Féjoz~\cite{Fe13}.

To measure how close the outer ellipse is from the inner ellipses when they are in opposition, the author defines
\begin{equation*}
\Delta = \max_{(\lambda_1, \lambda_2, g) \in \T^3}\max\{\sigma_0, \sigma_1\}{|X_1| \over |X_2|} = \max\{\sigma_0, \sigma_1\}{a_1(1 + e_1) \over a_2(1-e_2)}. 
\end{equation*}
He assumes that $\Delta <1$. This means that the outer ellipse does not meet the other two, whatever the difference $g$ of the arguments of the pericenters. Moreover, the eccentricity $e_2$ of the outer ellipse cannot be arbitrarily close to $1$. He also assumes that the eccentricity of the inner ellipses is upper bounded from $1$. 

Let $\mathfrak{P}$ be the reduced symplectically by translations phase space and $\mathfrak{M}$ be the space described by the three masses of the planets $m_0$, $m_1$ and $m_2$. 

\begin{definition}
\label{perregJ}
For a positive parameter $\delta$ and a non negative integer $k$, the perturbing region $\Pi^k_\delta$ of parameters $\delta$ and $k$ is the open subset of $\mathfrak{P} \times \mathfrak{M}$ defined by the following inequality
\begin{equation}
\label{PRFejoz}
\max\left\{{m_2 \over M_1} \left({a_1 \over a_2}\right)^{3 \over 2}, {\mu_1 \sqrt{M} \over M_1^{3 \over 2}} \left({a_1 \over a_2}\right)^2\right\}{1 \over \epsilon_2^{3(2+k)}(1 - \Delta)^{2k+1}}<\delta,
\end{equation}
where $M_1 = m_0 + m_1$, $M = m_0 + m_1 + m_2$ and $\mu_1 = {m_0 m_1 \over m_0 + m_1}$.
\end{definition}
Féjoz writes in his work that this inequality is not optimal and the given powers are not meaningful. He justifies this definition by proving that, inside the perturbing region, the perturbating function is $\delta$-small in a suitable $C^k$-norm. 
Concerning the case of Theorem \ref{Thmcomet}, where the masses are fixed, we point out that inequality~\eqref{PRFejoz} may be satisfied merely by assuming that ${a_1, \over a_2} \ll 1$ and $e_2 \le \mathrm{Cst}<1$. This is the so-called lunar or hierarchical regime.

In order to get rid of the degeneracy of the Keplerian Hamiltonian and hence apply the well-known KAM theorem, the secular Hamiltonian is introduced. Let $d$ and $k$ be suitable positive integers. On a suitable open set $\dot \Pi_\delta^k$ of $\Pi_\delta^k$, the author proves the existence of a $C^\infty$-symplectomorphism $\phi^d$ which is $\delta$-close to the identity in a suitable $C^k$-norm. The Hamiltonian  $H_0 \circ \phi^d$ can be split as follows
\begin{equation*}
H_0 \circ \phi^d = H_\pi^d + H^d_{\mathrm{comp}},
\end{equation*}
where $H^d_{\mathrm{comp}}$ is of size $\mathcal{O}(\delta^{d+1})$ on $\dot \Pi_\delta^{k + d(\tau + 4)}$. The secular Hamiltonian is $H_\pi^d$, which is Pöschel integrable~\cite{Po82}. It can be split into an integrable part $H^d_{\mathrm{int}}$ and a resonant part $H^d_{\mathrm{res}}$ of size $\mathcal{O}(\delta)$. The infinite jet of $H^d_{\mathrm{res}}$ vanishes along a suitable Cantor set. The previous splitting is obtained by an averaging process. It consists in averaging along the Keplerian ellipses parametrized by the mean anomalies $\lambda_1$ and $\lambda_2$ of the two fictitious Kepler problems where the Keplerian frequencies are non-resonant.  

After the reduction by the symmetry of rotation and far from elliptic singularities, the phase space of the secular Hamiltonian contains a positive measure of Lagrangian diophantine invariant tori.  The claim relies on a sophisticated version of KAM theorem, which is proved using a normal form theorem due to Herman. More specifically Féjoz proved the following theorem
\begin{theorems}
\label{Jacques}
In a rotating frame of reference, there are integers $k \ge 1$ and $d \ge 1$ and real numbers $\delta >0$ and $\tau \ge 1$ such that inside the perturbing region $\dot \Pi_\delta^{k + d(\tau + 4)}$ a positive measure of quasiperiodic Lagrangian tori of $H^d_\pi$ survive in the dynamics of the planar three-body problem.\footnote{Moreover, for any fixed masses, the set of quasiperiodic Lagrangian tori has positive Lebesgue measure.}
\end{theorems}

In the following section, we will use Theorem \ref{Jacques} as a black box in the proof of the existence of weakly asymptotically quasiperiodic solutions associated with the Hamiltonian of the planar three-body problem plus celestial body (see Theorem~\ref{Thmcomet}).


\section{Proof of Theorem \ref{Thmcomet}}\label{ProofTC}

We recall that we are considering three points of fixed masses $m_0$, $m_1$, and $m_2$ undergoing gravitational attraction in the plane.  These three points are perturbed by the passage of a celestial body with fixed mass $m_c$, whose motion is a given continuous function $c(t)$.  Only the planetary system is influenced by $c(t)$.  Moreover, we assume that 
\begin{equation*}
\lim_{t \to +\infty}|c(t)| = + \infty,  \quad \lim_{t \to +\infty}{d \over dt}|c(t)| = v_c>0.
\end{equation*}
Let $J = [1, +\infty)$ and $0<\varepsilon \le {1 \over 2}$ be the positive parameter in Theorem \ref{Thmcomet}. The phase space is the space
\begin{equation*}
\left\{((x_i,y_i) _{0 \le i\le 2}, t) \in \left(\R^2 \times \R^{2*}\right)^3 \times J \hspace{1mm}\Big|\hspace{1mm} \begin{matrix}\forall 0\le i < j \le 2, \hspace{1mm} x_i \ne x_j  \\
\forall 0\le i \le 2, \hspace{1mm}  {|x_i| \over |c(t)|} < \varepsilon \end{matrix}\right\}
\end{equation*}
where $(y_0, y_1, y_2)$ are the linear momentum covectors  and $(x_0, x_1, x_2)$ the position vectors of each body. The Hamiltonian in~\eqref{Hc}, that describes this system, is given by 
\begin{equation*}
H(x,y,t) = H_0(x,y) +H_c(x,t),
\end{equation*}
with
\begin{equation*}
H_0(x,y) = \sum_{i=0}^2 {|y_i|^2 \over 2m_i} - G\sum_{0 \le i<j\le 2}{m_i m_j \over |x_i - x_j|}, \quad H_c(x,t) = - G\sum_{i=0}^2 {m_i m_c \over |x_i - c(t)|}.
\end{equation*}
We may suppose $G$ equal to $1$.  We recall that $H_0$ is the Hamiltonian of the planar three-body problem, and $H_c$ is the Hamiltonian of the interaction between the planets and the celestial body $c(t)$. For $|c(1)|$ and $v_c$ sufficiently large and for $\varepsilon$ small enough, we will prove the existence of an open set of initial points giving rise to weakly asymptotically quasiperiodic solutions (see Definition \ref{asymsolI}) associated with the above Hamiltonian $H$. 

\subsection{Outline of the proof of Theorem \ref{Thmcomet}}

The proof of Theorem \ref{Thmcomet} is divided into five steps (we refer to Sections \ref{P1C}, \ref{SecQPDKC},  \ref{P3C}, \ref{P4C}, and \ref{P5C}). The first two parts (Sections \ref{P1C} and \ref{SecQPDKC},) are dedicated to the Hamiltonian of the planar three-body problem $H_0$.  In Section \ref{P1C},  we introduce a linear symplectic change of variable $\phi_0$ (see~\eqref{phi0}).  Letting $(X_i, Y_i)_{i=0,1,2}$ be the new variables, which should not be confused with the Jacobi coordinates introduced in Section \ref{QPMP3BP}, we can split the Hamiltonian $H_0$ in such a way that
\begin{equation*}
H_0 \circ  \phi_0(X, Y) = {|Y_0|^2 \over 2M} + K(X_1, X_2, Y_1, Y_2),
\end{equation*}
where $K$ is the Hamiltonian of the planar three-body problem after the symplectic reduction by the symmetry of translations, $X_0$ is the center of mass of the planetary system, $Y_0$ is the linear momentum of the planetary system, and $M = m_0 + m_1 + m_2$.  We stress that we introduce the above-mentioned variables in order to split the
dynamics into the absolute motion of the center of mass and the relative motion of the three bodies. 

For the second step of this proof (Section \ref{SecQPDKC}), we recall that Theorem \ref{Jacques} ensures the existence of Lagrangian $4$-dimensional invariant tori for the Hamiltonian $K$ in the phase space after the symplectic reduction by the symmetry of translations. As mentioned before, Theorem \ref{Jacques} proves the existence of quasiperiodic solutions with three frequencies for the Hamiltonian of the planar three-body problem in a rotating frame of reference. The additional frequency is given by the angular speed of the simultaneous rotations of the three ellipses.  
Let $B^n\subset \R^n$ be a $n$-dimensional ball centered at the origin with an unspecified radius that, conventionally, we will take greater than $1$.  In a neighborhood of one of the above-mentioned $4$-dimensional Lagrangian tori, we introduce a suitable symplectic change of coordinates $\phi$ (see~\eqref{phiC}) such that in these new symplectic variables, we can rewrite $H_0$ in the following form
\begin{align*}
& H_0 \circ \phi : \T^4 \times \R^2 \times B^4 \times B^2 \longrightarrow \R, \\
& H_0 \circ \phi (\theta, \xi, r, \eta) = e + \omega \cdot r + R_0(\theta, r) \cdot r^2 + {|\eta|^2 \over 2M}
\end{align*}
for some $e \in \R$, and $\omega \in \R^4$, whereas $R_0(\theta, r) \cdot r^2$ stands for the vector $r$ given twice as an argument of the symmetric bilinear form $R_0(\theta, r)$ (we refer to Lemma \ref{phiF}).  We stress that $\xi = X_0$ and $\eta = Y_0$. 

The third part of the proof (Section \ref{P3C}) is dedicated to the perturbing function $H_c$. We introduce a suitable subset $\mathcal{U}$ of the phase space $\T^4 \times \R^2 \times B^4 \times B^2 \times J$.  It is defined as the set of points that are sufficiently far from $c(t)$. We refer to~\eqref{U} for the definition of $\mathcal{U}$. For all $(\theta, \xi, r, \eta,t) \in \T^4 \times \R^2 \times B^4 \times B^2 \times J$, we define $\tilde \phi(\theta, \xi, r, \eta, t) = (\phi(\theta, \xi, r, \eta), t)$.  In this part, we will see that, for $|c(1)|$ and $v_c$ large enough,   
\begin{equation*}
H_c \circ \tilde \phi : \mathcal{U} \to \R
\end{equation*}
is well-defined (see Lemma \ref{Ux}), and it has a good decay in time (see Lemma \ref{LemmaHcstime}). 
Unfortunately,  we will see that the Hamiltonian $H \circ \tilde \phi = H_0 \circ \phi + H_c \circ \tilde \phi :\mathcal{U} \to \R$ does not satisfy the hypotheses of Theorem \ref{MD}.  
The reason is that $H \circ \tilde \phi$ is defined and satisfies good time-dependent estimates only on $\mathcal{U}$ and not on the entire phase space $\T^4 \times \R^2 \times B^4 \times B^2 \times J$. 

To solve this problem, in the fourth part of this proof (Section \ref{P4C}), we introduce a smooth extension \begin{equation*}
H_{ex} : \T^4 \times \R^2 \times B^4 \times B^2\times J \longrightarrow \R
\end{equation*}
of $H_c \circ \tilde \phi $ (see~\eqref{HexC}). This extension $H_{ex}$ coincides with $H_c \circ \tilde \phi$ on a suitable subset $\mathcal{U}_{1\over 2}$ of $\mathcal{U}$ where one expects the motions to take place.  The definition of $\mathcal{U}_{1\over 2}$ is provided by~\eqref{Uemezzo}.  In Lemma \ref{lemmaHex},  we establish that $H_{ex}$   exhibits the same decay in time as $H_c \circ \tilde \phi$ not only within $\mathcal{U}_{1\over 2}$ but on the entire phase space $ \T^4 \times \R^2 \times B^4 \times B^2\times J$.  
We conclude this part by verifying that the Hamiltonian 
\begin{equation*}
\hat H = H_0 \circ \phi + H_{ex}:  \T^4 \times \R^2 \times B^4 \times B^2\times J \longrightarrow \R
\end{equation*}
satisfies the hypotheses of Theorem \ref{MD} (see Lemma \ref{tildeHccarina}). Hence,  Theorem \ref{MD} provides the existence of a $C^1$-weakly asymptotic cylinder $\varphi^t$ for the Hamiltonian $\hat H$.  Let $\psi_{t_0, \hat H}^t$ be the flow at time $t$ with initial time $t_0$ of $\hat H$. Thanks to Proposition \ref{propasymsolI}, for all $z \in \varphi^1(\T^4 \times \R^2)$, it holds that $\psi_{1, \hat H}^t(z)$ is a weakly asymptotically quasiperiodic solution for $\hat H$.

In the final step of the proof  (Section \ref{P5C}), we prove the existence of an open subset of initial conditions giving rise to weakly asymptotically quasiperiodic solutions for the starting Hamiltonian $H$. For this purpose, we define the following set 
\begin{equation*}
B_1^2/2 = \left \{\xi \in \R^2 : |\xi| < {\varepsilon \over 6}|c(1)| \right\} \subset \R^2.
\end{equation*}
In Lemma \ref{Superlemmafinalecomet}, for $v_c$ large enough, we establish that for all $z \in \varphi^1 \left(\T^4 \times B_1^2/2\right)$,
\begin{equation*}
\psi_{1, \hat H}^t(z) \in \mathcal{U}_{1 \over 2}
\end{equation*}
for all $t \in J$.  The reason is that $c(t)$ diverges as  $|c(t)| \sim v_ct$, while the dynamic of the center of mass of the planetary system $\xi(t) = X_0(t)$ diverges as $|X_0(t)| \sim \ln t$.  Roughly speaking,  if the planes are far from $c(t)$ when $t=1$, then they remain far for all $t \in J$.  Thanks to the latter and the definition of $\hat H$,  we conclude that for all $z \in  \varphi^1 \left(\T^4 \times B_1^2/2\right)$,
\begin{equation*}
\psi_{1, \hat H}^t(z) = \psi_{1, H \circ \tilde \phi}^t(z)
\end{equation*}
for all $t \in J$.  Hence, using that $\phi$ is symplectic, we can verify that for all \newline$w \in \mathcal{W} = \phi \circ \varphi^1 \left(\T^4 \times B_1^2/2\right)$,
\begin{equation*}
\psi_{1,  H}^t(w) 
\end{equation*}
 is a weakly asymptotically quasiperiodic solution associated to $(X_H, X_{H_0}, \phi \circ \varphi_0)$.

\subsection{Splitting}\label{P1C}
The phase space and the Hamiltonian of the planar three-body problem are respectively
\begin{equation*}
\left\{(x_i,y_i) _{0 \le i\le 2} \in \left(\R^2 \times \R^{2*}\right)^3 \hspace{1mm}|\hspace{1mm} 0\le i < j \le 2, \hspace{1mm} x_i \ne x_j \right\},
\end{equation*}
and 
\begin{equation}
\label{H0Proof}
H_0(x,y) = \sum_{i=0}^2 {|y_i|^2 \over 2m_i} - \sum_{0 \le i<j\le 2}{m_i m_j \over |x_i - x_j|}.
\end{equation}

In this section, we would like to split the dynamics into the absolute motion of the center of mass and the relative motion of the three bodies. For this purpose, let us introduce the following linear symplectic change of coordinates 
\begin{equation}
\label{phi0}
\phi_0: (X_0,Y_0, X_1, Y_1, X_2, Y_2) \longrightarrow (x_0, y_0, x_1, y_1, x_2, y_2)
\end{equation}
such that
\begin{equation}
\label{XY}
\begin{cases}
X_0 = {m_0 \over M}x_0 + {m_1 \over M}x_1 + {m_2 \over M}x_2\\
X_1 = x_0 - x_1\\
X_2 = x_0 - x_2
\end{cases}
\begin{cases}
Y_0 = y_0 + y_1 + y_2\\
Y_1 = {m_1 \over M}y_0 - {m_0 + m_2\over M}y_1 + {m_1 \over M}y_2\\
Y_2 = {m_2 \over M}y_0 + {m_2 \over M}y_1 - {m_0 + m_1\over M}y_2 
\end{cases}
\end{equation}
where $M = m_0 + m_1 + m_2$. The left-hand side of the latter recalls the well-known heliocentric coordinates (see for example~\cite{Las89,Che89}). In these new variables, the Hamiltonian $H_0$ is split into two components: one dependent on $Y_0$ and another dependent on the variables $\{X_i, Y_i\}_{i = 1,2}$ (see~\eqref{H0phi0} below).  More specifically,  one can see that if $X_1 \ne 0$, $X_2 \ne 0$ and $X_1 \ne X_2$, the Hamiltonian $H_0 \circ \phi_0$ is equal to 
\begin{eqnarray}
\label{H0phi0}
H_0 \circ \phi_0(X,Y) &=& {|Y_0|^2 \over 2M}\\
&+& \underbrace{\left({|Y_1|^2 \over 2\mu_1} - {m_0 m_1 \over |X_1|}\right) + \left({|Y_2|^2 \over 2\mu_2} - {m_0 m_2 \over |X_2|}\right) + \left({Y_1 \cdot Y_2 \over m_0} - {m_1 m_2 \over |X_2 -X_1|}\right)}_{K(X_1,Y_1, X_2, Y_2)} \nonumber
\end{eqnarray}
where $\mu_1 = {m_0m_1 \over m_0 + m_1}$ and $\mu_2 = {m_0m_2 \over m_0 + m_2}$. We observe that $H_0\circ \phi_0$ is the sum of two independent Hamiltonians. The first ${|Y_0|^2 \over 2M}$ is responsible for the motion of the center of mass $X_0$ and the linear momentum $Y_0$, whereas $K$ is the Hamiltonian of the planar three-body problem after the reduction by the symmetry of translations.

\subsection{Quasiperiodic dynamics associated with $K$}\label{SecQPDKC}

In this section, we will use Theorem \ref{Jacques} in order to introduce symplectic coordinates and rewrite the Hamiltonian $K$, and hence $H_0$, in a more suitable form (see Lemma \ref{phiF} and Lemma \ref{H0carina}).  

For suitable integers $k \ge 1$, $d \ge 1$ and real numbers $\delta >0$ and $\tau \ge 1$, inside the perturbing region $\dot \Pi_\delta^{k + d(\tau + 4)}$ (see Definition \ref{perregJ}), Theorem \ref{Jacques} proves the existence of three-dimensional invariant tori for the Hamiltonian of the planar three-body problem $K$ in a rotating frame of reference. As mentioned above, these quasiperiodic motions have one additional frequency before the symplectic reduction by the symmetry of rotations. 
Here, we fix $m_0$, $m_1$ and $m_2$ as in Theorem \ref{Thmcomet} and we introduce the slice
 
\begin{equation*}
\dot \Pi_{\delta, m}^{k + d(\tau + 4)} = \dot \Pi_\delta^{k + d(\tau + 4)} \Big |_{m_0, m_1, m_2} \subset \mathfrak{P}.
\end{equation*}
We recall that $\mathfrak{P}$ is the phase space after the symplectic reduction by the symmetry by translations. In other words, $\dot \Pi_{\delta, m}^{k + d(\tau + 4)}$ is the subset of $\mathfrak{P}$ obtained by $\dot \Pi_\delta^{k + d(\tau + 4)}$ once we have fixed the masses $m_0$, $m_1$ and $m_2$. We refer to Section \ref{QPMP3BP} for the above notation.  
The following lemma provides suitable symplectic coordinates for the Hamiltonian of the planar three-body problem $K$ with the frame of reference attached to the center of mass of the planetary system.  We recall that $M_n$ is the set of $n$-dimensional matrices. Additionally,  $B^n \subset \R^n$ stands for a $n$-dimensional ball centered at the origin with an unspecified radius that, conventionally, we will take greater than $1$. 

\begin{lemma}
\label{phiF}
There exists a $C^\infty$ symplectic transformation 
\begin{equation*}
\phi_1 :(\theta, r)\longrightarrow(X_1, Y_1, X_2, Y_2)
\end{equation*}
defined on $\T^4 \times B^4$ with values in $\mathfrak{P}$ such that 
\begin{equation*}
K \circ \phi_1: \T^4 \times B^4 \to \R, \quad K \circ \phi_1(\theta, r) = e + \omega \cdot r + R_0(\theta, r) \cdot r^2
\end{equation*}
for some $e \in \R$, $\omega \in \R^4$ and $R_0:\T^4 \times B^4  \to M_4$.
\end{lemma}
\begin{proof}
By Theorem \ref{Jacques}, there exists a $4$-dimensional Lagrangian invariant torus $\mathcal{T} \subset \dot \Pi_{\delta, m}^{k + d(\tau + 4)}$ for $K$ supporting quasiperiodic dynamics. These tori form a set of positive Lebesgue measure, but we use only one such torus.

Due to the Weinstein Lagrangian neighbourhood theorem (see e.g. McDuff-Salamon~\cite{MDS}), there exists a neighbourhood $N(\mathcal{T})$ of $\mathcal{T}$ and a $C^\infty$ symplectomorphism
\begin{equation*}
\phi_1 : \T^4 \times B^4 \longrightarrow N(\mathcal{T}) \hspace{5mm} \mbox{such that} \hspace{5mm} \phi_1 (\T^4 \times \{0\})= \mathcal{T}.
\end{equation*}
We observe that $\phi_1(\T^4 \times \{0\}) = \mathcal{T}$ is a Lagrangian invariant torus for $K$. Hence
\begin{equation*}
K \circ \phi_1(\theta, 0) = c
\end{equation*}
for all $\theta \in \T^4$ and a suitable constant $c \in \R$. Moreover, $\phi_1(\T^4 \times \{0\}) = \mathcal{T}$ support a quasiperiodic dynamics with some frequency vector $\omega \in \R^4$, that is
\begin{equation*}
\partial_r\left(K \circ \phi_1\right)(\theta,0)=\omega
\end{equation*}
for all $\theta \in \T^4$.
\end{proof}

In order to obtain suitable symplectic coordinates for the Hamiltonian $H_0$ of the planar three-body problem, we need to lift the symplectic change of variables $\phi_1$ introduced in the previous lemma. 
For this purpose, we consider the following symplectic transformation 
\begin{equation*}
\bar \phi_1 : (\theta, \xi, r, \eta) \longrightarrow (X_0,Y_0, X_1, Y_1, X_2, Y_2)
\end{equation*}
defined on $\T^4 \times \R^2 \times B^4 \times B^2$ such that for all $(\theta, \xi, r, \eta) \in \T^4 \times \R^2 \times B^4 \times B^2$
\begin{equation}
\label{tildephiF}
\bar \phi_1(\theta, \xi, r, \eta) = (\xi, \eta, \phi_1(\theta, r))
\end{equation}
with $\xi = X_0$ and $\eta = Y_0$. We recall that $(X_0, Y_0)$ are the center of mass and the linear momentum of the planetary system introduced in the previous section (see~\eqref{XY}).  Furthermore,  we define
\begin{equation}
\label{phiC}
\phi =  \phi_0 \circ \bar \phi_1: (\theta, \xi, r, \eta) \longrightarrow (x_0,y_0, x_1, y_1, x_2, y_2)
\end{equation}
where $\phi_0$ is the linear symplectic transformation defined by~\eqref{phi0}.  We claim that, for all $k \in \Z$ with $k \ge 0$, $|D\phi|_{C^k}<\infty$, where $D\phi$ stands for the differential of $\phi$. This is because $\phi_0$ is a linear map, $\bar \phi_1$ is the identity with respect to $(\xi, \eta)$, and $\phi$ is $C^\infty$ and $1$-periodic with respect to $\theta_j$ for all $0 \le j \le 4$.  For this reason, for all $k \in \Z$ with $k \ge 0$, we define 
\begin{equation}
\label{UpsilonPhi}
\Upsilon^k_{\phi} = |D\phi|_{C^k}.
\end{equation}
In the following lemma, we rewrite $H_0$ in a more suitable form using the symplectic variables $(\theta, \xi, r, \eta)$.

\begin{lemma}
\label{H0carina}
We can write the Hamiltonian $H_0 \circ \phi$ in the following form
\begin{align*}
\label{H0bella}
& H_0 \circ \phi : \T^4 \times \R^2 \times B^4 \times B^2 \longrightarrow \R, \nonumber\\
& H_0 \circ \phi (\theta, \xi, r, \eta) = e + \omega \cdot r + R_0(\theta, r) \cdot r^2 + {|\eta|^2 \over 2M}
\end{align*}
with $e$, $\omega$ and $R_0$ as in Lemma \ref{phiF}.  
Moreover,  $H_0 \circ \phi \in C^\infty(\T^4 \times \R^2 \times B^4 \times B^2)$ and for all $z \in \Z$ with $k \ge 0$
\begin{equation}
\label{Upsilon2}
\Upsilon^k_{1,\phi} = |\partial^2_{(r,\eta)} \left(H_0 \circ \phi \right)|_{C^{k+1}} < \infty,
\end{equation}
where $\partial^2_{(r,\eta)} $ stands for the partial derivatives of order $2$ with respect to the variables $(r, \eta)$.
\end{lemma}
\begin{proof}
The first part of this lemma is a consequence of~\eqref{H0phi0},~\eqref{phiC}, and Lemma \ref{phiF}.  Concerning the second part, we observe that $H_0 \circ \phi \in C^\infty(\T^4 \times \R^2 \times B^4 \times B^2)$ because $\phi \in C^\infty(\T^4 \times \R^2 \times B^4 \times B^2)$ and takes values in a subset of the phase space where $H_0$ is also $C^\infty$.  The condition~\eqref{Upsilon2} is satisfied because $H_0 \circ \phi$, and consequently $\partial^2_{(r,\eta)} \left(H_0 \circ \phi \right)$, does not depend on the variable $\xi \in \R^2$, and $\partial^2_{(r,\eta)} \left(H_0 \circ \phi \right)$ is a $C^\infty$ function on the domain $\T^4 \times B^4 \times B^2$ and $1$-periodic in $\theta_j$ for all $1 \le j \le 4$. 
\end{proof}

\subsection{Perturbing function}\label{P3C}

This section deals with the study of the time-dependent perturbation $H_c$ (see~\eqref{Hc}).  We will introduce a suitable neighborhood $\mathcal{U}$ of $\T^4 \times \{0\} \times B^4 \times B^2 \times J \subset \T^4 \times \R^2 \times B^4 \times B^2 \times J$ characterized by the points that, at each time $t \in J$, are sufficiently far from $c(t)$ (see~\eqref{U} below).  Firstly, we will establish that $H_c$ is well-defined on $\mathcal{U}$ (see Lemma \ref{Ux}). It imposes a first restriction on $|c(1)|$.  Finally, we will verify that on $\mathcal{U}$, the perturbation $H_c$ satisfies good time-dependent estimates (see Lemma \ref{LemmaHcstime}). This verification necessitates an additional constraint on the parameter $v_c$,  and a second restriction on $|c(1)|$.

Let us be more precise and introduce that above-mentioned neighborhood $\mathcal{U}$ of $\T^4 \times \{0\} \times B^4 \times B^2 \times J \subset \T^4 \times \R^2 \times B^4 \times B^2 \times J$. For this purpose,  for all fixed $t \in J$, we define 
\begin{equation}
\label{Bt2}
B^2_t = \left\{ \xi \in \R^2 : {|\xi| \over |c(t)|} < {\varepsilon \over 3}\right\},
\end{equation}
where $\varepsilon$ is the positive parameter in Theorem \ref{Thmcomet}.
Let $\mathcal{U}$ be the following subset of $\T^4 \times \R^2 \times B^4 \times B^2 \times J$,
\begin{equation}
\label{U}
\mathcal{U} =  \bigcup_{t \in J}\left(\T^4 \times B^2_t \times B^4 \times B^2 \times \{t\}\right).
\end{equation}
In order to rewrite the perturbation $H_c$ in terms of the symplectic variables $(\theta, \xi, r, \eta)$, we define the following transformation $\tilde \phi$ such that
\begin{equation}
\label{tildephi}
\tilde \phi (\theta, \xi, r, \eta, t) = (\phi(\theta, \xi, r, \eta), t),
\end{equation}
for all $ (\theta, \xi, r, \eta, t) \in \T^4 \times \R^2 \times B^4 \times B^2 \times J$, where $\phi$ is the symplectic transformation defined by~\eqref{phiC}. In the following lemma, we verify that  $H_c \circ \tilde \phi : \mathcal{U} \to \R$ is well defined.  Specifically, our aim is to establish that $\tilde \phi \left(\mathcal{U}\right)$ is contained in the phase space~\eqref{xyC}.  This involves proving that 
\begin{equation*}
{|x_i(\theta, \xi, r)|  \over |c(t)|} < \varepsilon 
\end{equation*}
for all $t \in J$, $(\theta, \xi, r, \eta) \in \T^4 \times B^2_t \times B^4 \times B^2$, and $i=0,1,2$. We point out that thanks to~\eqref{XY} and Lemma \ref{phiF},  the variable $x_i$ depends on $(\theta, \xi, r) \in \T^4 \times \R^2 \times B^4$ for all $i=0, 1,2$. First,  some definitions and considerations.  We define the constant 
\begin{equation}
\label{k1}
\Upsilon_{2, \phi} = \max_{i=1,2}\sup_{(\theta, r) \in \T^4 \times B^4}|X_i(\theta, r)|.
\end{equation}
We observe that, thanks to Lemma \ref{phiF}, the variables $X_i$ with $i=1$,$2$ depend on the symplectic coordinates $(\theta, r) \in \T^4 \times B^4$.  Furthermore,  $\Upsilon_{2, \phi} < \infty$ because, for all $i = 1,2$, $X_i$ is a continuous function of $(\theta, r) \in \T^4 \times B^4$ and is $1$-periodic with respect to $\theta_j$ for all $0 \le j \le 4$. We also recall that the motion of the celestial body is a given continuous function $c(t)$ satisfying
\begin{equation*}
\lim_{t \to +\infty}|c(t)| = + \infty,  \quad \lim_{t \to +\infty}{d \over dt}|c(t)| = v_c>0.
\end{equation*}
 We observe that the latter implies the existence of $t_0 \gg 0$ such that
\begin{equation}
\label{ddtcv}
{v_c \over 2} \le {d \over dt}|c(t)| \le 2 v_c
\end{equation}
for all $t \ge t_0$. Without loss of generality, we can set $t_0 = 1$ by replacing $t$ with $t + t_0 - 1$.  Moreover,  by~\eqref{ddtcv} and the fundamental theorem of calculus 
\begin{equation}
\label{ccondition}
|c(1)| +{v_c \over 2} (t-1) \le |c(t)| \le |c(1)| + 2 v_c(t-1)
\end{equation}
for all $t \ge 1$.  Now, we have everything we need to prove the following
\begin{lemma}
\label{Ux}
We assume that
\begin{equation}
\label{ccomet1}
|c(1)| >{3\Upsilon_{2,\phi} \over \varepsilon},
\end{equation}
where $\Upsilon_{2, \phi}$ is the constant in~\eqref{k1} and $\varepsilon$ the positive parameter introduced in Theorem \ref{Thmcomet}.  Then,  for all $i = 0,1,2$
\begin{equation*}
{|x_i(\theta, \xi, r)|  \over |c(t)|} < \varepsilon 
\end{equation*}
for all $t \in J$ and  $(\theta, \xi, r, \eta) \in \T^4 \times B^2_t \times B^4 \times B^2$.
\end{lemma}
\begin{proof}
Because of~\eqref{XY},~\eqref{tildephiF}, and Lemma \ref{phiF}, we can rewrite the cartesian coordinates $(x_0, x_1, x_2)$ as follows
\begin{equation}
\label{xXthetaxir}
\begin{cases}
x_0(\theta, \xi, r) = \xi + {m_1 \over M}X_1(\theta, r) + {m_2 \over M}X_2(\theta, r)\\
x_1(\theta, \xi, r) = \xi - {m_0 + m_2 \over M}X_1(\theta, r) + {m_2 \over M}X_2(\theta, r)\\
x_2(\theta, \xi, r)= \xi + {m_1 \over M}X_1(\theta, r) - {m_0 + m_1 \over M}X_2(\theta, r),
\end{cases}
\end{equation}
for all $t \in J$ and  $(\theta, \xi, r, \eta) \in \T^4 \times B^2_t \times B^4 \times B^2$. By the latter,~\eqref{Bt2},~\eqref{U}, ~\eqref{k1},~\eqref{ccondition},  and~\eqref{ccomet1}
\begin{eqnarray*}
{|x_0(\theta, \xi, r) | \over |c(t)|} &\le& {\left| \xi + {m_1 \over M}X_1(\theta, r) + {m_2 \over M}X_2(\theta, r) \right|  \over |c(t)|}\\
&\le& {\left|\xi \right| \over |c(t)|} + {m_1 \over M}{\left|X_1(\theta, r)\right| \over |c(t)|} + {m_2 \over M}{\left|X_2(\theta, r)\right| \over |c(t)|}\\
&\le& {|\xi| \over |c(t)|} + {|X_1(\theta, r)| \over |c(1)|} + {|X_2(\theta, r)| \over |c(1)|} \\
&\le& {\varepsilon \over 3} + {2 \Upsilon_{2, \phi} \over |c(1)|} < \varepsilon
\end{eqnarray*}
for all $t \in J$ and $(\theta, \xi, r, \eta) \in \T^4 \times B^2_t \times B^4 \times B^2$. We recall that $M = m_0 + m_1 + m_2$, which implies ${m_i \over M} \le 1$ for all $0 \le i \le 2$.  We stress that in the third inequality of the latter, we used the fact that~\eqref{ccondition} implies $|c(t)| \ge |c(1)| +{v \over 2} (t-1) \ge |c(1)|$ for all $t \ge 1$. Similarly, the same claim holds for $x_1(\theta, \xi, r)$ and $x_2(\theta, \xi, r)$.
\end{proof}

The previous lemma ensures that $H_c \circ \tilde \phi$ is well defined on $\mathcal{U}$. Moreover, by~\eqref{xXthetaxir}, it is straightforward to verify that $H_c \circ \tilde \phi$ does not depend on the variable $\eta$.
In the second part of this section, we want to provide good estimates for $H_c \circ \tilde \phi$.  First, we need to prove the following

\begin{proposition}
\label{PcvComet}
Given $\varepsilon>0$,  if 
\begin{equation}
\label{propcv}
|c(1)| > { 1 \over \varepsilon}, \quad v_c > {2 \over \varepsilon},
\end{equation}
then
\begin{equation}
\label{perthyp2}
\sup_{t \ge 1}{t \over |c(t)|} < \varepsilon.
\end{equation}
\end{proposition}
\begin{proof}
By~\eqref{ccondition},  we can estimate ${t \over |c(t)|}$ as follows
\begin{equation*}
{t \over |c(t)|} \le {1 + (t -1) \over |c(1)| + {v_c \over 2}(t-1)},
\end{equation*}
for all $t \ge 1$. 
Thanks to~\eqref{propcv}, the conclusion~\eqref{perthyp2} is verified for $t=1$.  Now, we want to prove that the right-hand side of the latter is smaller than $\varepsilon$ also for all $t>1$. To this end, we suppose that there exists $t_0 >1$ such that 
\begin{equation*}
{1 + (t_0 -1) \over |c(1)| + {v_c \over 2}(t_0-1)} \ge \varepsilon.
\end{equation*}
We can rewrite the latter in the following form
\begin{equation*}
1 - \varepsilon |c(1)| \ge \left({\varepsilon \over 2}v_c -1 \right)(t_0 -1)
\end{equation*}
and this is a contradiction because by~\eqref{propcv}, 
\begin{equation*}
1 - \varepsilon |c(1)| <0, \quad \mbox{and} \quad \left({\varepsilon \over 2} v_c -1 \right)(t_0 -1) >0.
\end{equation*}
This concludes the proof of this proposition. 
\end{proof}

Now, we have all the elements to provide suitable estimates for the time-dependent perturbation $H_c$.  We recall that we use the notation $C$ to denote constants that we do not want to keep track of.  These constants may assume different values at different places. If we need to specify that a constant depends on particular parameters, we use the notation $C(\cdot)$, placing the parameters in brackets.

\begin{lemma}
\label{LemmaHcstime}
We assume that
\begin{equation}
\label{cvHyp2}
|c(1)| > { \max\{1, 3\Upsilon_{2, \phi}\} \over \varepsilon}, \quad v_c > {2 \over \varepsilon}
\end{equation}
where $\varepsilon$ is the positive parameter in Theorem \ref{Thmcomet} and $\Upsilon_{2, \phi}$ the constant in~\eqref{k1}. Then, for all $k \in \Z$ with $k \ge 0$ 
\begin{eqnarray}
\label{stimaHC0}
\sup_{t \ge 1}\left|H_c^t\right|_{C^0}t &\le& CMm_c \varepsilon, \\
\label{stimaHC}
\sup_{t \ge 1}\left|\partial_x H_c^t\right|_{C^{k}}t^2 &\le& C(k)Mm_c \varepsilon.
\end{eqnarray}
For all $t \in J$, the above norms $|\cdot|_{C^k}$ are taken on the domain $\phi\left(\T^4 \times B^2_t \times B^2 \times B^4\right)$, where $\phi$ is the symplectic transformation defined by~\eqref{phiC}.
\end{lemma}
\begin{proof}
For all $t \in J$ and  $(\theta, \xi, r,\eta) \in \T^4 \times B^2_t \times B^2 \times B^4$, let us rewrite the Hamiltonian $H_c \circ \tilde \phi$ in the following form
\begin{equation*}
\label{lemmaHvecchiaC}
H_c \circ \tilde \phi(\theta, \xi, r, \eta, t) = -\sum_{i=0}^2 {m_i m_c \over |x_i(\theta, \xi, r) - c(t)|}.
\end{equation*}
For the rest of this proof, $x_i = x_i(\theta, \xi,r)$ for all $0 \le i \le 2$. We drop the coordinates $(\theta, \xi,r)$ to obtain a more elegant form. 

Using Legendre polynomials,  for all $0 \le i\le2$,
\begin{eqnarray}
\label{LPC}
{1 \over \left|x_i - c(t) \right|} &=& {1 \over \left|c(t) \right|} \sum_{n\ge 0}\mathcal{P}_n(\cos \widehat{x_i c(t)}) \left({|x_i|\over|c(t)|}\right)^n\nonumber \\
&=&{1 \over \left|c(t) \right|}\left(1 + \sum_{n\ge 1}\mathcal{P}_n(\cos \widehat{x_i c(t)}) \left({|x_i|\over|c(t)|}\right)^n   \right)
\end{eqnarray}
where $\widehat{x_i c(t)}$ stands for the angle between the vectors $x_i$ and $c(t)$, whereas $\mathcal{P}_n$ is the $n$th Legendre polynomial.  We point out that the above expansion holds if ${|x_i|\over|c(t)|} < 1$, and this prerequisite is verified thanks to~\eqref{cvHyp2} and Lemma \ref{Ux}.

Now, thanks to~\eqref{cvHyp2},~\eqref{LPC}, Lemma \ref{Ux} and Proposition \ref{PcvComet}, for all $t \in J$ and $(x,y) =(x_0, x_1,x_2,y_0,y_1,y_2) \in \phi(\T^4 \times B^2_t \times B^4 \times B^2)$, 
\begin{eqnarray*}
|H^t_c(x)|t &=& \sum_{i=0}^2 {m_i m_ct \over |x_i - c(t)|} \\
&=& \sum_{i=0}^2{m_i m_ct \over |c(t)|}\left(1 + \sum_{n\ge 1}\mathcal{P}_n(\cos \widehat{x_i c(t)}) \left({|x_i|\over|c(t)|}\right)^n   \right)\\
&\le&\sum_{i=0}^2{m_i m_c t\over |c(t)|}\left(1 + \sum_{n \ge 1}\varepsilon^n\right) \le CMm_c \varepsilon.
\end{eqnarray*}
In the last line of the latter,  we used the fact that $\mathcal{P}_n(\cos \widehat{x_i c(t)}) \le 1$. Furthermore, given that $0<\varepsilon\le {1 \over 2}$, we can estimate $ \left(1 + \sum_{n \ge 1}\varepsilon^n\right)$ with a suitable constant $C$.  Whereas, by Proposition \ref{PcvComet},  ${t\over |c(t)|} \le \varepsilon$ for all $t \ge 1$.
Taking the sup over $\phi(\T^4 \times B^2_t \times B^4 \times B^2)$ and then for all $t \in J$ on the left-hand side of the latter, we obtain
\begin{equation*}
\sup_{t \ge 1}\left|H_c^t\right|_{C^0}t \le CMm_c \varepsilon.
\end{equation*}
This proves~\eqref{stimaHC0}. It remains to verify~\eqref{stimaHC}. For all $0 \le i \le 2$ and $k \ge 0$,  we denote by $\partial^{k+1}_{x_i}$ the derivative with respect to the variable $x_i$ of order $k+1$.  For all $t \in J$ and $(x,y) =(x_0, x_1,x_2,y_0,y_1,y_2) \in \phi(\T^4 \times B^2_t \times B^4 \times B^2)$,  we use~\eqref{cvHyp2},~\eqref{LPC}, Lemma \ref{Ux} and Proposition \ref{PcvComet}, to establish, similar to the previous case, that
\begin{eqnarray*}
|\partial^{k+1}_{x_i} H_c(x)|t^2 &\le& C(k){m_i m_c t^2 \over |x_i - c(t)|^{k+2}}\\
&=& C(k)m_i m_c{t^2 \over |c(t)|^{k+2}} \left( 1 + \sum_{n\ge 1}\mathcal{P}_n(\cos \widehat{x_i c(t)}) \left({|x_i|\over|c(t)|}\right)^n\right)^{k+2}\\
&\le& C(k)m_i m_c{t^{k+2} \over |c(t)|^{k+2}} \left( 1 + \sum_{n\ge 1}\mathcal{P}_n(\cos \widehat{x_i c(t)}) \left({|x_i|\over|c(t)|}\right)^n\right)^{k+2}\\
&\le& C(k)\varepsilon^{k+2} m_i m_c  \left(1 + \sum_{n \ge 1}\varepsilon^n\right)^{k+2} \le \varepsilon^{k+2} C(k) m_i m_c,\\
&\le& \varepsilon^{k+2} C(k) M m_c,
\end{eqnarray*}
where $t \ge 1$ and $k+2 \ge 2$ imply $t^2 \le t^{k+2}$ in the third line.  Also in this case, because of $0<\varepsilon\le{1\over2}$, we can estimate $ \left(1 + \sum_{n \ge 1}\varepsilon^n\right)^{k+2}$ by a suitable constant depending on $k$. Taking the max for all $i=0,1,2$ on the left-hand side of the latter we obtain
\begin{equation*}
\sup_{t \ge 1}\left|\partial^{k+1}_x H_c^t\right|_{C^0}t^2 \le C(k)Mm_c \varepsilon^{k+2}.
\end{equation*}
Now,  thanks to the above estimate and remembering the definition of the $C^k$ norm (see~\eqref{Cknorm} in Appendix \ref{A}),  one can conclude the proof of~\eqref{stimaHC}.
\end{proof}

\subsection{Smooth extension of the perturbing function}\label{P4C}

The previous section established that the time-dependent perturbation $H_c\circ \tilde \phi$ satisfies good estimates exclusively within $\mathcal{U}$ rather than in the entire phase space $\T^4 \times \R^2 \times B^4 \times B^2 \times J$. In order to apply Theorem \ref{MD} and prove the existence of weakly asymptotically quasiperiodic solutions for the Hamiltonian $H$ in~\eqref{Hc}, we need to introduce a suitable smooth extension $H_{ex}$ of $H_c \circ \tilde \phi$ defined on $\T^4 \times \R^2 \times B^4 \times B^2 \times J$ satisfying the same estimates in the entire phase space. 

This section is divided into two parts. First, we prove the existence of the above-mentioned smooth extension $H_{ex}$ (see Lemma \ref{lemmaHex}). In the second part, we verify that the new Hamiltonian $\hat H = H_0 \circ \phi + H_{ex}$ satisfies the hypotheses of Theorem \ref{MD} (see Lemma \ref{tildeHccarina}). 
To this end, we need to introduce a suitable subset $\mathcal{U}_{1 \over 2}$ of $\mathcal{U}$ where we expect the motions to take place. For all fixed $t \in J$, we define 
\begin{equation}
\label{Bt2emezzo}
B^2_t / 2 = \{ \xi \in \R^2 : {|\xi| \over |c(t)|} < {\varepsilon \over 6}\} \subset B^2_t,
\end{equation}
where $B^2_t$ is the set in~\eqref{Bt2} and $\varepsilon$ is the positive parameter in Theorem \ref{Thmcomet}.
Let $\mathcal{U}_{1 \over 2}$ be the following subset of $\mathcal{U}$,
\begin{equation}
\label{Uemezzo}
\mathcal{U}_{1 \over 2} =  \bigcup_{t \in J}\left(\T^4 \times \left(B^2_t/ 2\right)\times B^4 \times B^2 \times \{t\}\right) \subset \mathcal{U}.
\end{equation}
We denote by $C^\infty_0$ the set of the $C^\infty$ compactly supported functions. 
In order to define the above-mentioned smooth extension,  for all fixed $t \in J$, we consider the following family of functions
\begin{equation*}
\begin{cases}
\rho^t :\R^2 \longrightarrow B_t^2,\\
\rho^t \in C^\infty_0(\R^2),\\
\rho^t(\xi) = \xi \hspace{2mm} \mbox{for all $|\xi| \le \varepsilon {|c(t)| \over 6}$},\\
\rho^t(\xi) = 0 \hspace{2mm} \mbox{for all $|\xi| \ge \varepsilon {|c(t)| \over 3}$}.\\
\end{cases}
\end{equation*}
and we introduce the following map 
\begin{equation}
\begin{aligned}
\label{pi}
&\pi :\T^4 \times \R^2 \times B^4 \times B^2 \times J\to \T^4 \times \R^2 \times B^4 \times B^2\\
&\pi (\theta, \xi, r, \eta, t) = (\theta, \rho^t (\xi), r, \eta).
\end{aligned}
\end{equation}
Now, in agreement with the notation~\eqref{tildef} introduced in Section \ref{FS}, we denote by $\tilde \pi$ the following transformation
\begin{equation*}
\pi :\T^4 \times \R^2 \times B^4 \times B^2 \times J\to \T^4 \times \R^2 \times B^4 \times B^2\times J, \quad \tilde \pi (\theta, \xi, r, \eta, t) = (\pi(\theta, \xi, r, \eta, t),t).
\end{equation*}
It is straightforward to verify that
\begin{equation}
\label{PropPi}
\tilde \pi \left(\T^4 \times \R^2 \times B^4 \times B^2 \times J\right) \subset \mathcal{U}\hspace{5mm}\mbox{and}\hspace{5mm}\tilde \pi \big |_{\mathcal{U}_{1 \over 2}} = \mathrm{Id},
\end{equation}
where $\tilde \pi \big |_{\mathcal{U}_{1 \over 2}} = \mathrm{Id}$ means that $\tilde \pi (\theta, \xi, r, \eta, t)= (\theta, \xi, r, \eta, t)$ for all $(\theta, \xi, r, \eta, t) \in \mathcal{U}_{1 \over 2}$.  Moreover, if we assume that $|c(1)| < {1 \over \varepsilon}$, then for all fixed $t \in J$ and $k \in \Z$ with $k \ge 0$
\begin{equation}
\label{rhot}
|\rho^t|_{C^k} \le {C(k) \over \left(\varepsilon|c(t)|\right)^k} \le C(k).
\end{equation}
We define
\begin{equation}
\label{UpsilonPi}
\Upsilon_\pi^k=\sup_{t \in J}|D \pi^t|_{C^k}, 
\end{equation}
where $D\pi^t$ stands for the differential of $\pi^t$ with respect to the variables $(\theta, \xi, r, \eta)$ (we refer to~\eqref{ft} for the notation $\pi^t$).  Thanks to~\eqref{rhot},  if $|c(1)| < {1 \over \varepsilon}$, then for all $k \in \Z$ with $k \ge 0$ one can prove that 
\begin{equation}
\label{UpsilonPiBounded}
\Upsilon_\pi^k \le C(k).
\end{equation}

We define the following Hamiltonian 
\begin{equation}
\label{HexC}
H_{ex} = H_c \circ \tilde \phi \circ \tilde\pi : \T^4 \times \R^2 \times B^4 \times B^2 \times J \longrightarrow \R.
\end{equation}
It is straightforward to verify that, for all fixed $t \in J$, $H^t_{ex} \in C^\infty(\T^4 \times \R^2 \times B^4 \times B^2)$. Moreover, by construction, $\partial^k_{(\theta, \xi, r, \eta)} H_{ex} \in C(\T^4 \times \R^2 \times B^4 \times B^2 \times J)$ for all $k \in \Z$ with $k \ge 0$.  We point our that $\partial^k_{(\theta, \xi, r, \eta)}$ denotes the partial derivatives of order $k$ with respect to the variables $(\theta, \xi, r, \eta)$. 
In the following lemma, we will see that~\eqref{HexC} smoothly extends $H_c \circ \tilde \phi$ to the entire phase space 
$\T^4 \times \R^2 \times B^4 \times B^2 \times J$.
\begin{lemma}
\label{lemmaHex}
Let $H_{ex} $ be as in~\eqref{HexC}. Then,  $H_{ex}$ does not depend on $\eta$ and 
\begin{equation}
\label{HexHcC}
 H_{ex} \big |_{\mathcal{U}_{1 \over 2}} = H_c \circ \tilde \phi.
\end{equation}
That is, $H_{ex} (\theta, \xi, r,\eta,t) =  H_c \circ \tilde \phi(\theta, \xi, r,\eta,t)$ for all $(\theta, \xi, r, \eta,t) \in \mathcal{U}_{1 \over 2}$.
Moreover, if
\begin{equation}
\label{cvHyp3}
|c(1)| > { \max\{1, 3\Upsilon_{2, \phi}\} \over \varepsilon}, \quad v_c > {2 \over \varepsilon},
\end{equation}
then for all $k \in \Z$ with $k \ge 0$, we have the following estimates
\begin{eqnarray}
\label{stimaHex1}
\sup_{t \ge 1}| H_{ex}^t|_{C^0}t &\le& CMm_c \varepsilon, \\
\label{stimaHex2}
\sup_{t \ge 1}|\partial_{(\theta,\xi, r)} H_{ex}^t|_{C^{k}}t^2 &\le& C(k, \Upsilon^k_\phi)Mm_c \varepsilon,
\end{eqnarray}
where $\partial_{(\theta,\xi, r)}$ stands for the partial derivatives of order $1$ with respect to the variables $(\theta, \xi, r)$.  We recall that the constants $\Upsilon_{2, \phi}$ and $\Upsilon^k_\phi$ are defined by~\eqref{k1} and~\eqref{UpsilonPhi}, respectively. 
\end{lemma}
\begin{proof}
Using the definition of  $\pi$, especially~\eqref{PropPi}, verifying that~\eqref{HexHcC} is satisfied is straightforward. Now, we assume~\eqref{cvHyp3} and we prove~\eqref{stimaHex1} and~\eqref{stimaHex2}. We define $B^6 = B^4 \times B^2 \subset \R^4 \times \R^2$ to shorten the notation.  For the sake of clarity,  in this proof, we will specify the domain where the Hölder norms are taken (see Appendix \ref{A}).  For all $t \in J$, we observe that by~\eqref{stimaHC0} and~\eqref{HexC} 
\begin{eqnarray*}
| H_{ex}^t|_{C^0(\T^4 \times \R^2 \times B^6)}t &=& \left|\left(H_c \circ \tilde \phi \circ \tilde\pi\right)^t\right|_{C^0(\T^4 \times \R^2 \times B^6)}t\\
&\le&  |H^t_c |_{C^0\left(\phi\left(\T^4 \times B_t^2 \times B^6\right)\right)}t \le C Mm_c\varepsilon.
\end{eqnarray*}
Taking the sup for all $t \in J$ on the left-hand side of the latter, we prove~\eqref{stimaHex1}. Concerning the second estimate~\eqref{stimaHex2},  for all $t \in J$ and $k \in \Z$ with $k \ge 0$
\begin{align}
&| \partial_{(\theta,\xi, r)}H_{ex}^t|_{C^k(\T^4 \times \R^2 \times B^6)} t^2= \left|\partial_{(\theta,\xi, r)}\left(H_c \circ \tilde \phi \circ \tilde \pi\right)^t\right|_{C^k(\T^4 \times \R^2 \times B^6)}t^2\nonumber\\
&= \left|\partial_{(\theta,\xi, r)}\left(H^t_c \circ \phi \circ \pi^t\right)\right|_{C^k(\T^4 \times \R^2 \times B^6)}t^2\nonumber\\
&\le C(k) |\partial_x H^t_c \circ \phi \circ \pi^t|_{C^k(\T^4 \times \R^2 \times B^6)}t^2 |D\phi|_{C^k(\T^4 \times \R^2 \times B^2_t \times B^4)}|D\pi^t|_{C^k(\T^4 \times \R^2 \times B^6)}\nonumber\\
\label{StimeHex2Partial}
&\le C(k, \Upsilon^k_\phi) |\partial_x H^t_c \circ \phi \circ \pi^t|_{C^k(\T^4 \times \R^2 \times B^6)}t^2.
\end{align}
We note that the first inequality (third line) of the latter is due to the chain rule and property \textit{2} of Proposition \ref{Holder}.  In the last line, we use~\eqref{UpsilonPhi},~\eqref{UpsilonPi} and~\eqref{UpsilonPiBounded}.  We want to point out that in this proof we are using the notation~\eqref{ft}.

We need to study the cases $k=0$ and $k \ge 1$ separately.  In both cases,~\eqref{StimeHex2Partial} will be our starting point. Hence, letting $k=0$, thanks to~\eqref{StimeHex2Partial} we have that
\begin{align*}
&| \partial_{(\theta,\xi, r)}H_{ex}^t|_{C^0(\T^4 \times \R^2 \times B^6)} t^2 \le C(\Upsilon^0_\phi) |\partial_x H^t_c \circ \phi \circ \pi^t|_{C^0(\T^4 \times \R^2 \times B^6)}t^2\\
&\le C(\Upsilon^0_\phi) |\partial_x H^t_c|_{C^0(\phi(\T^4 \times B^2_t \times B^6)}t^2 \le C(\Upsilon^0_\phi) Mm_c\varepsilon
\end{align*}
for all $t \in J$.  We note that the last inequality of the latter is a consequence of~\eqref{stimaHC}.  Taking the sup for all $t \in J$, we prove~\eqref{stimaHex2} in the case $k=0$.  It remains to verify~\eqref{stimaHex2} when $k \ge 1$.  Therefore, for all $k \in \Z$ with $k \ge 1$, by~\eqref{StimeHex2Partial} we observe that
\begin{align}
&| \partial_{(\theta,\xi, r)}H_{ex}^t|_{C^k(\T^4 \times \R^2 \times B^6)} t^2 \le C(k, \Upsilon^k_\phi) |\partial_x H^t_c \circ \phi \circ \pi^t|_{C^k(\T^4 \times \R^2 \times B^6)}t^2\nonumber\\
\label{StimeHex2Final1}
&\le C(k, \Upsilon^k_\phi) |\partial_x H^t_c \circ \phi|_{C^k(\T^4 \times B^2_t \times B^6)}t^2 \left(1 + \left(\Upsilon^0_\pi\right)^k + \Upsilon_\pi^{k-1}\right)\\
\label{StimeHex2Final2}
&\le C(k, \Upsilon^k_\phi) |\partial_x H^t_c|_{C^k(\phi(\T^4 \times B^2_t \times B^6))}t^2 \left(1 + \left(\Upsilon^0_\phi\right)^k + \Upsilon_\phi^{k-1}\right)\\
&\le C(k, \Upsilon^k_\phi) Mm_c\varepsilon
\end{align}
for all $t \in J$. At line~\eqref{StimeHex2Final1} we use property \textit{5} of Proposition \ref{Holder} and~\eqref{UpsilonPi}. Similarly, line~\eqref{StimeHex2Final2} is a consequence of property \textit{5} of Proposition \ref{Holder},~\eqref{UpsilonPhi} and~\eqref{UpsilonPiBounded}.
Whereas, in the last line of the latter we use~\eqref{stimaHC}. Taking the sup for all $t \in J$,  on both sides of the above inequality we conclude the proof of~\eqref{stimaHex2}. 
\end{proof}
The rest of this section is dedicated to showing that the Hamiltonian
\begin{equation}
\label{HatH}
\hat H: \T^4 \times \R^2 \times B^4 \times B^2 \longrightarrow \R,\quad \hat H  = H_0 \circ  \phi +  H_{ex}
\end{equation}
satisfies the hypothesis of Theorem \ref{MD}, where $H_0$ is the Hamiltonian of the planar three-body problem~\eqref{H0Proof}.  In the following lemma, we will use the notation and the Banach spaces introduced in Section \ref{FS}.

\begin{lemma}
\label{tildeHccarina}
Let $\hat H$ be as in~\eqref{HatH}. We can rewrite $\hat H$ in the following form 
\begin{eqnarray}
\hat H(\theta, \xi, r,\eta,t) &=&  e +  \omega \cdot r + a(\theta, \xi,t) + b(\theta, \xi,t) \cdot r \nonumber\\
\label{HatH2}
&+& R(\theta, \xi,r,t)\cdot r^2 + {|\eta|^2 \over 2M}
\end{eqnarray}
for all $(\theta, \xi, r, \eta, t) \in \T^4 \times \R^2 \times B^4 \times B^2 \times J$. We note that $R(\theta, \xi,r,t)\cdot r^2$ stands for the vector $r$ given twice as the argument of the symmetric bilinear form $R$. Moreover,  if
\begin{equation*}
|c(1)| > { \max\{1, 3\Upsilon_{2, \phi}\} \over \varepsilon}, \quad v_c > {2 \over \varepsilon},
\end{equation*}
then, for all $k \in \Z$ with $k \ge 1$
\begin{equation}
\label{regabhatH}
a \in \mathcal{S}_{k,(0,2)}, \quad b \in \mathcal{S}_{k+1,1}, \quad \partial_{(r,\eta)}^2 \hat H \in \mathcal{S}_{k+1, 0}
\end{equation}
with 
\begin{eqnarray}
|a|_{k, (0,2)} &\le&C_1(k, \Upsilon^{k}_{\phi})Mm_c\varepsilon ,\nonumber\\
\label{StimeHccarina}
| b|_{k+1, 1} &\le& C_2(k,  \Upsilon^{k+1}_{\phi})Mm_c\varepsilon ,\\
\left|\partial_{(r,\eta)}^2 \hat H\right|_{k+1, 0} &\le& \Upsilon^{k+1}_{1, \phi} + C_3(k,  \Upsilon^{k+2}_{\phi})Mm_c\varepsilon.\nonumber
\end{eqnarray}
We point out that $C_1$, $C_2$, and $C_3$ are constants depending on the parameters in brackets that we need to keep track of. Whereas, $\Upsilon^{k+1}_{1, \phi}$,  and $\Upsilon^{k}_{\phi}$ are defined by~\eqref{Upsilon2} and~\eqref{UpsilonPhi}, respectively.
\end{lemma}
\begin{proof}
To establish the first part of this lemma,  we write $H_0 \circ \phi$ as in Lemma \ref{H0carina}, and we expand $H_{ex}$ in a small neighborhood of $r=0$,  in such a way that 
\begin{eqnarray*}
\hat H(\theta, \xi, r, \eta, t) &=& e + \omega \cdot r + R_0(\theta, r) \cdot r^2 + {|\eta|^2 \over 2M}\\
&+& H_{ex} (\theta, \xi,0,t) + \partial_r H_{ex} (\theta, \xi,0,t) \cdot r + \int_0^1 (1 -\tau)\partial^2_r H_{ex} (\theta, \xi,\tau r, t)d\tau \cdot r^2
\end{eqnarray*}
for all $(\theta, \xi, r, \eta, t) \in \T^4 \times \R^2 \times B^4 \times B^2 \times J$, where $e$, $\omega$ and $R_0$ are those in Lemma \ref{H0carina} and we recall that $H_{ex}$ does not depend on $\eta$ (see Lemma \ref{lemmaHex}).  We define 
\begin{eqnarray*}
 a(\theta, \xi,t) &=&  H_{ex} (\theta, \xi,0,t),  \\
 b(\theta, \xi,t) &=& \partial_r H_{ex} (\theta, \xi,0,t), \\
 R(\theta, \xi,r,t) &=& R_0(\theta, r)  + \int_0^1 (1 -\tau)\partial^2_r H_{ex} (\theta, \xi,\tau r, t)d\tau,
\end{eqnarray*}
for all $(\theta, \xi, r, \eta, t) \in \T^4 \times \R^2 \times B^4 \times B^2 \times J$, and this concludes the proof of~\eqref{HatH2}.  

In the proof of the second part of this lemma, we will widely use Lemma \ref{H0carina} and Lemma \ref{lemmaHex}, especially the bounds~\eqref{Upsilon2},~\eqref{stimaHex1}, and~\eqref{stimaHex2}. Hence, remembering the notation introduced in Section \ref{FS}, we observe that
\begin{eqnarray}
|a|_{k+1, (0,2)} &=& |a|_{k+1, 0} + |\partial_{(\theta, \xi)} a|_{k,2} \le \sup_{t\ge 1}|H_{ex}^t|_{C^{k+1}} + \sup_{t\ge 1}|\partial_{(\theta, \xi)}H_{ex}^t|_{C^{k}}t^2\nonumber\\
\label{HatHStime1}
&\le&\sup_{t\ge 1}|H_{ex}^t|_{C^{0}}t + 2\sup_{t\ge 1}|\partial_{(\theta, \xi, r)}H_{ex}^t|_{C^{k}}t^2\\
&\le& C(k, \Upsilon^{k}_{\phi})Mm_c\varepsilon \nonumber\\
\label{HatHStime2}
|b|_{k+1, 1} &\le& \sup_{t\ge 1}|\partial_{(\theta, \xi)}H_{ex}^t|_{C^{k+1}}t \le \sup_{t\ge 1}|\partial_{(\theta, \xi)}H_{ex}^t|_{C^{k+1}}t^2\\
&\le& C(k, \Upsilon^{k+1}_{\phi})Mm_c\varepsilon \nonumber\\
|\partial^2_{(r,\eta)} \hat H|_{k+1, 0} &=& \sup_{t\ge1} |\partial^2_{(r,\eta)} \left(H_0 \circ \phi\right) + \partial_{(r,\eta)}^2 H^t_{ex}|_{C^{k+1}}\nonumber\\
\label{HatHStime3}
&\le& |\partial^2_{(r,\eta)} \left(H_0 \circ \phi\right)|_{C^{k+1}} + \sup_{t\ge 1}|\partial_rH_{ex}^t|_{C^{k+2}}t^2\\
&\le& \Upsilon^{k+1}_{1, \phi} + C(k, \Upsilon^{k+2}_{\phi})Mm_c\varepsilon.\nonumber
\end{eqnarray}
We stress that in line~\eqref{HatHStime1} we used $$\sup_{t\ge 1}|H_{ex}^t|_{C^{k+1}} \le \sup_{t\ge 1}|H_{ex}^t|_{C^{0}}t + \sup_{t\ge 1}|\partial_{(\theta, \xi, r)}H_{ex}^t|_{C^{k}}t^2$$ and the trivial estimate $\sup_{t\ge 1}|\partial_{(\theta, \xi)}H_{ex}^t|_{C^{k}}t^2 \le \sup_{t\ge 1}|\partial_{(\theta, \xi, r)}H_{ex}^t|_{C^{k}}t^2$.  Whereas, in lines~\eqref{HatHStime2} and~\eqref{HatHStime3} we used the following estimates $$ \sup_{t\ge 1}|\partial_{(\theta, \xi)}H_{ex}^t|_{C^{k+1}}t \le \sup_{t\ge 1}|\partial_{(\theta, \xi)}H_{ex}^t|_{C^{k+1}}t^2, \quad \sup_{t\ge 1}|\partial_{(r,\eta)}^2 H^t_{ex}|_{C^{k+1}} \le \sup_{t\ge 1}|\partial_rH_{ex}^t|_{C^{k+2}}t^2$$
due to the definition of the norm $|\cdot|_{C^k}$ and to the fact that $t \ge 1$.  This concludes the proof of the estimates~\eqref{StimeHccarina}.  Concerning the regularity of $a$, $b$ and $ \partial_{(r,\eta)}^2 \hat H$, specifically as stated in~\eqref{regabhatH}, it is a straightforward consequence of~\eqref{StimeHccarina}, the regularity of $H_0 \circ \phi$ (see Lemma \ref{H0carina}) and the construction of $H_{ex}$ (see~\eqref{pi} and~\eqref{HexC}).
\end{proof}


The above Lemma is the main ingredient in establishing that the Hamiltonian $\hat H$ in~\eqref{HatH} satisfies the hypotheses of Theorem \ref{MD}.  Let $s$, $\beta$, $\alpha$, $\lambda$ and $\rho$ be the following parameters
\begin{equation*}
s = 8, \quad \beta= {31 \over 24}, \quad  \alpha = {7 \over 6}, \quad \lambda = {15 \over 4}, \quad \rho=1. 
\end{equation*}
The above parameters are obtained by~\eqref{parametriconsminimo} with $s=8$. 
One can see that $s$, $\beta$, $\alpha$, $\lambda$ and $\rho$ satisfy~\eqref{parametersMD}.  Now, we define the following constant 
\begin{equation}
\label{UpsilonPhiPi}
\Upsilon_{\phi,M, m_c} = \max\left\{C_1(\Upsilon^8_{\phi})Mm_c, C_2(\Upsilon^{9}_{\phi})Mm_c, \Upsilon^{9}_{1, \phi} + C_3(\Upsilon^{10}_{\phi})Mm_c\right\}
\end{equation}
where $C_1$, $C_2$ and $C_3$ are the constants in~\eqref{StimeHccarina}. We fix the parameters $\delta$ and $\Upsilon$ in~\eqref{HMD} as follows
\begin{equation*}
\delta = 0, \quad \Upsilon = \Upsilon_{\phi,M, m_c}
\end{equation*}
and we assume
\begin{equation*}
|c(1)| > { \max\{1, 3\Upsilon_{2, \phi}\} \over \varepsilon}, \quad v_c > {2 \over \varepsilon}. 
\end{equation*}
Furthermore, we define the following Hamiltonian
\begin{equation*}
\hat h :  \T^4 \times \R^2 \times B^4 \times B^2 \to \R,  \quad \hat h(\theta, \xi, r, \eta) = e +  \omega \cdot r + R(\theta, \xi,r,t)\cdot r^2 + {|\eta|^2 \over 2M}
\end{equation*}
and the trivial embedding 
\begin{equation*}
\varphi_0 :\T^4 \times \R^2 \to \T^4 \times \R^2 \times B^4 \times B^2 , \quad \varphi_0 (\theta, \xi) = (\theta, \xi, 0, 0).
\end{equation*}
Thanks to Lemma \ref{tildeHccarina}, one can see that the hypotheses of Theorem \ref{MD} are satisfied. Then, for $\varepsilon$ small enough with respect to $\Upsilon$, $m_c$ and $M$, there exist $v$, $\Gamma:\T^4 \times \R^2 \times J \to \R^{6}$, with $v$, $\Gamma \in \mathcal{S}_{1, 1}$ such that letting
\begin{equation}
\label{WAChatH}
\varphi^t:\T^4 \times \R^2 \to \T^4 \times \R^2 \times B^4 \times B^2, \quad \varphi^t(\theta, \xi) = (\theta, \xi, v^t(\theta, \xi))
\end{equation}
for all $t \in J$, the above defined family of embeddings $\varphi^t$ is a $C^1$-weakly asymptotic cylinder associated to $(X_{\hat H}, X_{\hat h}, \varphi_0)$, where $\Gamma$ is the disturbing term associated to $\varphi^t$ (see Definition \ref{weakasymcyl}). Moreover
\begin{equation}
\label{StimeVGamma}
|v|_{1, 1}<1, \quad |\Gamma|_{1, 1} <1.
\end{equation}

\subsection{Existence of weakly asymptotically quasiperiodic solutions}\label{P5C}

In the previous section, we proved the existence of a $C^1$ weakly asymptotic cylinder $\varphi^t$ associated to $(X_{\hat H}, X_{\hat h}, \varphi_0)$ (see~\eqref{WAChatH}).  Building upon this result, Proposition \ref{propasymsolI} ensures the existence of weakly asymptotically quasiperiodic solutions associated to $(X_{\hat H}, X_{\hat h}, \varphi_0)$. In the present section, for $v_c$ large enough, we show the existence of an open set of initial points giving rise to weakly asymptotically quasiperiodic solutions associated to $(X_H, X_{H_0},  \phi \circ \varphi_0)$, where $H$ and $H_0$ are the Hamiltonians defined by~\eqref{Hc}, whereas $\phi$ is the symplectic transformation in~\eqref{phiC}.  For this purpose, we introduce the following notation.  Let $v_1$, $\Gamma_1 : \T^4 \times \R^2 \times J \to \R^4$ and  $v_2$, $\Gamma_2 : \T^4 \times \R^2 \times J \to \R^2$ be the components of the functions $v$ and $\Gamma$ at the end of the previous section. That is $v = (v_1, v_2)$ and $\Gamma = (\Gamma_1, \Gamma_2)$. 
We recall that $B_t^2/2$ is defined by~\eqref{Bt2emezzo}, and letting $\psi_{1, H}^t$ be the flow at time $t$ with initial time $1$ of $H$, we have the following 
\begin{lemma}
\label{Superlemmafinalecomet}
We assume
\begin{equation}
\label{hypvfinaleComet}
v_c>{12 \over \varepsilon}.
\end{equation}
Then, for all $w \in \mathcal{W} = \phi \circ \varphi^1 \left(\T^4 \times \left(B_1^2/2\right)\right)$, $\psi_{1, H}^t(w)$ is a weakly asymptotically quasiperiodic solution associated to $(X_H, X_{H_0}, \phi \circ \varphi_0)$.
\end{lemma}
\begin{proof}
Let $\psi_{1, \hat H}^t$ be the flow at time $t$ with initial time $1$ of $\hat H$.  Thanks to Proposition \ref{propasymsolI}, we know that, for all $(\theta, \xi) \in \T^4 \times \R^2$, $\psi^t_{1, \hat H} \circ \varphi^1(\theta, \xi)$ is a weakly asymptotically quasiperiodic solution associated to $(X_{\hat H}, X_{\hat h}, \varphi_0)$. 

In the first part of this proof, we will prove that, 
\begin{equation}
\label{HatHPsiU1/2}
\psi_{1, \hat H}^t \circ \varphi^1(\theta, \xi) \in \mathcal{U}_{1 \over 2}
\end{equation}
for all $(\theta, \xi) \in \T^4 \times \left(B_1^2/2\right)$,  and $t \in J$, where $ \mathcal{U}_{1 \over 2}$ is defined by~\eqref{Uemezzo}.  To this end, for all $(\theta, \xi) \in \T^4 \times \left(B_1^2/2\right)$ and $t \in J$, we define the following families of functions
 \begin{equation*}
 \theta_1^t:\T^4 \times \R^2 \to \R^4, \hspace{3mm}   \xi_1^t:\T^4 \times \R^2 \to \R^2, \hspace{3mm}  r_1^t:\T^4 \times \R^2 \to \R^4, \hspace{3mm}  \eta_1^t:\T^4 \times \R^2 \to \R^2
 \end{equation*}
 in such a way that 
\begin{equation}
\label{notazioneflussostrana}
( \theta^t_1(\theta, \xi), \xi^t_1(\theta, \xi),  r^t_1(\theta, \xi), \eta^t_1(\theta, \xi)) = \psi_{1, \hat H}^t \circ \varphi^1(\theta, \xi)
\end{equation}
for all $(\theta, \xi) \in \T^4 \times \left(B_1^2/2\right)$ and $t \in J$. We observe that, by the definition of $\mathcal{U}_{1 \over 2}$ and the above notation,~\eqref{HatHPsiU1/2} is equivalent to show that
\begin{equation}
\label{retaxi}
(r^t_1(\theta, \xi), \eta^t_1(\theta, \xi))\in B^4 \times B^2, \hspace{7mm}  {| \xi^t_{1}(\theta, \xi)| \over |c(t)|} <{\varepsilon \over 6}
\end{equation}
for all $(\theta, \xi) \in \T^4 \times \left(B_1^2/2\right)$ and $t \in J$. We recall that $B^2 \subset \R^2$ and $B^4 \subset \R^4$ are balls centered at the origin with an unspecified radius that, conventionally, we will take greater than $1$. 

We recall the notation $\bar \omega = (\omega, 0) \in \R^6$ (see~\eqref{BarOmega}) and we denote by $\psi_{1, \bar \omega + \Gamma}^t$ the flow at time $t$ with initial time $1$ of $\bar \omega + \Gamma$.  In order to prove~\eqref{retaxi}, and hence~\eqref{HatHPsiU1/2},  we recall that if $\varphi^t$ is a $C^1$-weakly asymptotic cylinder associated to $(X_{\hat H}, X_{\tilde h}, \varphi_0)$, then 
\begin{equation}
\label{remindhyp1Cbiss}
\psi_{1, \hat H}^t \circ \varphi^1(\theta, \xi) = \varphi^t \circ \psi^t_{1, \bar \omega + \Gamma}(\theta, \xi)
\end{equation}
for all $(\theta, \xi)\in \T^4 \times \R^2$ and $t \in J$ (we refer to Proposition \ref{hyp1WACbiss}). Using the notation~\eqref{notazioneflussostrana}, the latter equation~\eqref{remindhyp1Cbiss} implies that 
\begin{equation}
\label{notazioneflussostrana+relazionetori}
( \theta^t_1(\theta, \xi), \xi^t_1(\theta, \xi),  r^t_1(\theta, \xi), \eta^t_1(\theta, \xi)) = \psi_{1, \hat H}^t \circ \varphi^1(\theta, \xi) =  \varphi^t \circ \psi^t_{1, \bar \omega + \Gamma}(\theta, \xi)
\end{equation}
for all $(\theta, \xi) \in \T^4 \times \left(B_1^2/2\right)$ and $t \in J$. 
However, we observe that $\varphi^t$ is the identity with respect to $(\theta, \xi)$ (see~\eqref{WAChatH}). Then, thanks to~\eqref{notazioneflussostrana+relazionetori} and the special form of $\varphi^t$, one can see that
\begin{equation*}
  r^t_1(\theta, \xi)= v^t_1( \theta^t_1(\theta, \xi), \xi^t_1(\theta, \xi)), \quad  \eta^t_1(\theta, \xi) = v^t_2(\theta^t_1(\theta, \xi), \xi^t_1(\theta, \xi))
\end{equation*}
for all $(\theta, \xi) \in \T^4 \times \left(B_1^2/2\right)$ and $t \in J$.  Hence, thanks to~\eqref{StimeVGamma}
\begin{equation*}
| r^t_1(\theta, \xi) | \le |v_1^t|_{C^1} < {1 \over t} \le 1, \quad |\eta^t_1(\theta, \xi) | \le |v_2^t|_{C^1} < {1 \over t} \le 1
\end{equation*}
for all $(\theta, \xi) \in \T^4 \times \left(B_1^2/2\right)$ and $t \in J$. This proves the first part of~\eqref{retaxi}.

Concerning the second part of~\eqref{retaxi}, by~\eqref{notazioneflussostrana+relazionetori} and the special form of $\varphi^t$, one can see that $ \xi^t_{1}(\theta, \xi)$ is the unique solution of the following system 
\begin{equation*}
\begin{cases}
\partial_t {\xi}^t_{1}(\theta, \xi) = \Gamma_2(\theta^t_1(\theta, \xi), \xi^t_1(\theta, \xi), t)\\
\xi^1_{1}(\theta, \xi) = \xi
\end{cases}
\end{equation*}
where we recall that $\Gamma_2$ is the second component of the disturbing term $\Gamma$ associated to $\varphi^t$, $\theta \in \T^4$ is fixed, and $\xi \in \left(B^2_1/2\right)$.  Using the fundamental theorem of calculus and~\eqref{StimeVGamma}
\begin{equation}
\label{Primastimaxi}
|\xi^t_{1}(\theta, \xi)| \le |\xi^1_{1}(\theta, \xi)| + \int_1^t |\Gamma^\tau|_{C^1} d\tau \le |\xi| + \ln t 
\end{equation}
for all $t \in J$. We claim that
\begin{equation}
\label{condizionefinaleprimaparte}
{|\xi^t_{1}(\theta, \xi)| \over |c(t)|} \le {|\xi| + \ln t \over  |c(1)| +{v_c \over 2} (t-1)} < {\varepsilon \over 6}
\end{equation}
for all $t \in J$, where the first inequality of the latter is a consequence of~\eqref{Primastimaxi} and~\eqref{ccondition}. 
The second inequality of~\eqref{condizionefinaleprimaparte} is true for $t=1$. However, if we suppose the existence of  $t_0 >1$ in such a way that 
\begin{equation*}
 {|\xi| +\ln t_0 \over  |c(1)| +{v_c \over2} (t_0-1)} \ge {\varepsilon \over 6},
\end{equation*}
then we can rewrite the latter in the following way
\begin{equation*}
|\xi| - {\varepsilon \over 6}|c(1)| \ge  {\varepsilon \over 12}v_c(t_0 -1) - \ln t_0.
\end{equation*}
We observe that this is a contradiction because $\xi \in  \left(B^2_1/2\right)$ implies $|\xi| - {\varepsilon \over 6}|c(1)| < 0$ and thanks to~\eqref{hypvfinaleComet} one can prove $ {\varepsilon \over 12}v_c(t -1) - \ln t \ge 0$ for all $t \in J$.  Indeed, letting
\begin{equation*}
f:[1, +\infty) \to \R, \quad f(t) =  {\varepsilon \over 12}v_c(t -1) - \ln t
\end{equation*}
one can see that $f(1)=0$, and using~\eqref{hypvfinaleComet}, one can prove that $f'(t) >0$ for all $t >1$, where $f'$ stands for the derivative of $f$.  Hence,~\eqref{condizionefinaleprimaparte} is verified for all $t \in J$, and this concludes the proof of~\eqref{retaxi} and thus~\eqref{HatHPsiU1/2}.

In this second part of this proof, we observe that~\eqref{HatHPsiU1/2} and Lemma \ref{lemmaHex} ensure that for all $(\theta, \xi) \in \T^4 \times \left(B_1^2/2\right)$, 
\begin{equation}
\label{UltimaEqFlussoMD}
\psi_{1, \hat H}^t \circ \varphi^1(\theta, \xi) = \psi_{1, H \circ \tilde \phi}^t \circ \varphi^1(\theta, \xi)
\end{equation}
for all $t \in J$. 
 
Now, for all $w \in \mathcal{W} = \phi \circ \varphi^1 \left(\T^4 \times \left(B_1^2 /2\right)\right)$, we want to prove that $\psi_{1, H}^t(w)$ is a weakly asymptotically quasiperiodic solution associated to $(X_H, X_{H_0}, \phi \circ \varphi_0)$. For this purpose, we observe that, for all $w \in \mathcal{W} = \phi \circ \varphi^1 \left(\T^4 \times \left(B_1^2 /2\right)\right)$, there exists $(\theta, \xi) \in \T^4 \times \left(B_1^2 /2\right)$ such that $w = \phi \circ \varphi^1 (\theta, \xi)$. We know that $\phi$ is symplectic, then by~\eqref{remindhyp1Cbiss} and~\eqref{UltimaEqFlussoMD}, we can rewrite $\psi_{1, H}^t(w)$ in the following form
\begin{eqnarray*}
\psi_{1, H}^t(w) &=& \psi_{1, H}^t \circ \phi \circ \varphi^1 (\theta, \xi) = \phi \circ \psi_{1, H \circ \tilde \phi}^t \circ \varphi^1 (\theta, \xi) = \phi \circ \psi_{1, \hat H}^t \circ \varphi^1 (\theta, \xi)\\
&=& \phi \circ \varphi^t \circ \psi^t_{1, \bar \omega + \Gamma} (\theta, \xi).
\end{eqnarray*}
Moreover, for all $t \in J$ 
\begin{eqnarray*}
\left|\psi_{1, H}^t(w)  - \phi \circ \varphi_0\circ \psi^t_{1, \bar \omega + \Gamma} (\theta, \xi)\right| &\le& \left|\phi \circ \varphi^t \circ \psi^t_{1, \bar \omega + \Gamma} (\theta, \xi)  - \phi \circ \varphi_0\circ \psi^t_{1, \bar \omega + \Gamma} (\theta, \xi)\right|\\
&\le& \left|\phi \circ \varphi^t  - \phi \circ \varphi_0\right|_{C^0}\\
&\le& C\left|D \phi\right|_{C^1}\left|\varphi^t  -  \varphi_0\right|_{C^0}\\
&\le& C\Upsilon^1_\phi\left|\varphi^t  -  \varphi_0\right|_{C^0}
\end{eqnarray*}
where the constant $\Upsilon^1_\phi$ is defined by~\eqref{UpsilonPhi}. Taking the limit for $t \to +\infty$, thanks to~\eqref{hyp2WAC}, we conclude the proof of this lemma.
\end{proof}

This concludes the proof of Theorem \ref{Thmcomet}.

\appendix
\section{Hölder classes of functions}
\label{A}

This part is dedicated to recalling the definition of Hölder classes of functions and some well-known properties. To this end, let $E$ be an open subset of $\R^n$ and $k \ge 0$ a positive integer. We denote by $C^k(E)$ the spaces of functions $f: E \to \R$ with continuous partial derivatives $\partial^\alpha f \in C^0(E)$ for all $\alpha \in \N^n$ with $|\alpha|=\alpha_1+...+\alpha_n \le k$ and verifying
\begin{equation}
\label{Cknorm}
|f|_{C^k(E)} = \sup_{|\alpha|\le k}|\partial^\alpha f|_{C^0}=  \sup_{|\alpha|\le k}\sup_{x \in E}|\partial^\alpha f(x)| <\infty.
\end{equation}
Given $\sigma=k+\mu$, with $k \in \Z$, $k \ge 0$ and $0 < \mu <1$, we define the Hölder spaces $C^\sigma(E)$ as the spaces of functions $f\in C^k(E)$ verifying
\begin{equation}
\label{Holdernorm}
|f|_{C^\sigma(E)} = \sup_{|\alpha|\le k}|\partial^\alpha f|_{C^0(E)} + \sup_{x,y\in E\\, x\ne y, |\alpha| = k}{|\partial^\alpha f(x) - \partial^\alpha f(y)| \over |x-y|^\mu}<\infty.
\end{equation}
We will use the same notations for vector-valued functions or matrices. More specifically,  in the case of functions $f=(f_1,...,f_n)$ with values in $\R^n$, we set $|f|_{C^\sigma(E)} = \max_{1 \le i \le n} |f_i|_{C^\sigma(E)}$. Moreover, in agreement with the convention made above, if $M =\{m_{ij}\}_{1 \le i,j \le n}$ is a $n \times n$ matrix, we set $|M|_{C^\sigma(E)}= \max_{1 \le i,j \le n} |m_{ij}|_{C^\sigma(E)}$.  In what follows, and generally, if there are no doubts about the domain where the norms are considered, we prefer to omit the domain $E$ in the notation $|\cdot|_{C^\sigma}$ of the norm.


Now, we present a series of properties widely used in this paper.  First, we recall that $C(\cdot)$ stands for constants depending on  $n$ and other parameters in brackets. The following proposition provides a convexity property of the above-mentioned norms. 
\begin{proposition}
\label{convexity}
For all $f\in C^{\sigma_1}(E)$, then 
\begin{equation*}
|f|^{\sigma_1-\sigma_0}_{C^\sigma}\le C(\sigma_1) |f|^{\sigma_1-\sigma}_{C^{\sigma_0}}|f|^{\sigma-\sigma_0}_{C^{\sigma_1}} \hspace{3mm} \mbox{for all $0 \le \sigma_0 \le \sigma \le \sigma_1$}.
\end{equation*}
\end{proposition}
\begin{proof}
We refer to~\cite{Hor76} for the proof. 
\end{proof}

The next proposition deals with the composition and product of Hölder functions 
\begin{proposition}
\label{Holder}
We consider $f$, $g \in C^\sigma(E)$ and $\sigma \ge 0$.
\begin{enumerate}
\item For all $\beta \in \N^{n}$ and $s\ge 0$, if $|\beta| + s \le \sigma$ then  $\left|{\partial^{|\beta|} \over \partial{x_1}^{\beta_1}... \partial{x_n}^{\beta_n}} f \right|_{C^s} \le C|f|_{C^\sigma}$.\\
\item  $|fg|_{C^\sigma} \le C(\sigma)\left(|f|_{C^0}|g|_{C^\sigma} + |f|_{C^\sigma}|g|_{C^0}\right)$. 
\end{enumerate}
Now, we consider composite functions.  Let $E_1$ be an open subset of $\R^n$ and $z:E_1\to E$ a function taking values in the domain of $f$.  In what follows $\partial z$ stands for the partial derivatives of $z$.
\begin{enumerate}
\item[3.] If $\sigma < 1$, $f \in C^1(E)$, $z \in C^\sigma (E_1)$ then $f\circ z \in C^\sigma(E_1)$ and
$$
|f \circ z|_{C^\sigma} \le C(|f|_{C^1}|z|_{C^\sigma}+ |f|_{C^0}).
$$
\end{enumerate} 
\begin{enumerate}
\item[4.] If $\sigma < 1$, $f \in C^\sigma(E)$, $z \in C^1 (E_1)$ then $f\circ z \in C^\sigma(E_1)$ 
$$
|f \circ z|_{C^\sigma} \le C(|f|_{C^\sigma}|\partial z|^\sigma_{C^0}+ |f|_{C^0}).
$$
\end{enumerate}
\begin{enumerate}
\item[5.] If $\sigma \ge 1$ and $f \in C^\sigma (E)$, $z \in C^\sigma (E_1)$ then $f\circ z \in C^\sigma(E_1)$ 
$$
|f \circ z|_{C^\sigma} \le C(\sigma) \left(|f|_{C^\sigma}|\partial z|^\sigma_{C^0} + |f|_{C^1}|\partial z|_{C^{\sigma-1}}+ |f|_{C^0}\right).
$$
\end{enumerate}
\end{proposition}
\begin{proof}
We refer to~\cite{Hor76} for the proof of \textit{1}, \textit{2}, and \textit{3}.  We begin by proving \textit{4}.  For this purpose, we claim that 
\begin{equation}
\label{RobaPer4}
{|f\circ z(x) - f\circ z(y)| \over |x-y|^\sigma} \le |f|_{C^\sigma}|\partial z|^\sigma_{C^0}
\end{equation}
for all $x$, $y \in E_1$ with $x \ne y$. It is obvious that the latter is satisfied if  $z(x) = z(y)$. Now, we assume $z(x) \ne z(y)$,  and we observe that
\begin{eqnarray}
{|f\circ z(x) - f\circ z(y)| \over |x-y|^\sigma} &=& {|f\circ z(x) - f\circ z(y)| \over |z(x)-z(y)|^\sigma}\left( {|z(x) -  z(y)| \over |x-y|} \right)^\sigma \nonumber\\
\label{AppProp2_4_1}
&\le&|f|_{C^\sigma}\left( \int_0^1 |\partial z(y+ \tau(x-y))|d\tau \right)^\sigma\\
&\le& |f|_{C^\sigma}|\partial z|^\sigma_{C^0} \nonumber
\end{eqnarray}
where the inequality~\eqref{AppProp2_4_1} is a consequence of the definition~\eqref{Holdernorm}, and Taylor formula. This proves the claim. The proof of \textit{4}, now is a consequence of~\eqref{RobaPer4} and the definition~\eqref{Holdernorm}.

It remains to prove \textit{5.} By~\eqref{Holdernorm},   
\begin{eqnarray}
\label{proofHoldercomp}
|f \circ z|_{C^\sigma} &\le& |f|_{C^0} + |(\partial f\circ z)^T \partial z |_{C^{\sigma -1}},
\end{eqnarray}
where $\partial f$ stands for the partial derivatives of $f$,  and $T$ for the transpose. Thanks to property \textit{2}
\begin{eqnarray*}
|f \circ z|_{C^\sigma} &\le& |f|_{C^0} + |(\partial f\circ z)^T \partial z |_{C^{\sigma -1}}\\
&\le& |f|_{C^0} +C(\sigma)|\partial f\circ z|_{C^{\sigma -1}}|\partial z |_{C^0} + C(\sigma)|\partial f\circ z|_{C^0}|\partial z |_{C^{\sigma -1}}.
\end{eqnarray*} 

We observe that $|\partial f\circ z|_{C^0}|\partial z |_{C^{\sigma -1}} \le |f|_{C^1}|\partial z |_{C^{\sigma -1}}$, it remains to estimate $|\partial f\circ z|_{C^{\sigma -1}}|\partial z |_{C^0}$. If $\sigma \le 2$, $|\partial f \circ z|_{C^{\sigma-1}} \le | f|_{C^\sigma}|\partial z|^{\sigma-1}_{C^0}+ |f|_{C^1}$ thanks to \textit{4}. Then 
\begin{eqnarray*}
|\partial f \circ z|_{C^{\sigma-1}} |\partial z|_{C^0}&\le& C(\sigma) \left(| f|_{C^\sigma}|\partial z|^\sigma_{C^0}+ |f|_{C^1}|\partial z|_{C^0}\right)\\
&\le& C(\sigma) \left(| f|_{C^\sigma}|\partial z|^\sigma_{C^0}+ |f|_{C^1}|\partial z|_{C^{\sigma -1}}\right),
\end{eqnarray*}
whence the property holds in this case. If $\sigma >2$, assuming that \textit{5} is already proven for $\sigma -1$, we find
\begin{eqnarray*}
|\partial f \circ z|_{C^{\sigma-1}} |\partial z|_{C^0}&\le& C(\sigma)\left(|\partial f|_{C^{\sigma-1}}|\partial z|^\sigma_{C^0}+ |f|_{C^2}|\partial z|_{C^{\sigma-2}}|\partial z|_{C^0} + |f|_{C^1}|\partial z|_{C^0}\right)\\
&\le&C(\sigma)\left(|f|_{C^\sigma}|\partial z|^\sigma_{C^0}+ |f|_{C^2}|\partial z|_{C^{\sigma-2}}|\partial z|_{C^0} + |f|_{C^1}|\partial z|_{C^{\sigma -1}}\right).
\end{eqnarray*}
It remains to find a good estimate for the central term of the last line of the latter. By Proposition \ref{convexity}
\begin{eqnarray*}
|f|_{C^2}|\partial z|_{C^{\sigma-2}}|\partial z|_{C^0} &\le& C(\sigma) \left(|f|^{\sigma-2 \over \sigma-1}_{C^1}|f|_{C^\sigma}^{1\over \sigma -1}\right) \left(|\partial z|^{1\over \sigma -1}_{C^0}|\partial z|_{C^{\sigma-1}}^{\sigma-2\over \sigma -1}\right)|\partial z|_{C^0} \\
&\le&C(\sigma)  \left(|f|_{C^1}|\partial z|_{C^{\sigma-1}}\right)^{\sigma-2 \over \sigma-1}  \left(|f|_{C^\sigma}|\partial z|_{C^0}^\sigma\right)^{1 \over \sigma-1},
\end{eqnarray*}
since $a^\lambda b^{1-\lambda} \le C(a+b)$ for $0 < \lambda < 1$, we have the claim. 
\end{proof}


\section*{Acknowledgement}

\textit{These results are part of my PhD thesis, which I prepared at Université Paris-Dauphine. I want to thank my thesis advisors, Abed Bounemoura and Jacques Féjoz. Without their advice and support, this work would not exist. Thank you to  Inmaculada Baldomá Barraca for the interesting discussions we had.}

\textit{This project has received funding from the European Union’s Horizon 2020 research and innovation programme under the Marie Skłodowska-Curie grant agreement No 754362.} \includegraphics[scale=0.01]{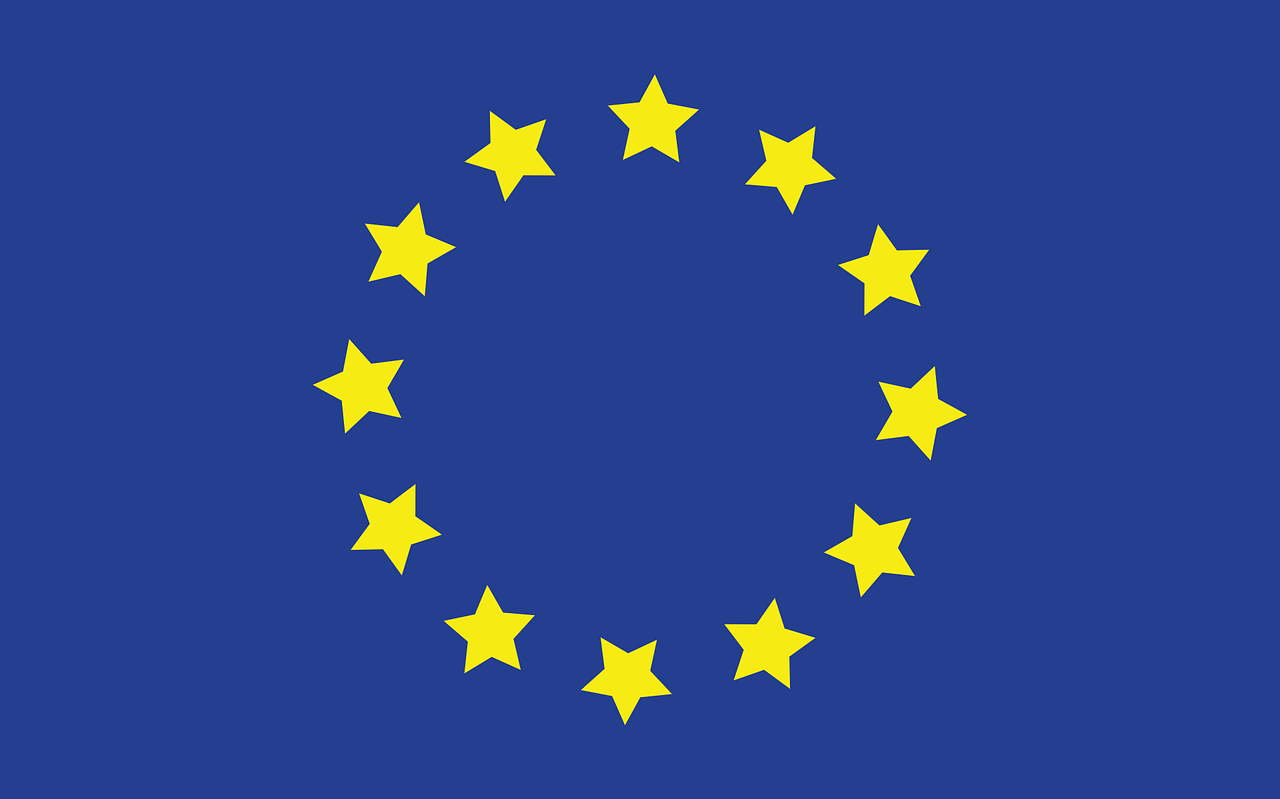}

\bibliographystyle{plain}
\bibliography{ref}

\begin{thebibliography}{10}

\bibitem{Amann}
Herbert Amann.
\newblock {\em Ordinary differential equations}, volume~13 of {\em De Gruyter
  Studies in Mathematics}.
\newblock Walter de Gruyter \& Co., Berlin, 1990.
\newblock An introduction to nonlinear analysis, Translated from the German by
  Gerhard Metzen.

\bibitem{Arn63a}
V.~I. Arnold.
\newblock Proof of a theorem of {A}. {N}. {K}olmogorov on the preservation of
  conditionally periodic motions under a small perturbation of the
  {H}amiltonian.
\newblock {\em Uspehi Mat. Nauk}, 18(5 (113)):13--40, 1963.

\bibitem{Arn63b}
V.~I. Arnold.
\newblock Small denominators and problems of stability of motion in classical
  and celestial mechanics.
\newblock {\em Uspehi Mat. Nauk}, 18(6 (114)):91--192, 1963.

\bibitem{BFMn+1}
Inmaculada Baldom\'{a}, Ernest Fontich, and Pau Mart\'{\i}n.
\newblock Whiskered parabolic tori in the planar {$(n+1)$}-body problem.
\newblock {\em Comm. Math. Phys.}, 374(1):63--110, 2020.

\bibitem{BGGS84}
G.~Benettin, L.~Galgani, A.~Giorgilli, and J.-M. Strelcyn.
\newblock A proof of {K}olmogorov's theorem on invariant tori using canonical
  transformations defined by the {L}ie method.
\newblock {\em Nuovo Cimento B (11)}, 79(2):201--223, 1984.

\bibitem{BB15}
Massimiliano Berti and Philippe Bolle.
\newblock A {N}ash-{M}oser approach to {KAM} theory.
\newblock In {\em Hamiltonian partial differential equations and applications},
  volume~75 of {\em Fields Inst. Commun.}, pages 255--284. Fields Inst. Res.
  Math. Sci., Toronto, ON, 2015.

\bibitem{BdlL11}
Daniel Blazevski and Rafael de~la Llave.
\newblock Time-dependent scattering theory for {ODE}s and applications to
  reaction dynamics.
\newblock {\em J. Phys. A}, 44(19):195101, 26, 2011.

\bibitem{BDT}
Alberto Boscaggin, Walter Dambrosio, and Susanna Terracini.
\newblock Scattering parabolic solutions for the spatial {$N$}-centre problem.
\newblock {\em Arch. Ration. Mech. Anal.}, 223(3):1269--1306, 2017.

\bibitem{bost1986tores}
Jean-Beno{\^\i}t Bost.
\newblock Tores invariants des systemes dynamiques hamiltoniens.
\newblock {\em S{\'e}minaire Bourbaki}, 1984:85, 1986.

\bibitem{CdlL15}
Marta Canadell and Rafael de~la Llave.
\newblock K{AM} tori and whiskered invariant tori for non-autonomous systems.
\newblock {\em Phys. D}, 310:104--113, 2015.

\bibitem{Chazy22}
Jean Chazy.
\newblock Sur l'allure du mouvement dans le probl\`eme des trois corps quand le
  temps cro\^{i}t ind\'{e}finiment.
\newblock {\em Ann. Sci. \'{E}cole Norm. Sup. (3)}, 39:29--130, 1922.

\bibitem{Che89}
A.~Chenciner.
\newblock Intégration du problème de kepler par la méthode de
  hamilton-jacobi: coordonnées action-angle de delaunay.
\newblock {\em Notes scientifiques et techniques du Bureau des Longitudes},
  S026, 1989.

\bibitem{Chi03}
Luigi Chierchia.
\newblock K{AM} lectures.
\newblock In {\em Dynamical systems. {P}art {I}}, Pubbl. Cent. Ric. Mat. Ennio
  Giorgi, pages 1--55. Scuola Norm. Sup., Pisa, 2003.

\bibitem{PC11}
Luigi Chierchia and Gabriella Pinzari.
\newblock The planetary {$N$}-body problem: symplectic foliation, reductions
  and invariant tori.
\newblock {\em Invent. Math.}, 186(1):1--77, 2011.

\bibitem{FW14}
Alessandro Fortunati and Stephen Wiggins.
\newblock Persistence of {D}iophantine flows for quadratic nearly integrable
  {H}amiltonians under slowly decaying aperiodic time dependence.
\newblock {\em Regul. Chaotic Dyn.}, 19(5):586--600, 2014.

\bibitem{Fe02}
Jacques Féjoz.
\newblock Quasiperiodic motions in the planar three-body problem.
\newblock {\em J. Differential Equations}, 183(2):303--341, 2002.

\bibitem{Fe04}
Jacques Féjoz.
\newblock D\'{e}monstration du `th\'{e}or\`eme d'{A}rnold' sur la stabilit\'{e}
  du syst\`eme plan\'{e}taire (d'apr\`es {H}erman).
\newblock {\em Ergodic Theory Dynam. Systems}, 24(5):1521--1582, 2004.

\bibitem{Fe13}
Jacques Féjoz.
\newblock On action-angle coordinates and the {P}oincar\'{e} coordinates.
\newblock {\em Regul. Chaotic Dyn.}, 18(6):703--718, 2013.

\bibitem{Fe16}
Jacques Féjoz.
\newblock Introduction to {KAM} theory with a view to celestial mechanics.
\newblock In {\em Variational methods}, volume~18 of {\em Radon Ser. Comput.
  Appl. Math.}, pages 387--433. De Gruyter, Berlin, 2017.

\bibitem{Hor76}
Lars Hörmander.
\newblock The boundary problems of physical geodesy.
\newblock {\em Arch. Rational Mech. Anal.}, 62(1):1--52, 1976.

\bibitem{kawai2007transition}
Shinnosuke Kawai, Andr{\'e}~D Bandrauk, Charles Jaff{\'e}, Thomas Bartsch,
  Jesus Palacian, and T~Uzer.
\newblock Transition state theory for laser-driven reactions.
\newblock {\em The Journal of chemical physics}, 126(16):164306, 2007.

\bibitem{Kol54}
A.~N. Kolmogorov.
\newblock On conservation of conditionally periodic motions for a small change
  in {H}amilton's function.
\newblock {\em Dokl. Akad. Nauk SSSR (N.S.)}, 98:527--530, 1954.

\bibitem{Las89}
J.~Laskar.
\newblock Les variables de poincaré et le développement de la fonction
  pertubatrice.
\newblock {\em Notes scientifiques et techniques du Bureau des Longitudes},
  S026, 1989.

\bibitem{MV}
E.~Maderna and A.~Venturelli.
\newblock Globally minimizing parabolic motions in the {N}ewtonian {$N$}-body
  problem.
\newblock {\em Arch. Ration. Mech. Anal.}, 194(1):283--313, 2009.

\bibitem{M18}
Jessica~Elisa Massetti.
\newblock A normal form \`a la {M}oser for diffeomorphisms and a generalization
  of {R}\"{u}ssmann's translated curve theorem to higher dimensions.
\newblock {\em Anal. PDE}, 11(1):149--170, 2018.

\bibitem{M19}
Jessica~Elisa Massetti.
\newblock Normal forms for perturbations of systems possessing a {D}iophantine
  invariant torus.
\newblock {\em Ergodic Theory Dynam. Systems}, 39(8):2176--2222, 2019.

\bibitem{MDS}
Dusa McDuff and Dietmar Salamon.
\newblock {\em Introduction to symplectic topology}.
\newblock Oxford Graduate Texts in Mathematics. Oxford University Press,
  Oxford, third edition, 2017.

\bibitem{moser1962invariant}
J.~Moser.
\newblock On invariant curves of area-preserving mappings of an annulus.
\newblock {\em Nachr. Akad. Wiss. G\"{o}ttingen Math.-Phys. Kl. II},
  1962:1--20, 1962.

\bibitem{Po82}
J\"{u}rgen P\"{o}schel.
\newblock Integrability of {H}amiltonian systems on {C}antor sets.
\newblock {\em Comm. Pure Appl. Math.}, 35(5):653--696, 1982.

\bibitem{Posc01}
Jurgen Pöschel.
\newblock A lecture on the classical {KAM} theorem.
\newblock In {\em Smooth ergodic theory and its applications ({S}eattle, {WA},
  1999)}, volume~69 of {\em Proc. Sympos. Pure Math.}, pages 707--732. Amer.
  Math. Soc., Providence, RI, 2001.

\bibitem{Sca22T}
Donato Scarcella.
\newblock Asymptotically quasiperiodic solutions for time-dependent
  hamiltonians.
\newblock {\em Dynamical Systems [math.DS]. Université Paris sciences et
  lettres, 2022. English. ⟨NNT : 2022UPSLD047⟩. ⟨tel-03963784v2⟩},
  2022.

\bibitem{Sca22c}
Donato Scarcella.
\newblock Asymptotic motions converging to arbitrary dynamics for
  time-dependent hamiltonians.
\newblock {\em Nonlinear Analysis}, 243:113528, 2024.

\bibitem{Sca22}
Donato Scarcella.
\newblock Asymptotically quasiperiodic solutions for time-dependent
  hamiltonians.
\newblock {\em Nonlinearity}, 37(6):065005, 2024.

\bibitem{Sca22b}
Donato Scarcella.
\newblock Biaymptotically quasiperiodic solutions for time-dependent
  hamiltonians.
\newblock {\em Regular and Chaotic Dynamics}, 29:620--653, 2024.

\bibitem{thieme1995asymptotically}
H~Thieme and C~Castillo-Chavez.
\newblock {\em Asymptotically autonomous epidemic models. Mathematical
  Population Dynamics: Analysis of Heterogeneity, vol. I}.
\newblock Wuerz publishing, 1995.

\bibitem{Zeh76}
E.~Zehnder.
\newblock Generalized implicit function theorems with applications to some
  small divisor problems. {I}.
\newblock {\em Comm. Pure Appl. Math.}, 28:91--140, 1975.

\bibitem{Zeh75}
E.~Zehnder.
\newblock Generalized implicit function theorems with applications to some
  small divisor problems. {II}.
\newblock {\em Comm. Pure Appl. Math.}, 29(1):49--111, 1976.

\end{thebibliography}

\end{document}